\documentclass[preprint,reqno,addressatend,12pt,onehalfspacing]{imsart}
\RequirePackage[OT1]{fontenc}
\RequirePackage{amsthm,amsmath}
\RequirePackage[colorlinks,citecolor=blue,urlcolor=blue]{hyperref}
\RequirePackage{natbib}
\usepackage{a4wide,xcolor,eucal,enumerate,mathrsfs,graphics,caption,subfig,booktabs,verbatim}
\usepackage{pdfsync}
\usepackage{amsmath,amssymb,epsfig,amsthm}
\usepackage[latin1]{inputenc}
\usepackage{dsfont}
\startlocaldefs

\theoremstyle{definition}
\newtheorem{de}{Definition}[section]

\theoremstyle{plain}
\newtheorem{theo}[de]{Theorem}
\newtheorem{lemma}[de]{Lemma}

\theoremstyle{remark}
\newtheorem{re}[de]{Remark}

\newcommand{\R}{\mathbb{R}}
\newcommand{\Z}{\mathbb{Z}}
\newcommand{\p}{\mathbb{P}}

\newcommand{\E}{\mathbb{E}}

\newcommand{\N}{\mathbb{N}}

\newcommand{\F}{\mathcal{F}}
\newcommand{\B}{\mathcal{B}}
\newcommand{\G}{\mathcal{G}}

\newcommand{\A}{\mathcal{A}}
\newcommand{\1}{\mathds{1}}
\newcommand{\M}{\mathcal{M}}
\newcommand{\m}{\mathbb{M}}
\newcommand{\w}{\mathbb{B}}

\def\argmin{\mathop{\mathrm{argmin}}}

\endlocaldefs

\begin{document}

\begin{frontmatter}
\title{Smoothed isotonic estimators of a monotone baseline hazard in the Cox model}
\runtitle{Smoothed isotonic hazard estimators}

\begin{aug}
\author{\fnms{Hendrik P.} \snm{Lopuha\"a}\ead[label=e1]{h.p.lopuhaa@tudelft.nl}},
\author{\fnms{Eni} \snm{Musta}\ead[label=e2]{e.musta@tudelft.nl}}

\runauthor{H.P.~Lopuha\"a and E.~Musta}

\affiliation{Delft Institute of Applied Mathematics, Delft University of Technology}  

\address{Delft Institute of Applied Mathematics\\
Delft University of Technology\\
van Mourik Broekmanweg 6, 2628 XE Delft, Netherlands\\
\printead{e1}\\
\printead{e2}
}
\end{aug}

\begin{abstract}
We consider the smoothed maximum likelihood estimator and the smoothed Grenander-type estimator for a monotone baseline hazard rate
$\lambda_0$ in the Cox model.
We analyze their asymptotic behavior and show that they are asymptotically normal at rate $n^{m/(2m+1)}$,
when~$\lambda_0$ is $m\geq 2$ times continuously differentiable,
and that both estimators are asymptotically equivalent.
Finally, we present numerical results on pointwise confidence intervals that 
illustrate the comparable behavior of the two methods.
\end{abstract}


\end{frontmatter}

\section{Introduction}
The semi-parametric Cox regression model is a very popular method  in survival analysis that allows incorporation of covariates when studying lifetime distributions in the presence of right censored data.
The ease of interpretation, resulting from the formulation in terms of the hazard rate as well as the proportional effect of the covariates, and the fact that the regression coefficients (parametric component) can be estimated while leaving the baseline distribution unspecified, favor the wide use of this framework, especially in medical applications.
On the other hand, since its first introduction (see~\cite{Cox72}), much effort has been spent on giving a firm mathematical basis to this approach.
Initially, the attention was on the derivation of large sample properties of the maximum partial likelihood estimator of the regression coefficients and of the Breslow estimator for the cumulative baseline hazard (e.g., see~\cite{Efron72}, \cite{Cox75}, \cite{Tsiatis81}).
Although the most attractive property of this approach is that it does not assume any fixed shape on the hazard curve, there are several cases
where order restrictions better match the practical expectations (e.g., see~\cite{Geloven13} for an example of a decreasing hazard in a large clinical trial for patients with acute coronary syndrome).
Estimation of the baseline hazard function under monotonicity constraints has been studied in~\cite{CC94} and~\cite{LopuhaaNane2013}.

Traditional isotonic estimators, such as maximum likelihood estimators and Grenander-type estimators
are step functions that exhibit a non normal limit distribution at rate $n^{1/3}$.
On the other hand, a long stream of research has shown that, if one is willing to assume more regularity on the function of interest,
smooth estimators might be preferred to piecewise constant ones because they can be used to achieve a faster rate of convergence to a Gaussian distributional law and to estimate derivatives.
Isotonized smooth estimators, obtained either by a least squares projection, maximum likelihood, or penalization,
are studied in~\cite{mukerjee1988}, \cite{ramsay1998}, \cite{eggermont-lariccia2000}, \cite{vdvaart-vdlaan2003},
and in~\cite{mammen1991},
who also compares isotonized kernel estimators with smoothed isotonic estimators in the regression setting.
Smoothed maximum likelihood estimators for monotone functions have also been investigated
by~\cite{DGL13}, to bootstrap from a smooth decreasing density estimate,
by~\cite{GJW10} for the current status model, together with a maximum smoothed likelihood estimator, and
by~\cite{GJ13} for estimating a monotone hazard rate,
together with a penalized least squares estimator.
Other references for combining shape constraints and smoothness can be found in Chapter~8 in~\cite{GJ14}.
Distribution theory was first studied by~\cite{mukerjee1988}, who established asymptotic normality for a kernel smoothed least squares regression estimator,
but this result is limited to a rectangular kernel and the rate of convergence is slower than the usual rate for kernel estimators.
In~\cite{vdvaart-vdlaan2003} it is shown that the isotonized kernel density estimator has the same limit normal distribution
at the usual rate $n^{m/(2m+1)}$ as the ordinary kernel density estimator, when the density is $m$ times continuously differentiable.
Similar results were obtained by~\cite{GJW10} for the smoothed maximum likelihood estimator and
the maximum smoothed likelihood estimator, and by~\cite{GJ13} for a smoothed Grenander-type estimator.

Smooth estimation under monotonicity constraints for the baseline hazard in the Cox model was introduced in~\cite{Nane}. 
By combining an isotonization step with a smoothing step and alternating the order of smoothing and isotonization, four different estimators can be constructed. Two of them are kernel smoothed versions of the maximum likelihood estimator and the Grenander-type estimator from~\cite{LopuhaaNane2013}. The third estimator is a maximum smoothed likelihood estimator obtained by first smoothing the loglikelihood of the Cox model and then finding the maximizer of the smoothed likelihood among all decreasing baseline hazards. The forth one is a Grenander-type estimator based on the smooth Breslow estimator for the cumulative hazard. 
{Three of these} estimators were shown to be consistent in~\cite{Nane}. Moreover, the last two methods have been studied in~\cite{LopuhaaMusta2016} and were shown to be asymptotically normal at the usual rate $n^{m/(2m+1)}$, where $m$ denotes the level of smoothness of the baseline hazard. The main interest of the present paper is to investigate the asymptotic behavior of the first two methods, the smoothed maximum likelihood estimator and a smoothed Grenander-type estimator.

This is particularly challenging for the Cox model, because the existing approaches to these type of problems
for smoothed isotonic estimators do not apply to the Cox model.
The situation is different from isotonized smooth estimators, such as the maximum smoothed likelihood estimator and a Grenander-type estimator based on the smoothed Breslow estimator, which are studied in~\cite{LopuhaaMusta2016}. 
In the latter paper, the main idea is that the isotonic smooth estimator can be represented as a least squares projection of a naive estimator (smooth but not monotone), 
which is simpler to analyze and asymptotically equivalent to  the isotonic smooth estimator.

{The smoothed Grenander-type estimator in the ordinary right censoring model without covariates was investigated by~\cite{LopuhaaMustaSN17}.
Following the approach in~\cite{GJ13}, asymptotic normality was established by using a Kiefer-Wolfowitz type of result, recently derived in~\cite{DL14}.
Unfortunately, the lack of a Kiefer-Wolfowitz type of result for the Breslow estimator
provides a strong limitation towards extending the previous approach to the more general setting of the Cox model.
Recently, \cite{GJ14} developed a different method for finding the limit distribution of smoothed isotonic estimators,
which is mainly based on uniform $L_2$-bounds on the distance between the non-smoothed isotonic estimator and the true function,
and also uses that the maximal distance between succeeding points of jump of the isotonic estimator is of the order~$O_p(n^{-1/3}\log n)$.
A sketch of proof in the right censoring model is given in Section~11.6 of~\cite{GJ14}.
However, these two key ingredients heavily depend on having exponential bounds for tail probabilities of the so-called inverse process,
or rely on a strong embedding for the relevant sum process.
Exponential bounds for tail probabilities of the inverse process are difficult to obtain in the Cox model and a strong embedding for the Breslow estimator is not available.
Nevertheless, inspired by the approach in~\cite{GJ14}, we obtain polynomial bounds, which will suffice for obtaining uniform $L_2$-bounds,
and we avoid using the maximal distance between succeeding points of jump of the non-smoothed isotonic estimator,
by establishing a sufficiently small bound on the expected supremum distance between the non-smoothed isotonic estimator and the true
baseline hazard.}

This leads to asymptotic normality at rate $n^{m/(2m+1)}$ of
the smoothed maximum likelihood estimator and the smoothed Grenander-type estimator,
which are also shown to be asymptotically equivalent.
By means of a small simulation we investigate the finite sample behavior in terms of
asymptotic confidence intervals corresponding to the limit normal distributions, as well as 
bootstrap confidence intervals based on a smooth bootstrap proposed by~\cite{Burr94} and~\cite{XSY2014}.
As expected, no estimator performs strictly better than the other.

The paper is organized as follows.
In Section~\ref{sec:model} we specify the Cox regression model and provide some background information that will be used  in the sequel.
The kernel smoothed versions of the Grenander-type estimator and of the maximum likelihood estimator of a non-decreasing
baseline hazard function are considered in Section~\ref{sec:SG}.
We only consider the case of a non-decreasing baseline hazard.
The same results can be obtained similarly for a non-increasing hazard.
{The results of a small simulation study are reported in Section~\ref{sec:conf-int}
	and we conclude with a brief discussion in Section~\ref{sec:discussion}.}
In order to keep the exposition clear and simple, most of the proofs are delayed until Section~\ref{sec:proofs},
and remaining technicalities have been put in the Supplemental Material~\cite{LopuhaaMustaSM2016}.

\section{The Cox regression model}\label{sec:model}
Let $X_1,\dots,X_n$ be an i.i.d.~sample representing the survival times of $n$ individuals, which can be observed only on time intervals $[0,C_i]$ for some i.i.d. censoring times $C_1,\dots,C_n$.
One observes i.i.d. triplets $(T_1,\Delta_1,Z_1),\dots,(T_n,\Delta_n,Z_n)$, where $T_i=\min(X_i,C_i)$ denotes the follow up time,  $\Delta_i=\1_{\{X_i\leq C_i\}}$ is the censoring indicator and $Z_i\in\R^p$ is a time independent covariate vector. Given the covariate vector $Z,$ the event time $X$ and the censoring time $C$ are assumed to be independent. Furthermore, conditionally on $Z=z,$ the event time is assumed to be a nonnegative r.v. with an absolutely continuous distribution function $F(x|z)$ and density $f(x|z).$ Similarly the censoring time is assumed to be a nonnegative r.v. with an absolutely continuous distribution function $G(x|z)$ and density $g(x|z).$ The censoring mechanism is assumed to be non-informative, i.e. $F$ and $G$ share no parameters.
Within the Cox model, the conditional hazard rate $\lambda(x|z)$ for a subject with covariate vector $z\in\R^p$, is related to the corresponding covariate by
\[
\lambda(x|z)=\lambda_0(x)\,\mathrm{e}^{\beta'_0z},\quad x\in\R^+,
\]
where $\lambda_0$ represents the baseline hazard function, corresponding to a subject with $z=0$, and $\beta_0\in\R^p$ is the vector of the regression coefficients.

Let $H$ and $H^{uc}$ denote respectively the distribution function of the follow-up time and the sub-distribution function of the uncensored observations, i.e.,
\begin{equation}
\label{eq:def Huc}
H^{uc}(x)=\p(T\leq x,\Delta=1)=\int \delta\1_{\{t\leq x\}}\,\mathrm{d}\mathbb{P}(t,\delta,z),
\end{equation}
where $\p$ is the distribution of $(T,\Delta,Z)$.
We also require the following assumptions, some of which are common in large sample studies of the Cox model (e.g. see~\cite{LopuhaaNane2013}):
\begin{itemize}
\item[(A1)]
Let $\tau_F,\,\tau_G$ and $\tau_H$ be the end points of the support of $F,\,G$ and $H$. Then
\[
\tau_H=\tau_G<\tau_F\leq\infty.
\]
\item[(A2)]
There exists $\epsilon>0$ such that
\[
\sup_{|\beta-\beta_0|\leq\epsilon}\E\left[|Z|^2\,\mathrm{e}^{2\beta'Z}\right]<\infty.
\]
\item[(A3)]
There exists $\epsilon>0$ such that
\[
\sup_{|\beta-\beta_0|\leq\epsilon}\E\left[|Z|^2\,\mathrm{e}^{4\beta'Z}\right]<\infty.
\]
\end{itemize}
Let us briefly comment on these assumptions.
While the first one tells us that, at the end of the study, there is at least one subject alive, the other
two are somewhat hard to justify from a practical point of view. 
One can think of (A2) and (A3) as conditions on the boundedness of the second moment of the covariates,
uniformly for $\beta$ in a neighborhood of~$\beta_0$. 

By now, it seems to be rather a standard choice estimating $\beta_0$ by $\hat{\beta}_n$, the maximizer of the partial likelihood function, as proposed by~\cite{Cox72}.
The asymptotic behavior was first studied by~\cite{Tsiatis81}.
We aim at estimating $\lambda_0$, subject to the constraint that it is increasing (the case of a decreasing hazard is analogous), on the basis of $n$ observations $(T_1,\Delta_1,Z_1),\dots,(T_n,\Delta_n,Z_n)$.  We  refer to the quantity
\[
\Lambda_0(t)=\int_0^t\lambda_0(u)\,\mathrm{d}u,
\]
as the cumulative baseline hazard and, by introducing
\begin{equation}
\label{eq:def Phi}
\Phi(x;\beta)=\int \1_{ \{t\geq x\}}\,\mathrm{e}^{\beta'z}\,\mathrm{d}\p(t,\delta,z),
\end{equation}
we have
\begin{equation}
\label{eqn:lambda0}
\lambda_0(x)
=
\frac{h(x)}{\Phi(x;\beta_0)},
\end{equation}
where $h(x)=\mathrm{d}H^{uc}(x)/\mathrm{d}x$
(e.g., see (9) in~\cite{LopuhaaNane2013}).
For $\beta\in\R^p$ and $x\in\R$, the function $\Phi(x;\beta)$ can be estimated by
\begin{equation}
\label{eq:def Phin}
\Phi_n(x;\beta)=\int \1_{\{t\geq x\}} \mathrm{e}^{\beta'z}\,\mathrm{d}\p_n(t,\delta,z),
\end{equation}
where $\p_n$ is the empirical measure of the triplets $(T_i,\Delta_i,Z_i)$ with $i=1,\dots,n.$  Moreover, in Lemma 4 of~\cite{LopuhaaNane2013} it is shown that
\begin{equation}
\label{eqn:Phi}
\sup_{x\in\R}|\Phi_n(x;\beta_0)-\Phi(x;\beta_0)|=O_p(n^{-1/2}).
\end{equation}
It will be often used throughout the paper that a stochastic bound of the same order holds also for the distance between the cumulative hazard  $\Lambda_0$ and  the Breslow estimator
\begin{equation}
\label{eq:Breslow}
\Lambda_n(x)=\int \frac{\delta\1_{\{ t\leq x\}}}{\Phi_n(t;\hat{\beta}_n)}\,\mathrm{d}\p_n(t,\delta,z),
\end{equation}
but only on intervals staying away of the right boundary, i.e.,
\begin{equation}
\label{eqn:Breslow}
\sup_{x\in[0,M]}|\Lambda_n(x)-\Lambda_0(x)|=O_p(n^{-1/2}), \qquad\text{for all }0<M<\tau_H,
\end{equation}
(see Theorem 5 in~\cite{LopuhaaNane2013}).

Smoothing is done by means of kernel functions.
We will consider kernel functions $k$ that are $m$-orthogonal, for some $m\geq 1$,
which means that~$\int |k(u)||u|^m\,\mathrm{d}u<\infty$ and
$\int k(u)u^j\,\mathrm{d}u=0$, for $j=1,\ldots,m-1$, if $m\geq 2$.
We assume that
\begin{equation}
\label{def:kernel}
\begin{split}
&
\text{$k$ has bounded support $[-1,1]$ and is such that $\int_{-1}^1 k(y)\,\mathrm{d}y=1$;}\\
&
\text{$k$ is differentiable with a uniformly bounded derivative.}
\end{split}
\end{equation}
We denote by $k_b$ its scaled version $k_b(u)=b^{-1}k(u/b)$.
Here $b=b_n$ is a bandwidth that depends on the sample size, in such a way that
$0<b_n\to 0$ and $nb_n\to\infty$, as $n\to\infty$.
From now on, we will simply write $b$ instead of $b_n$.
Note that if $m>2$, the kernel function $k$ necessarily attains negative values and as a result also the smooth estimators of the baseline hazard defined in
Sections~\ref{sec:SG} may be negative and monotonicity might not be preserved.
To avoid this, one could restrict oneself to $m=2$.
In that case, the most common choice is to let $k$ be a symmetric probability density.


\section{Smoothed isotonic estimators}
\label{sec:SG}
We consider smoothed versions of two isotonic estimators for $\lambda_0$,
i.e, the maximum likelihood estimator~$\hat\lambda_n$ and the Grenander-type estimator~$\tilde\lambda_n$,
introduced in~\cite{LopuhaaNane2013}.
The maximum likelihood estimator of a nondecreasing baseline hazard rate $\lambda_0$
can be characterized as the left derivative of the greatest convex minorant of the cumulative sum diagram
consisting of points
$P_0=(0,0)$ and
$P_j=\big(\hat{W}_n(T_{(j+1)}),V_n(T_{(j+1)})\big)$,
for $j=1,\ldots,n-1$,
where~$\hat{W}_n$ and $V_n$ are defined as
\begin{equation}
\label{eq:def Wn Vn}
\begin{split}
\hat{W}_n(x)
&=
\int
\left(
\mathrm{e}^{\hat{\beta}'_nz}
\int_{T_{(1)}}^x\1_{\{u\geq s\}}\,\mathrm{d}s \right)\,\mathrm{d}\p_n(u,\delta,z),\quad x\geq T_{(1)},\\
V_n(x)
&=
\int\delta\1_{\{u<x\}}\,\mathrm{d}\p_n(u,\delta,z),
\end{split}
\end{equation}
with $\hat\beta_n$ being the partial maximum likelihood estimator (see Lemma~1 in~\cite{LopuhaaNane2013}).
For a fixed $x\in[0,\tau_H]$, the smoothed maximum likelihood estimator $\hat{\lambda}^{SM}_n$ of a nondecreasing baseline hazard rate $\lambda_0$,
was defined in~\cite{Nane} by
\begin{equation}
\label{def:lambdaSM}
\hat{\lambda}_n^{SM}(x)
=
\int_{(x-b)\vee 0}^{(x+b)\wedge\tau_H} k_b(x-u)\,\hat{\lambda}_n(u)\,\mathrm{d}u.
\end{equation}

The Grenander-type estimator $\tilde{\lambda}_n$ of a nondecreasing baseline hazard rate $\lambda_0$
is defined as the left hand slope of the greatest convex minorant (GCM) $\tilde{\Lambda}_n$ of the Breslow estimator $\Lambda_n$.
For a fixed $x_0\in[0,\tau_H]$, we consider the smoothed Grenander-type estimator $\tilde{\lambda}_n^{SG}$, which is defined by
\begin{equation}
\label{def:lambdaSG}
\tilde{\lambda}_n^{SG}(x)
=
\int_{(x-b)\vee 0}^{(x+b)\wedge\tau_H} k_b(x-u)\tilde{\lambda}_n(u)\,\mathrm{d}u.
\end{equation}
Uniform strong consistency on compact intervals in the interior of the support $[\epsilon,M]\subset[0,\tau_H]$ is provided by Theorem 5.2 of~\cite{Nane},
\begin{equation}
\label{eq:unif conv SG}
\sup_{x\in[\epsilon,M]}
\left|
\tilde{\lambda}_n^{SG}(x)-\lambda_0(x)
\right|
\to
0,
\quad
\text{with probability one.}%
\end{equation}
Strong pointwise consistency of $\hat{\lambda}_n^{SM}$ in the interior of the support is established in Theorem~5.1 in~\cite{Nane}.
Under additional smoothness assumptions on $\lambda_0$, one can obtain uniform strong consistency for $\hat{\lambda}_n^{SM}$ similar to~\eqref{eq:unif conv SG}.
Inconsistency at the boundaries is a frequently encountered problem in such situations and can be partially avoided by using a boundary corrected kernel.
One possibility is to construct linear combinations of $k(u)$ and $uk(u)$ with coefficients depending on the value near the boundary
(e.g., see~\cite{zhangkarunamuni1998}, \cite{DGL13}, or~\cite{LopuhaaMustaSN17}).
Then, it can be proved, exactly as it is done in~\cite{LopuhaaMustaSN17}, that uniform consistency holds on $[0,M]\subset[0,\tau_H]$.

\begin{figure}[t]
\includegraphics[bb=28 124 932 544,width=\textwidth,clip=]{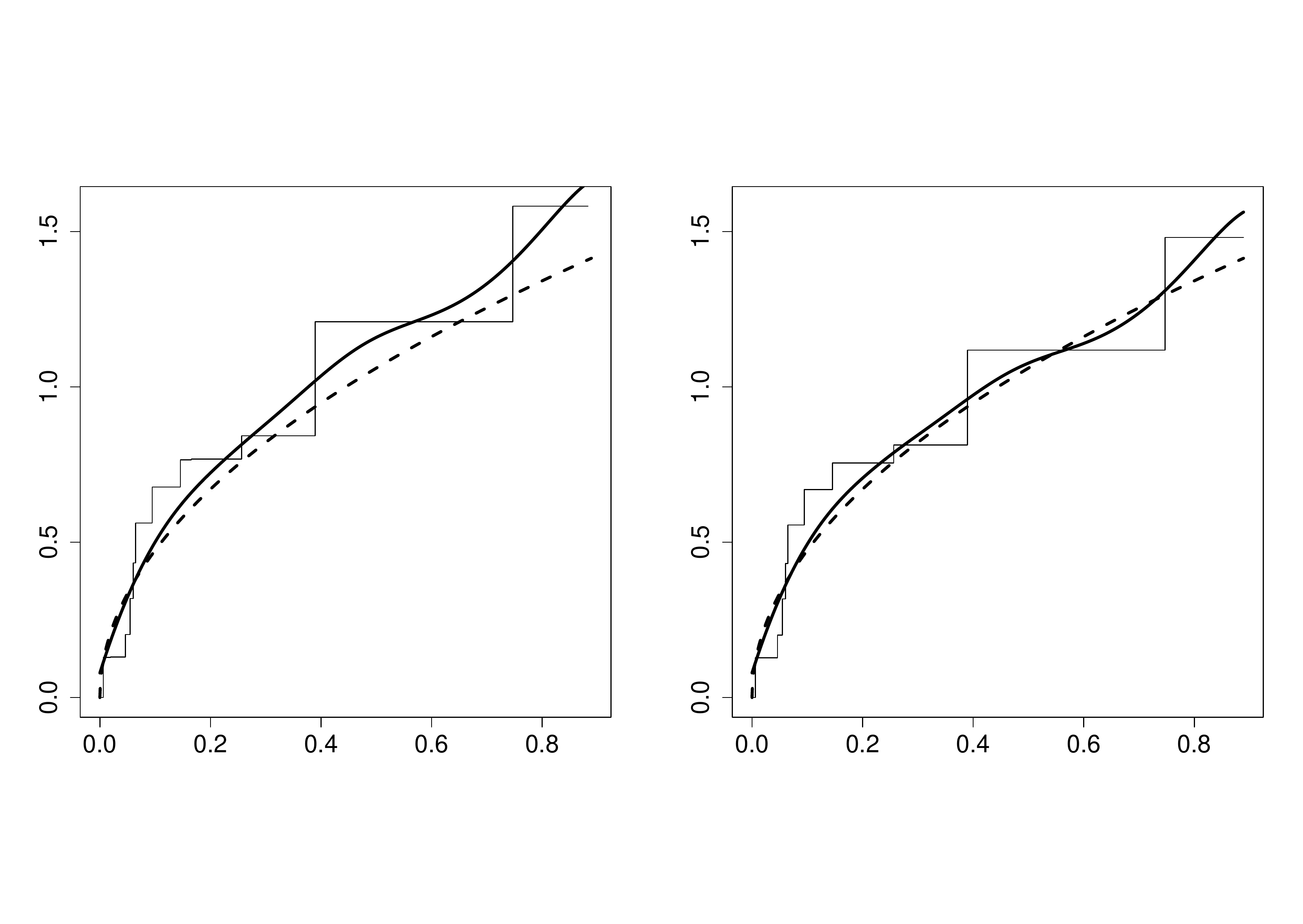}
\caption{Left panel: The MLE (piecewise constant solid line) of the baseline hazard (dashed) together with the smoothed MLE (solid).
Right panel: The Grenander estimator (piecewise constant solid line) of the baseline hazard (dashed)
together with the smoothed Grenander estimator (solid).}
\label{fig:SM}
\end{figure}
Figure~\ref{fig:SM} shows the smoothed maximum likelihood estimator (left)
and the smoothed Grenander-type estimator (right) for a sample of size $n=500$ from a Weibull baseline distribution with shape parameter $1.5$ and scale $1$.
For simplicity, we assume that the real valued covariate and the censoring times are uniformly $(0,1)$ distributed and we take $\beta_0=0.5$.
We used a boundary corrected triweight kernel function $k(u)=(35/32)(1-u^2)^3\1_{\{|u|\leq 1\}}$ and bandwidth $b=n^{-1/5}$.

In the remainder of this section we will derive the pointwise asymptotic distribution of both smoothed isotonic estimators,
in~\eqref{def:lambdaSM} and~\eqref{def:lambdaSG}.
As already mentioned, our approach is inspired by techniques introduced in Section 11.6 of~\cite{GJ14}.
We briefly describe this approach for the smoothed Grenander estimator, for which the computations are more complicated.
We start by writing
\begin{equation}
\label{eqn:lambda_SG_decomposition}
\tilde{\lambda}^{SG}_n(x)
=
\int k_b(x-u)\,\mathrm{d}\Lambda_0(u)+\int k_b(x-u)\,\mathrm{d}(\tilde{\Lambda}_n-\Lambda_0)(u).
\end{equation}
The first (deterministic) term on the right hand side of~\eqref{eqn:lambda_SG_decomposition} gives us the asymptotic bias. 
The method applied in~\cite{LopuhaaMustaSN17} for the right censoring model continues by decomposing the second term in two parts
\[
\int k_b(x-u)\,\mathrm{d}(\tilde{\Lambda}_n-\Lambda_n)(u)+\int k_b(x-u)\,\mathrm{d}(\Lambda_n-\Lambda_0)(u),
\]
and then uses the Kiefer-Wolfowitz type of result 
\begin{equation}
\label{eqn:KW}
\sup_{t\in[0,\tau_H]}|\tilde{\Lambda}_n(t)-\Lambda_n(t)|=O_P\left(n^{-2/3}(\log n)^{2/3} \right),
\end{equation}
to show that $\int k_b(x-u)\,\mathrm{d}(\tilde{\Lambda}_n-\Lambda_n)(u)$ converges to zero. 
{Finally}, results from empirical process theory are used to show the asymptotic normality of $\int k_b(x-u)\,\mathrm{d}(\Lambda_n-\Lambda_0)(u)$.
This approach {cannot} be followed in our case because of the lack of a Kiefer-Wolfowitz type of result as in~\eqref{eqn:KW} for the Cox model. 

{Alternatively, we proceed by describing the main steps of the $L_2$-bounds approach introduced in~\cite{GJ14}.}
On an event $E_n$ with probability tending to one, we will approximate
\begin{equation}
\label{eq:term2 SG}
\int k_b(x-u)\,\mathrm{d}(\tilde{\Lambda}_n-\Lambda_0)(u)
\end{equation}
by
$\int
\theta_{n,x}(u,\delta,z)
\,\mathrm{d}\mathbb{P}(u,\delta,z)$,
for some suitable function $\theta_{n,x}$ (see Lemma~\ref{le:2-1}),
whose piecewise constant modification $\overline{\theta}_{n,x}$ integrates to zero with respect to the empirical measure $\mathbb{P}_n$
(see Lemma~\ref{le:2-2}).
This enables us to approximate~\eqref{eq:term2 SG} by
\begin{equation}
\label{eq:decomp}
\int
\overline{\theta}_{n,x}(u,\delta,z)\,\mathrm{d}(\mathbb{P}_n-\mathbb{P})(u,\delta,z)
+
\int
\left(
\overline{\theta}_{n,x}(u,\delta,z)-\theta_{n,x}(u,\delta,z)
\right)\,\mathrm{d}\mathbb{P}(u,\delta,z).
\end{equation}
Then, the key step is to bound the second integral in~\eqref{eq:decomp} by means of $L_2$-bounds on the distance between the ordinary Grenander estimator and the true baseline hazard
(see Lemma~\ref{le:result4}).
The last step consists of replacing $\overline{\theta}_{n,x}$ by a deterministic function $\eta_{n,x}$
(see Lemma~\ref{le:result5}) and use empirical process theory to show that
\[
\int
\eta_{n,x}(u,\delta,z)
\,\mathrm{d}(\mathbb{P}_n-\mathbb{P})(u,\delta,z)
\]
is asymptotically normal.

Before we proceed to our first main result, we will formulate the steps described above in a series of lemmas.
Let $x\in(0,\tau_H)$ and $0<M<\tau_H$.
For $n$ sufficiently large, such that $0<x-b<x+b<M$, define
\begin{equation}
\label{eqn:an}
a_{n,x}(u)=\frac{k_b(x-u)}{\Phi(u;\beta_0)},
\quad
\text{for }u\leq x+b,
\end{equation}
where $\Phi(u;\beta_0)$ is defined in~\eqref{eq:def Phi},
and $a_{n,x}(u)=0$ for $u>x+b$.
We then have the following approximation for~\eqref{eq:term2 SG}.
The proof can be found in Section~\ref{sec:proofs}.
\begin{lemma}
\label{le:2-1}
Suppose that (A1)--(A2) hold.
Let $a_{n,x}$ be defined by~\eqref{eqn:an} and let $\hat\beta_n$ be the partial MLE for $\beta_0$.
There exists an event $E_n$, with $\1_{E_n}\to1$ in probability, such that for
\begin{equation}
\label{eqn:thetan}
\theta_{n,x}(u,\delta,z)
=
\1_{E_n}
\left\{
\delta\,a_{n,x}(u)-\mathrm{e}^{\hat{\beta}'_n\,z}\int_0^u a_{n,x}(v)\,\mathrm{d}\tilde{\Lambda}_n(v)
\right\},
\end{equation}
it holds
\[
\int\theta_{n,x}(u,\delta,z)\,\mathrm{d}\p(u,\delta,z)=-\1_{E_n}\int k_b(x-u)\,\mathrm{d}(\tilde{\Lambda}_n-\Lambda_0)(u)+O_p(n^{-1/2}).
\]
\end{lemma}

Next, we consider a piecewise constant modification $\overline{a}_{n,x}\overline{\Phi}_n$ of $a_{n,x}\Phi_n$, which is constant on the same intervals as $\tilde{\lambda}_n$. Let $\tau_0=x-b$, $\tau_{m+1}=x+b$ and let $(\tau_i)_{i=1}^m$ be successive points of jump of $\tilde{\lambda}_n$ in the interval $(x-b,x+b)$.
Then, for $u\in(\tau_i,\tau_{i+1}]$, we choose
\begin{equation}
\label{def:piecewise constant}
\overline{a}_{n,x}\overline{\Phi}_n(u;\hat{\beta}_n)=a_{n,x}(\hat{A}_n(u))\Phi_n(\hat{A}_n(u);\hat{\beta}_n),
\end{equation}
where for $u\in(\tau_i,\tau_{i+1}]$,
\begin{equation}
\label{eq:An}
\hat{A}_n(u)
=
\begin{cases}
\tau_i, & \text{ if }\lambda_0(t)>\tilde{\lambda}_n(\tau_{i+1}),\text{ for all }t\in(\tau_i,\tau_{i+1}],\\
s, & \text{ if }\lambda_0(s)=\tilde{\lambda}_n(s),\text{ for some }s\in(\tau_i,\tau_{i+1}],\\
\tau_{i+1}, & \text{ if }\lambda_0(t)<\tilde{\lambda}_n(\tau_{i+1}),\text{ for all }t\in(\tau_i,\tau_{i+1}].
\end{cases}
\end{equation}
Furthermore, let $E_n$ be the event from Lemma~\ref{le:2-1} and define
\begin{equation}
\label{def:Psibar}
\overline{\Psi}_{n,x}(u)=\frac{\overline{a}_{n,x}\overline{\Phi}_n(u;\hat{\beta}_n)}{\Phi_n(u;\hat{\beta}_n)}\,\1_{E_n}, \quad u\in[x-b,x+b],
\end{equation}
and $\overline{\Psi}_{n,x}(u)=0$, for $u\neq[x-b,x+b]$.
Note that, since $u\leq x+b<M<T_{(n)}$ on the event~$E_n$, we have $\Phi_n(u;\hat\beta_n)>0$
(see the proof of Lemma~\ref{le:2-1}),
and thus $\overline{\Psi}_{n,x}(u)$ is well defined.
Now, define the following piecewise constant modification of $\theta_{n,x}$, by
\begin{equation}
\label{eqn:thetabar}
\overline{\theta}_{n,x}(u,\delta,z)=\delta\,\overline{\Psi}_{n,x}(u)-\mathrm{e}^{\hat{\beta}'_n\,z}\int_0^u \overline{\Psi}_{n,x}(v)\,\mathrm{d}\tilde{\Lambda}_n(v).
\end{equation}
We then have the following property.
The proof can be found in Section~\ref{sec:proofs}.
\begin{lemma}
\label{le:2-2}
Let $\overline{\theta}_{n,x}$ be defined in~\eqref{eqn:thetabar}.
Then
\begin{equation}
\label{eqn:int}
\int\overline{\theta}_{n,x}(u,\delta,z)\,\mathrm{d}\p_n(u,\delta,z)=0.
\end{equation}
\end{lemma}
At this point it is important to discuss in some detail how we will obtain suitable bounds for the second integral in~\eqref{eq:decomp}.
In order to do so, we first introduce the inverse process~$\tilde U_n$.
It is defined by
\begin{equation}
\label{eq:def Un}
\tilde{U}_n(a)
=
\argmin_{x\in[0,T_{(n)}]} \left\{\Lambda_n(x)-ax\right\}.
\end{equation}
and it satisfies the switching relation, $\tilde{\lambda}_n(x)\leq a$ if and only if $\tilde{U}_n(a)\geq x$, for $x\leq T_{(n)}$.
In their analysis of the current status model, \cite{GJW10} encounter an integral that is similar to the second integral in~\eqref{eq:decomp}.
They bound this integral using that the maximal distance between succeeding points of jump of the isotonic estimator is of the order $O_p(n^{-1/3}\log n)$.
Such a property typically relies on the exponential bounds for the tail probabilities of~$\tilde{U}_n(a)$,
obtained either by using a suitable exponential martingale (e.g., see Lemma~5.9 in~\cite{groeneboom-wellner1992}),
or by an embedding of the relevant sum process into Brownian motion or Brownian bridge (e.g., see Lemma~5.1 in~\cite{DurotKulikovLopuhaa2012}).
Unfortunately, an embedding of the process~$\Lambda_n$ is not available and in our current situation
the martingale approach only yields to polynomial bounds for tail probabilities of~$\tilde{U}_n(a)$.
A polynomial bound was also found by~\cite{durot2007} (see her Lemma~2) leading to
\begin{equation}
\label{eq:sup E}
\sup_{x\in I_n}
\mathbb{E}
\left[
\big(\tilde\lambda_n(x)-\lambda_0(x)\big)^p
\right]
\leq Kn^{-p/3},
\end{equation}
for $p\in[1,2)$ and some interval $I_n$ (see her Theorem~1).
By intersecting with the event~$E_n$ from Lemma~\ref{le:2-1} we extend~\eqref{eq:sup E} to a similar bound for $p=2$.
\cite{GJ14} provide an alternative approach to bound the second integral in~\eqref{eq:decomp}, based on bounds for~\eqref{eq:sup E} with $p=2$.
Unfortunately, they still make use of the fact that the maximum distance between succeeding points of jump of the isotonic estimator is of the order~$O_p(n^{-1/3}\log n)$
to obtain a result similar to~\eqref{eqn:statement}.
Nevertheless, we do follow the approach in~\cite{GJ14}, but instead of using the maximum distance between succeeding points of jump of $\tilde\lambda_n$,
we use a bound on
\begin{equation}
\label{eq:E sup}
\mathbb{E}\left[
\sup_{x\in [\epsilon,M]}
\big(\tilde\lambda_n(x)-\lambda_0(x)\big)^2
\right],
\end{equation}
for $0<\epsilon<M<\tau_H$.
Exponential bounds for the tail probabilities of~$\tilde{U}_n(a)$ would yield the same bound for~\eqref{eq:E sup}
as the one in~\eqref{eq:sup E} apart from a factor $\log n$.
Since we can only obtain polynomial bounds on the tail probabilities of~$\tilde{U}_n(a)$,
we establish a bound for~\eqref{eq:E sup} of the order~$O(n^{-4/9})$.
This is probably not optimal, but this
will turn out to be sufficient for our purposes and leads to the following intermediate result,
of which the proof can be found in Section~\ref{sec:proofs}.
\begin{lemma}
\label{le:result4}
Suppose that (A1)--(A2) hold. 
Fix $x\in(0,\tau_h)$ and let $\theta_{n,x}$ and $\bar\theta_{n,x}$ be defined by~\eqref{eqn:thetan} and~\eqref{eqn:thetabar}, respectively.
Assume that $\lambda_0$ is differentiable, such that $\lambda'_0$ is uniformly bounded above and below by strictly positive constants.
Assume that $x\mapsto \Phi(x;\beta_0)$ is differentiable with a bounded derivative in a neighborhood of $x$
{and let} $k$ satisfy~\eqref{def:kernel}.
Then, it holds
\[
\int
\left\{
\overline{\theta}_{n,x}(u,\delta,z)-\theta_{n,x}(u,\delta,z)
\right\}\,\mathrm{d}\p(u,\delta,z)=
O_p(b^{-1}n^{-2/3}).
\]
\end{lemma}
The last step is to replace $\overline{\theta}_{n,x}$ in the first integral of~\eqref{eq:decomp} with a deterministic approximation.
This is done in the next lemma,
of which the proof can be found in Section~\ref{sec:proofs}.
\begin{lemma}
\label{le:result5}
Suppose that (A1)--(A3) hold.
Fix $x\in(0,\tau_h)$ and take $0<\epsilon<x<M'<M<\tau_H$.
Assume that $\lambda_0$ is differentiable, such that $\lambda'_0$ is uniformly bounded above and below by strictly positive constants.
Assume that $x\mapsto \Phi(x;\beta_0)$ is differentiable with a bounded derivative in a neighborhood of $x$.
Let $\bar\theta_{n,x}$ be defined in~\eqref{eqn:thetabar} and define
\begin{equation}
\label{eqn:eta}
\eta_{n,x}(u,\delta,z)=\1_{E_n}
\left(\delta\,a_{n,x}(u)-\mathrm{e}^{\beta'_0z}\,\int_0^u a_{n,x}(v)\,\mathrm{d}\Lambda_0(v)\right),\quad u\in[0,\tau_H].
\end{equation}
where $a_{n,x}$ is defined in~\eqref{eqn:an} and $E_n$ is the event from Lemma~\ref{le:2-1}.
Let $k$ satisfy~\eqref{def:kernel}.
Then, it holds
\begin{equation}
\label{eqn:statement}
\int
\left\{
\overline{\theta}_{n,x}(u,\delta,z)-\eta_{n,x}(u,\delta,z)
\right\}\,\mathrm{d}(\p_n-\p)(u,\delta,z)
=
O_p(b^{-3/2}n^{-13/18})+O_p(n^{-1/2})
\end{equation}
\end{lemma}
We are now in the position to state our first main result.
\begin{theo}
\label{theo:distr}
Suppose that (A1)--(A3) hold. 
Fix $x\in(0,\tau_h)$.
Assume that $\lambda_0$ is $m\geq 2$ times continuously differentiable in $x$,
such that $\lambda'_0$ is uniformly bounded above and below by strictly positive constants.
Moreover, assume that $t\mapsto \Phi(t;\beta_0)$ is differentiable with a bounded derivative in a neighborhood of~$x$
{and let} $k$ satisfy~\eqref{def:kernel}. 
Let~$\tilde{\lambda}^{SG}$ be defined in~\eqref{def:lambdaSG} and assume that $n^{1/(2m+1)}b\to c>0$.
Then, it holds
\[
n^{m/(2m+1)}
\left(
\tilde{\lambda}_n^{SG}(x)-\lambda_0(x)
\right)\xrightarrow{d}N(\mu,\sigma^2),
\]
where
\begin{equation}
\label{eqn:as.mean-var}
\mu=\frac{(-c)^m}{m!}\lambda^{(m)}_0(x)\int_{-1}^1 k(y)y^m\,\mathrm{d}y
\quad\text{ and }\quad
\sigma^2
=
\frac{\lambda_0(x)}{c\Phi(x;\beta_0)}
\int k^2(u)\,\mathrm{d}u.
\end{equation}
Furthermore,
\begin{equation}
\label{eq:equivalence SG SM}
n^{m/(2m+1)}
\left(
\tilde{\lambda}_n^{SG}(x)-\tilde{\lambda}_n^{SM}(x)
\right)
\to0,
\end{equation}
in probability, where $\tilde{\lambda}_n^{SM}(x)$ is defined in~\eqref{def:lambdaSM},
so that $\tilde{\lambda}_n^{SM}(x)$ has the same limiting distribution as~$\tilde{\lambda}_n^{SG}(x)$.
\end{theo}
\begin{proof}
Choose $0<\epsilon<x<M'<M<\tau_H$, so that for $n$ sufficiently large, we have $\epsilon<x-b\leq x+b\leq M'$.
Consider the event~$E_n$ from Lemma~\ref{le:2-1} and choose $\xi_1,\xi_2>0$ and~$\xi_3$, such that it satisfies~\eqref{eqn:xi}.
We write
\begin{equation}
\label{eq:decomposition lambdaSG}
\begin{split}
\tilde{\lambda}^{SG}_n(x)
&=
\int k_b(x-u)\,\mathrm{d}\tilde{\Lambda}_n(u)\\
&=
\int k_b(x-u)\,\mathrm{d}\Lambda_0(u)+\1_{E_n}\,\int k_b(x-u)\,\mathrm{d}(\tilde{\Lambda}_n-\Lambda_0)(u)\\
&\qquad+
\1_{E_n^c}\,\int k_b(x-u)\,\mathrm{d}(\tilde{\Lambda}_n-\Lambda_0)(u).
\end{split}
\end{equation}
Because $\1_{E_n^c}\to 0$ in probability, the third term on the right hand side tends to zero in probability.
For the first term, we obtain from a change of variable, a Taylor expansion, and the properties of the kernel:
\begin{equation}
\label{eqn:asymptotic-mean}
\begin{split}
&
n^{m/(2m+1)}
\left\{
\int k_b(x-u)\,\lambda_0(u)\,\mathrm{d}u-\lambda_0(x)
\right\}\\
&=
n^{m/(2m+1)}
\int_{-1}^1 k(y)
\left\{
\lambda_0(x-by)-\lambda_0(x)
\right\}\,\mathrm{d}y\\
&=
n^{m/(2m+1)}
\int_{-1}^1 k(y)
\left\{
-\lambda'_0(x)by+\cdots+\frac{\lambda^{(m-1)}_0(x)}{(m-1)!}(-by)^{m-1}+\frac{\lambda^{(m)}_0(\xi_n)}{m!}(-by)^m
\right\}\,\mathrm{d}y\\
&\to
\frac{(-c)^m}{m!}\lambda^{(m)}_0(x)\int_{-1}^1 k(y)y^m\,\mathrm{d}y,
\end{split}
\end{equation}
with $|\xi_n-x|<b|y|$.
Finally, for the second term on the right hand side of~\eqref{eq:decomposition lambdaSG},
Lemmas~\ref{le:2-1} to~\ref{le:result5} yield that
\begin{equation}
\label{eqn:key}
\begin{split}
&n^{m/(2m+1)}
\1_{E_n}
\int k_b(x-u)\,\mathrm{d}(\tilde{\Lambda}_n-\Lambda_0)(u)\\
&=
n^{m/(2m+1)}
\int\eta_{n,x}(u,\delta,z)\,\mathrm{d}(\p_n-\p)(u,\delta,z)+o_p(1).
\end{split}
\end{equation}
For the first term on the right hand side of~\eqref{eqn:key} we can write
\begin{equation}
\label{eq:decomp eta}
\begin{split}
&
n^{m/(2m+1)}
\int\eta_{n,x}(u,\delta,z)\,\mathrm{d}(\p_n-\p)(u,\delta,z)\\
&=
n^{m/(2m+1)}
\1_{E_n}
\int\frac{\delta k_b(x-u)}{\Phi(u;\beta_0)}\,\mathrm{d}(\p_n-\p)(u,\delta,z)\\
&\quad
-
n^{m/(2m+1)}
\1_{E_n}
\int
\mathrm{e}^{\beta'_0z}\int_0^u a_{n,x}(v)\,\mathrm{d}\Lambda_0(v)\,\mathrm{d}(\p_n-\p)(u,\delta,z).
\end{split}
\end{equation}
We will show that the first term on the right hand is asymptotically normal and the second term tends to zero in probability.
Define
$Y_{n,i}
=
n^{-(m+1)/(2m+1)}\Delta_i k_b(x-T_i)/\Phi(T_i;\beta_0)$,
so that the first term on the right hand side of~\eqref{eq:decomp eta} can be written as
\[
\1_{E_n}
n^{m/(2m+1)}
\int\frac{\delta k_b(x-u)}{\Phi(u;\beta_0)}\,\mathrm{d}(\p_n-\p)(u,\delta,z)
=
\1_{E_n}
\sum_{i=1}^n
\left(Y_{n,i}-\E\left[Y_{n,i}\right]
\right).
\]
Using~\eqref{eqn:lambda0}, together with a Taylor expansion and the boundedness assumptions on the derivatives of
$\lambda_0$ and $\Phi(x;\beta_0)$, we have
\begin{equation}
\label{eqn:asymptotic_variance}
\begin{split}
&
\sum_{i=1}^n \text{Var}(Y_{n,i})\\
&=
n^{-1/(2m+1)}
\left\{
\int\frac{k^2_b(x-u)}{\Phi(u;\beta_0)^2}\,\mathrm{d}H^{uc}(u)-\left(\int\frac{k_b(x-u)}{\Phi(u;\beta_0)}\,\mathrm{d}H^{uc}(u)\right)^2
\right\}\\
&=
n^{-1/(2m+1)}
\left\{
\frac{1}{b}\int_{-1}^1k^2(y)\frac{\lambda_0(x-by)}{\Phi(x-by;\beta_0)}\,\mathrm{d}y-\left(\int k_b(x-u)\,\lambda_0(u)\,\mathrm{d}u\right)^2
\right\}\\
&=
\frac{\lambda_0(x)}{c\Phi(x;\beta_0)}
\int_{-1}^1k^2(y)\,\mathrm{d}y
-
n^{-1/(2m+1)}
\int_{-1}^1yk^2(y)
\left[
\frac{\mathrm{d}}{\mathrm{d}x}
\frac{\lambda_0(x)}{\Phi(x;\beta_0)}
\right]_{x=\xi_y}
\,\mathrm{d}y+o(1)\\
&=
\frac{\lambda_0(x)}{c\Phi(x;\beta_0)}
\int_{-1}^1k^2(y)\,\mathrm{d}y+o(1).
\end{split}
\end{equation}
Moreover, $|Y_{n,i}|\leq n^{-(m+1)/(2m+1)}\Phi(M;\beta_0)^{-1}\sup_{x\in[-1,1]}k(x)$,
so that
$\sum_{i=1}^n
\E
\big[
|Y_{n,i}|^2\1_{\{|Y_{n,i}|>\epsilon\}}
\big]
\to 0$,
for any $\epsilon>0$,
since $\1_{\{|Y_{n,i}|>\epsilon\}}=0$, for $n$ sufficiently large.
Consequently, by Lindeberg central limit theorem, and the fact that $\1_{E_n}\to1$ in probability, we obtain
\begin{equation}
\label{eq:asymp norm term1}
\1_{E_n}
n^{m/(2m+1)}
\int\frac{\delta k_b(x-u)}{\Phi(u;\beta_0)}\,\mathrm{d}(\p_n-\p)(u,\delta,z)
\to
N(0,\sigma^2).
\end{equation}
For the second term on the right hand side of~\eqref{eq:decomp eta}, write
\[
n^{m/(2m+1)}\int\mathrm{e}^{\beta'_0z}\int_0^u a_{n,x}(v)\,\mathrm{d}\Lambda_0(v)\,\mathrm{d}(\p_n-\p)(u,\delta,z)
=
\sum_{i=1}^n
\left(\widetilde{Y}_{n,i}-\E [\widetilde{Y}_{n,i}]\right).
\]
where
 \[
\widetilde{Y}_{n,i}
=
n^{-(m+1)/(2m+1)}\mathrm{e}^{\beta'_0Z_i}\,\int_0^{T_i}\frac{k_b(x-v)}{\Phi(v;\beta_0)}\,\mathrm{d}\Lambda_0(v).
\]
We have
\[
\sum_{i=1}^n \text{Var}(\widetilde{Y}_{n,i})
\leq
\sum_{i=1}^n
\mathbb{E}\left[\widetilde{Y}_{n,i}^2\right]
\leq
n^{-1/(2m+1)}
\int \mathrm{e}^{2\beta'_0z}
\left(
\int_0^u\frac{k_b(x-v)}{\Phi(v;\beta_0)}\,\mathrm{d}\Lambda_0(v)
\right)^2\,\mathrm{d}\p(u,\delta,z),
\]
where the integral on the right hand side is bounded by
\[
\left(
\int_{x-b}^{x+b}\frac{k_b(x-v)}{\Phi(v;\beta_0)}\,\mathrm{d}\Lambda_0(v)
\right)^2
\Phi(0;2\beta_0)
\leq
\frac{\Phi(0;2\beta_0)}{\Phi^2(M;\beta_0)}
\left(
\int_{x-b}^{x+b}k_b(x-v)\,\mathrm{d}\Lambda_0(v)
\right)^2
=O(1).
\]
Hence, the second term on the right hand side of~\eqref{eq:decomp eta} tends to zero in probability.
Together with~\eqref{eq:decomposition lambdaSG}, \eqref{eqn:asymptotic-mean}, and~\eqref{eq:asymp norm term1},
this proves the first part of the theorem.

For the smoothed maximum likelihood estimator, we can follow the same approach and obtain similar results as those in Lemmas~\ref{le:2-1} to~\ref{le:result5}.
The arguments are more or less the same as those used to prove Lemmas~\ref{le:2-1} to~\ref{le:result5}. 
{We briefly sketch the main differences.
First, instead of $\tilde{\Lambda}_n$, we now use  
\[
\hat{\Lambda}_n(x)=\int_0^x\hat{\lambda}_n(u)\,\mathrm{d} u
\]
in~\eqref{eq:term2 SG}.}
Then, since the maximum likelihood estimator {is defined as the left slope of the greatest convex minorant of a cumulative sum diagram
that is different from the one corresponding to the Grenander-type estimator,} 
{Lemmas~\ref{le:2-1} and~\ref{le:2-2}} will hold with a different event $\widehat{E}_n$ and~$\overline{\Psi}_{n,x}$ will have a simpler form 
(see {Lemmas~\ref{le:2-1 SM}-\ref{le:2-2 SM} and definition}~\eqref{def:Psibar SM} in~\cite{LopuhaaMustaSM2016}). 
{
Similar to the proof of Lemma~\ref{le:result4},
the proof of its counterpart for the maximum likelihood estimator (see Lemma~\ref{le:result4 SM} in~\cite{LopuhaaMustaSM2016})
is quite technical and involves bounds on the tail probabilities of the inverse process corresponding to~$\hat\lambda_n$
(see Lemma~\ref{le:invSM}), used to obtain the analogue of~\eqref{eq:E sup}
(see Lemma~\ref{le:lambdaSM}).}
Moreover, the inverse process related to the maximum likelihood estimator is defined by 
\begin{equation}
\label{eq:def Un SM}
\hat{U}_n(a)
=
\argmin_{x\in[T_{(1)},T_{(n)}]}
\left\{V_n(x)-a\hat{W}_n(x)\right\},
\end{equation}
{where $V_n$ and $\hat{W}_n$ are defined in~\eqref{eq:def Wn Vn}},
and we get a slightly different bound on the tail probabilities of $\hat{U}_n$ (compare Lemma~\ref{le:inv} and Lemma~\ref{le:invSM} in~\cite{LopuhaaMustaSM2016}). The reason is that  the martingale decomposition of the process $V_n(t)-a\hat{W}_n(t)$ has a simper form.
{The counterpart of Lemma~\ref{le:result5} (see Lemma~\ref{le:result5 SM} in~\cite{LopuhaaMustaSM2016})
is established in the same way, replacing $\tilde{\lambda}_n$ by $\hat{\lambda}_n$.}
For details we refer to Section~\ref{sec:SMLE} in~\cite{LopuhaaMustaSM2016}. 

From~\eqref{eq:decomposition lambdaSG} and~\eqref{eqn:key}, we have that
\begin{equation}
\label{eq:expansion SG}
\begin{split}
n^{m/(2m+1)}\tilde{\lambda}^{SG}_n(x)
&=
n^{m/(2m+1)}\int k_b(x-u)\,\mathrm{d}\Lambda_0(u)\\
&\quad
+
n^{m/(2m+1)}
\int
\eta_{n,x}(u,\delta,z)\,\mathrm{d}(\p_n-\p)(u,\delta,z)+o_p(1)
\end{split}
\end{equation}
where $\eta_{n,x}$ is defined in~\eqref{eqn:eta}
and where
\begin{equation}
\label{eq:asymp norm SG}
n^{m/(2m+1)}
\int
\eta_{n,x}(u,\delta,z)\,\mathrm{d}(\p_n-\p)(u,\delta,z)
\to
N(0,\sigma^2).
\end{equation}
Similarly, from the results in Section~\ref{sec:SMLE} of~\cite{LopuhaaMustaSM2016}, we have that there
exists an event $\widehat{E}_n$, such that
\begin{equation}
\label{eq:expansion SM}
\begin{split}
n^{m/(2m+1)}
\tilde{\lambda}^{SM}_n(x)
&=
n^{m/(2m+1)}\int k_b(x-u)\,\mathrm{d}\Lambda_0(u)\\
&\quad+
n^{m/(2m+1)}
\int
\widehat{\eta}_{n,x}(u,\delta,z)\,\mathrm{d}(\p_n-\p)(u,\delta,z)+o_p(1)
\end{split}
\end{equation}
where $\widehat{\eta}_{n,x}$ is defined in~\eqref{eqn:eta} with $\widehat E_n$ instead of $E_n$,
where $\1_{\widehat{E}_n}\to 1$ in probability,
and where
\begin{equation}
\label{eq:asymp norm SM}
n^{m/(2m+1)}
\int
\widehat{\eta}_{n,x}(u,\delta,z)\,\mathrm{d}(\p_n-\p)(u,\delta,z)
\to
N(0,\sigma^2).
\end{equation}
Together with~\eqref{eq:asymp norm SG} and~\eqref{eq:asymp norm SM},
this means that
\[
\begin{split}
&
n^{m/(2m+1)}
\left(
\tilde{\lambda}_n^{SG}(x)-\hat{\lambda}_n^{SM}(x)
\right)\\
&=
\left(
\1_{\widehat{E}_n^c}\1_{E_n}-\1_{E_n^c}\1_{\widehat{E}_n}
\right)\\
&\qquad\times
n^{m/(2m+1)}
\int
\left\{
\delta a_{n,x}(u)-\mathrm{e}^{\beta'_0z}\int_0^u a_{n,x}(v)\,\mathrm{d}\Lambda_0(v)
\right\}\,\mathrm{d}(\p_n-\p){(u,\delta,z)}+o_p(1)\\
&=
\1_{\widehat{E}_n^c}O_p(1)
-
\1_{E_n^c}O_p(1)
+o_p(1)
=
o_p(1),
\end{split}
\]
because $\1_{\widehat{E}_n^c}\to 0$ and $\1_{E_n^c}\to0$ in probability.
\end{proof}
Note that in the special case $\beta_0=0$ and $m=2$,
we recover Theorem 3.2 in~\cite{LopuhaaMustaSN17} and Theorem 11.8 in~\cite{GJ14},
for the right censoring model without covariates.
The fact that $\tilde{\lambda}_n^{SG}(x)$ and~$\hat{\lambda}_n^{SM}(x)$ are asymptotically equivalent does not come as a surprise,
since for the corresponding isotonic estimators according to Theorem 2 in~\cite{LopuhaaNane2013}, for $x\in(0,\tau_H)$ fixed,
{$n^{1/3}
\big(
\tilde\lambda_n(x)-\hat\lambda_n(x)
\big)
\to0$,
in probability.}
However, we have not been able to exploit this fact, and we have established the asymptotic equivalence in~\eqref{eq:equivalence SG SM}
by obtaining the expansions in~\eqref{eq:expansion SG} and~\eqref{eq:expansion SM} separately for each estimator.

{\begin{re}
\label{rem:partial MLE}
The estimators considered in Theorem~\ref{theo:distr} are based on the partial maximum likelihood estimator $\hat\beta_n$, which
defines the Breslow estimator, see~\eqref{eq:Breslow}, and the cumulative sum diagram  from which the SMLE is determined,
see~\eqref{eq:def Wn Vn}.
However, Theorem~\ref{theo:distr} remains true, if $\hat\beta_n$ is any estimator that satisfies 
\begin{equation}
\label{eq:cond PMLE}
\hat\beta_n-\beta_0\to0, \text{ a.s.,}
\quad\text{and}\quad
\sqrt{n}(\hat\beta_n-\beta_0)=O_p(1)
\end{equation}
In particular, this holds for the partial MLE for $\beta_0$.
See, e.g.,  Theorems~3.1 and~3.2 in~\cite{Tsiatis81}. 
When proving consistency of the bootstrap, we are not able to establish bootstrap versions
of Theorems~3.1 and~3.2 in~\cite{Tsiatis81}, but, in view of this remark, it is sufficient to assume the bootstrap version of~\eqref{eq:cond PMLE}.
\end{re}
}



\section{Numerical results for pointwise confidence intervals}
\label{sec:conf-int}

In this section we illustrate the finite sample performance of the two estimators considered previously by constructing pointwise confidence intervals for the baseline hazard rate.
We consider two different procedures: the first {one} relies on the limit distribution and the second {one} is a bootstrap based method. 
In all the simulations we use the triweight kernel function,
which means that the degree of smoothness is $m=2$. 
The reason {for choosing a second-order kernel} is that higher order kernels {may also take} negative values, which then might lead to non monotone estimators for the baseline hazard.

\subsection{Asymptotic confidence intervals}
\label{subsec:asympconf}
From Theorem~\ref{theo:distr} it can be seen that the asymptotic $100(1-\alpha)\%$-confidence intervals at the point $x_0\in(0,\tau_H)$
are of the form
\[
\lambda_n^{SI}(x_0)
-n^{-2/5}
\left\{
\widehat{\mu}_n(x_0)\pm\widehat{\sigma}_n(x_0)q_{1-\alpha/2}
\right\},
\]
where $q_{1-\alpha/2}$ is the $(1-\alpha/2)$ quantile of the standard normal distribution,
$\lambda_n^{SI}(x_0)$ is the smooth isotonic estimator at hand (SG or SMLE),
and $\widehat{\sigma}_n(x_0)$, $\widehat{\mu}_n(x_0)$ are corresponding plug-in estimators of the asymptotic mean and standard deviation, respectively. 
However, from the expression of the asymptotic mean in Theorem~\ref{theo:distr} for $m=2$, it is obvious that obtaining the plug-in estimators requires estimation of the second derivative of~$\lambda_0$.
Since accurate estimation of derivatives is a hard problem, we choose to avoid it by using undersmoothing.
This procedure is {to be preferred above} bias estimation,
because it is computationally more convenient and leads to better results
(see also~\cite{Hall92}, \cite{GJ15}, \cite{CHT06}).
Undersmoothing consists of using a bandwidth of a smaller order than the optimal one (in our case~$n^{-1/5}$).
As a result, the bias of $n^{2/5}( \lambda_n^{SI}(x_0)-\lambda_0(x_0))$, which is of the order $n^{2/5}b^2$ (see~\eqref{eqn:asymptotic-mean}), will converge to zero. On the other hand, the asymptotic variance is $n^{-1/5}b^{-1}\sigma^2$ (see~\eqref{eqn:asymptotic_variance} with $m=2$).
For example, with $b=n^{-1/4}$, asymptotically $n^{2/5}( \lambda_n^{SI}(x_0)-\lambda_0(x_0))$ behaves like a normal distribution with mean of the order $n^{-1/10}$ and variance $n^{1/20}\sigma^2$. Hence, the confidence interval becomes
\begin{equation}
\label{def:confint2}
\lambda_n^{SI}(x_0)
\pm n^{-3/8}
\widehat{\sigma}_n(x_0)q_{1-\alpha/2},
\end{equation}
where
\begin{equation}
\label{eqn:plugin_sigma}
\widehat{\sigma}_n(x_0)=\frac{\lambda_n^{SI}(x_0)}{c \Phi_n(x_0;\hat{\beta}_n)}\int_{-1}^1k(y)^2\,\mathrm{d} y.
\end{equation}
Note that undersmoothing leads to confidence intervals of asymptotic length $O_P(n^{-3/8})$, while the optimal ones would be of length $O_P(n^{-2/5})$.
In our simulations, the event times are generated from a Weibull baseline distribution with shape parameter $1.5$ and scale $1$.
The real valued covariate and the censoring time are chosen to be uniformly distributed on the interval $(0,1)$ and we take $\beta_0=0.5$.
We note that this setup corresponds to around $35\%$ uncensored observations.
Confidence intervals are calculated at the point $x_0=0.5$ using 10\,000 sets of data and
we take bandwidth $b=cn^{-1/4}$, with $c=1$, and kernel function $k(u)=(35/32)(1-u^2)^3\1_{\{|u|\leq 1\}}$.

{It is important to note that the performance depends strongly on the choice of the constant~$c$, 
because the asymptotic length is inversely proportional to $c$ (see~\eqref{eqn:plugin_sigma}). 
This means that, by choosing a smaller $c$ we get wider confidence intervals and as a result higher coverage probabilities.
However, it is not clear which would be the optimal choice of such a constant.
This is actually a common problem in the literature (see for example~\cite{CHT06} and \cite{GCM96}).
As indicated in~\cite{mullerwang1990}, cross-validation methods that consider a trade-off between bias and variance
suffer from the fact that the variance of the estimator increases as one approaches the endpoint of the support.
This is even enforced in our setting, because the bias is also decreasing when approaching the endpoint of the support.
We tried a locally adaptive choice of the bandwidth, as proposed in~\cite{MW90},
by minimizing an estimator of the Mean Squared Error, but in our setting this method did not lead to better results.
A simple choice is to take~$c$ equal to the range of the data (see~\cite{GJ15}), which in our case corresponds to $c=1$.}

\begin{table}[t]
\[
\begin{tabular}{ccccccccccccc}
\toprule
      &    \multicolumn{2}{c}{SG}&   \multicolumn{2}{c}{SMLE} &    \multicolumn{2}{c}{SG$_0$}&   \multicolumn{2}{c}{SMLE$_0$} & \multicolumn{2}{c}{Kernel} & \multicolumn{2}{c}{Grenander} \\
$n$     & AL    & CP & AL    & CP    & AL    & CP & AL    & CP  & AL    & CP  & AL    & CP  \\
50    &  1.411 & 0.732 & 1.583 & 0.751 & 1.281 & 0.915 & 1.426 & 0.944 & 1.458 & 0.727 & 0.980 & 0.440 \\
100   &  0.996 & 0.740 & 1.101 & 0.796 & 0.984 & 0.941 & 1.057 & 0.958 & 1.055 & 0.756 & 0.757 & 0.500 \\
500   &  0.545 & 0.824 & 0.563 & 0.857 & 0.538 & 0.949 & 0.559 & 0.977 & 0.560 & 0.822 & 0.449 & 0.615 \\
1000  &  0.421 & 0.852 & 0.430 & 0.883 & 0.419 & 0.957 & 0.430 & 0.979 & 0.429 & 0.845 & 0.359 & 0.657 \\
5000  &  0.232 & 0.910 & 0.234 & 0.916 & 0.232 & 0.969 & 0.234 & 0.981 & 0.234 & 0.884 & 0.215 & 0.764 \\
\bottomrule
\end{tabular}
\]
\caption{The average length (AL) and the coverage probabilities (CP) for $95\%$ pointwise confidence intervals of the baseline hazard rate at the point $x_0=0.5$ based on the asymptotic distribution.
SG and SMLE use $\hat{\beta}_n$, while SG$_0$ and SMLE$_0$ use $\beta_0$.}
\label{tab:1}
\end{table}
Table~\ref{tab:1} shows the performance of the estimators. 
The four columns corresponding to SG and SMLE list the average length (AL) and the coverage probabilities (CP)
of the confidence intervals given in~\eqref{def:confint2} for various sample sizes. 
Results {indicate that the SMLE behaves slightly better, but as the sample size increases its behavior becomes comparable to that of the SG estimator.} 
Even though the coverage probabilities are {below} the nominal level of $95\%$, smoothing leads to significantly more accurate results in comparison with the non-smoothed Grenander-type estimator given in the last two columns of Table~\ref{tab:1}. The confidence intervals for the Grenander-type estimator are constructed on the basis of Theorem~2 in \cite{LopuhaaNane2013}, i.e., they are of the form
{$\tilde{\lambda}_n(x_0)\pm n^{-1/3}\hat{C}_n(x_0)q_{1-\alpha/2}(\Z)$,}
where
\[
\hat{C}_n(x_0)=\left(\frac{4\tilde{\lambda}_n(x_0)\tilde{\lambda}'_n(x_0)}{\Phi_n(x_0;\hat{\beta}_n)} \right)^{1/3},
\]
$q_\alpha(\Z)$ is the $\alpha$-quantile of the distribution of 
$\Z=\argmin_{t\in\R}\{W(t)+t^2\}$, 
with $W$ a standard two-sided Brownian motion starting from zero.
In particular, $q_{0.975}(\Z)=0.998181$.
The {main} advantage of using the non-smoothed Grenander-type estimator is that {it does not involve the choice of a tuning parameter.
However,} the performance is not satisfactory, because we still need to estimate the derivative of $\lambda_0$, 
which is difficult if the estimator of $\lambda_0$ is a step function. 
Here we use the slope of the segment $[\tilde{\lambda}_n(T_{(i)},\tilde{\lambda}_n(T_{i+1})]$ on the interval $[T_{(i)},T_{(i+1)}]$ that contains $x_0$.

We also compare the performance of the {SG estimator and the SMLE} with that of the ordinary (non-monotone) kernel estimator 
\[
\lambda_n^s(x_0)=\int k_b(x_0-u)\,\mathrm{d}\Lambda_n(u),
\]
which is shown in {columns 10-11} of Table~\ref{tab:1}. We note that the kernel estimator coincides with the naive estimator 
{that approximates the isotonized smoothed Breslow estimator, see Section~4 in~\cite{LopuhaaMusta2016}.} 
In their proof of Theorem~4.3 it is shown that~$\lambda_n^s$ exhibits a limit distribution which coincides with the  one of the smooth estimators in Theorem~\ref{theo:distr}. 
Also the numerical results in Table~\ref{tab:1} confirm that the performance of the kernel estimator is comparable with that of the {smoothed} isotonic estimators. 
However, we notice that the latter ones have slightly  better coverage probabilities and shorter confidence intervals.

\begin{figure}[t]
\centering
\subfloat[][SG]
{\includegraphics[bb=189 44 770 625,width=0.45\textwidth,clip=]{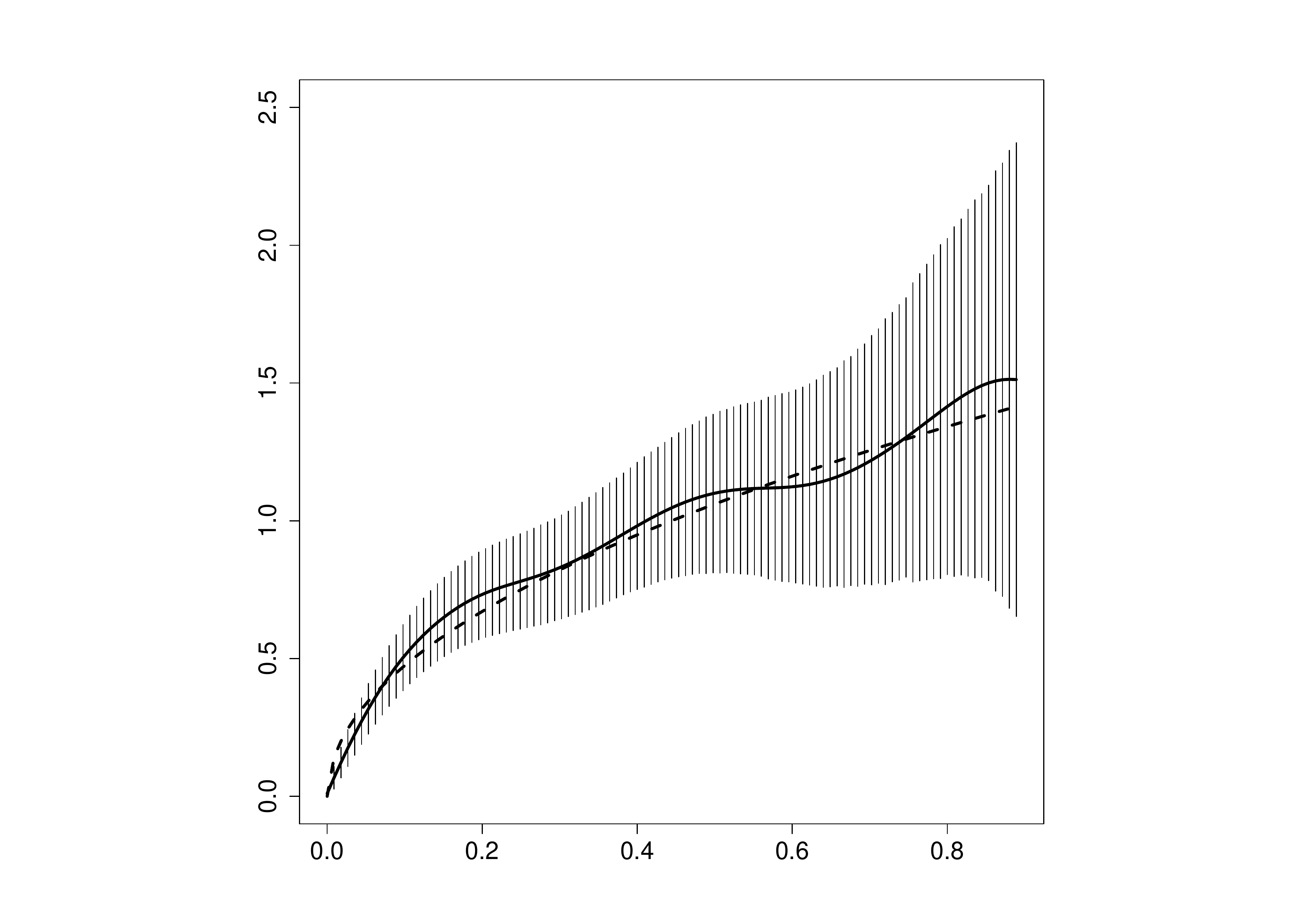}} \quad
\subfloat[][SMLE]
{\includegraphics[bb=189 44 770 625,width=0.45\textwidth,clip=]{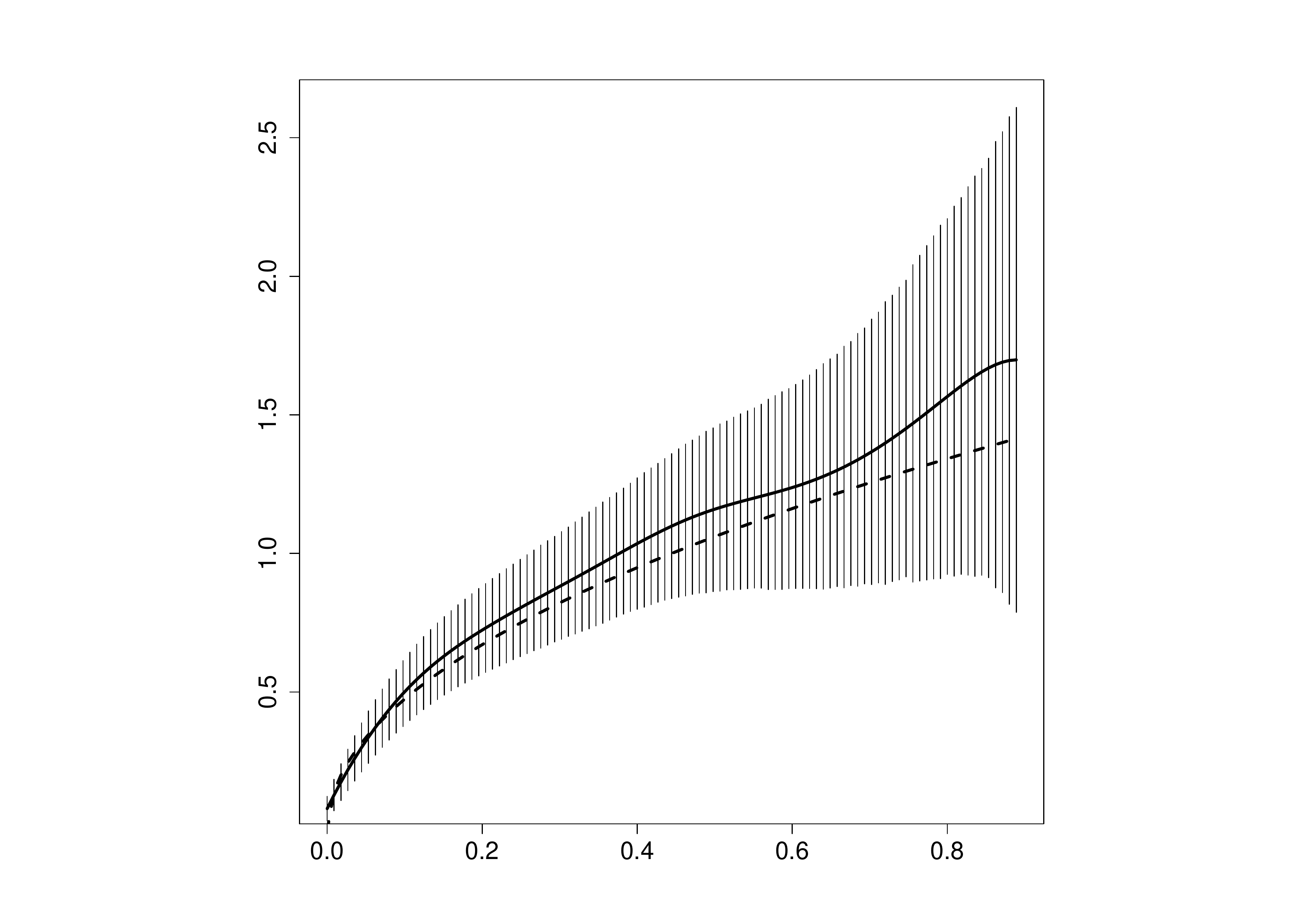}}
\caption{$95\%$ pointwise confidence intervals based on the asymptotic distribution for the baseline hazard rate using undersmoothing.}
\label{fig:subfig4}
\end{figure}
Moreover, as noticed in~\cite{LopuhaaMusta2016}, estimation of the parameter $\beta_0$ also has a great effect on the accuracy of the results. The columns 6-9 of Table~\ref{tab:1} show that if we use the true value of $\beta_0$ in the computation of the estimators, the coverage probabilities increase significantly but in this case the confidence intervals for the SMLE become too conservative.
Although the partial ML estimator $\hat{\beta}_n$ is a standard estimator for the regression coefficients,
the efficiency results are only asymptotic.
As pointed out in~\cite{CO84}  and ~\cite{RZ11}, for finite samples the use of the partial likelihood leads to a loss of accuracy.
Recently, \cite{RZ11} introduced
the~MLE for $\beta_0$ obtained by joint maximization of  the loglikelihood in the Cox model over both $\beta$ and $\lambda_0$.
It was shown that for small and moderate sample sizes, the joint MLE for $\beta_0$ performs better than~$\hat{\beta}_n$. However, in our case, using this estimator instead of $\hat{\beta}_n$, does not bring any essential difference in the coverage probabilities.

Pointwise confidence intervals, for a fixed sample size~$n=500$, at different points of the support are illustrated in Figure~\ref{fig:subfig4}.
The results are again comparable and the common feature is that the length increases as we move to the left boundary.
This is due to the fact that the length is proportional to the asymptotic standard deviation, which in this case turns out to be increasing, 
$\sigma^2(x)=1.5\sqrt{x}/(c\Phi(x;\beta_0)).$ Note that $\Phi(x;\beta_0)$ defined in~\eqref{eq:def Phi} is decreasing.

\subsection{Bootstrap confidence intervals}
\label{subsec:BSconf}
In an attempt to improve the coverage probabilities, we also construct bootstrap confidence intervals.
Studies on bootstrap confidence intervals in the Cox model are investigated in~\cite{Burr94} and~\cite{XSY2014}.
{In the latter paper, the authors investigate several bootstrap procedures for the Cox model.
We will use one (method M5) of the two proposals for a smooth bootstrap that had the best performance and were recommended by the authors.}

We fix the covariates and we generate the event time $X_i^*$ from a smooth estimate for the cdf of~$X$ conditional on $Z_i$:
\[
\hat{F}_n(x|Z_i)
=
1-\exp\left\{-\Lambda^s_n(x)\mathrm{e}^{\hat{\beta}'_nZ_i} \right\},
\]
where $\Lambda^s_n$ is the smoothed Breslow estimator 
\[
\Lambda^s_n(x)=\int k_b(x-u)\Lambda_n(u)\,\mathrm{d}u.
\]
The censoring times~$C^*_i$ are generated from the Kaplan-Meier estimate $\hat{G}_n$.
Then we take~$T^*_i=\min(X^*_i,C^*_i)$ and $\Delta^*_i=\1_{\{X^*_i\leq C^*_i\}}$.
For constructing the confidence intervals, we take $1000$ bootstrap samples $(T^*_i,\Delta^*_i,Z_i)$  and for each
bootstrap sample we compute the smoothed Grenander-type estimate~$\tilde{\lambda}^{SG,*}_n(x_0)$ 
and the smoothed maximum likelihood estimate $\hat{\lambda}^{SM,*}_n(x_0)$.
Here the kernel function is the same as before and the bandwidth is taken to be $b=n^{-1/5}$.
Then, the $100(1-\alpha)\%$ confidence interval is given by
\begin{equation}
\label{eq:CI bootstrap}
\left[
q^*_{\alpha/2}(x_0),q^*_{1-\alpha/2}(x_0)
\right],
\end{equation}
where $q^*_{\alpha}(x_0)$ is the $\alpha$-percentile of the $1000$ values of the estimates $\tilde{\lambda}^{SG,*}_n(x_0)$ or $\hat{\lambda}^{SM,*}_n(x_0)$.

The average length and the empirical coverage for $1000$ iterations and different sample sizes are reported in
Table~\ref{tab:bootstrap1}.
\begin{table}[h]
\begin{tabular}{ccccccc}
		\toprule
		&&    \multicolumn{2}{c}{SMLE}&&   \multicolumn{2}{c}{SG} \\
		\cline{3-4} \cline{6-7}
		\\[-5pt]
		$n$     && AL    & CP   && AL &CP \\
		100   && 1.870 & 0.948  && 1.376 & 0.899 \\
		500   && 0.730 & 0.942 && 0.660 & 0.892 \\
		1000  &&0.521& 0.960  && 0.487 & 0.902 \\
		5000 && 0.247 & 0.957 && 0.239 & 0.938 \\
		\bottomrule\\
\end{tabular}
	\caption{The average length (AL) and the coverage probabilities (CP) for the $95\%$ bootstrap confidence intervals
		of the baseline hazard rate at the point $x_0=0.5$, {using the tri-weight kernel and $b=n^{-2/5}$}.}
	\label{tab:bootstrap1}
\end{table}
We observe that bootstrap confidence intervals behave better {than} confidence intervals {in Table~\ref{tab:1}},
i.e., the coverage probabilities are closer to the nominal level of $95\%$.
Comparing also with the two alternative estimators considered in \cite{LopuhaaMusta2016} we notice that the SMLE and the 
maximum smoothed likelihood estimator (MSLE) have better coverage probabilities
{than the smoothed Grenander-type and isotonized Breslow estimator, respectively.

In order to provide some theoretical evidence for the consistency of the method,
we would like to establish that,
given the data $(T_1,\Delta_1,Z_1),\dots,(T_n,\Delta_n,Z_n)$, it holds
\begin{equation}
\label{eqn:asymptotic_normality_bootstrap}
n^{2/5}\left(\lambda_n^{SI,*}(x)-\lambda_n^{SI}(x)\right)\xrightarrow{d} N(\tilde{\mu},\sigma^2),
\end{equation}
for some $\tilde{\mu}\in\R$ (possibly different from $\mu$ in Theorem~\ref{theo:distr}) and $\sigma^2$ as in~\eqref{eqn:as.mean-var},
where $\lambda_n^{SI}$ is one of the smooth isotonic estimators at hand and $\lambda_n^{SI,*}$ is the same estimator computed for the bootstrap sample. 
A detailed investigation of~\eqref{eqn:asymptotic_normality_bootstrap} is beyond the scope of this paper.
Nevertheless, in view of Remark~\ref{rem:partial MLE}, we are able to obtain~\eqref{eqn:asymptotic_normality_bootstrap} for the smoothed Grenander estimator,
if $\hat\beta_n^*-\hat{\beta}_n\to 0$, 
for almost all sequences $(T_i^*,\Delta_i^*,Z_i)$, $i=1,2,\ldots$, conditional on the sequence  $(T_i,\Delta_i,Z_i)$, $i=1,2,\ldots$,
and $\sqrt{n}(\hat\beta_n^*-\hat{\beta}_n)=O_p^*(1)$.
By the latter we mean that for all $\epsilon>0$, there exists $M>0$ such that
\[
\limsup_{n\to\infty}
P_n^*
\left(
\sqrt{n}
|\hat\beta_n^*-\hat{\beta}_n|
>M
\right)
<
\epsilon,
\qquad
\p-\text{almost surely}.
\]
where $P_n^*$ is the measure corresponding to the distribution of
$(T^*,\Delta^*,Z)$ conditional on the data $(T_1,\Delta_1,Z_1)$, $\ldots$ ,$(T_n,\Delta_n,Z_n)$, with
$T^*=(\min(X^*,C^*)$ and $\Delta^*=\1_{\{X^*\leq C^*\}},Z)$,
where~$X^*$ conditional on $Z$ has distribution function $\hat F_n(x\mid Z)$ and $C^*$ has distribution function~$\hat{G}_n$.
To prove~\eqref{eqn:asymptotic_normality_bootstrap}, we mimic the proof of Theorem~\ref{theo:distr}, 
which means that one needs to establish the bootstrap versions of Lemmas~\ref{le:2-1}-\ref{le:result5}.
A brief sketch of the arguments is provided in Appendix C of~~\cite{LopuhaaMustaSM2016}.

Then, we can approximate the distribution of $n^{2/5}(\lambda_0(x_0)-\lambda_n^{SI}(x_0))$ 
by the distribution of $n^{2/5}(\lambda_n^{SI,*}(x_0)-\lambda_n^{SI}(x_0))-(\tilde{\mu}+\mu)$. 
Consequently, we can write
\[
\begin{split}
&
P_n^*\{q_{\alpha/2}^*(x_0)\leq \lambda_n^{SI,*}(x)\leq q_{1-\alpha/2}^*(x_0)\}\\
&=
P^*_n\left\{\lambda_0(x_0)\in\left[q_{\alpha/2}^*(x_0)-n^{-2/5}(\tilde{\mu}+\mu),q_{1-\alpha/2}^*(x_0)-n^{-2/5}(\tilde{\mu}+\mu)\right] \right\}
\end{split}
\]
This means that we should actually take $[q_{\alpha/2}^*(x_0),q_{1-\alpha/2}^*(x_0)]-n^{-2/5}(\tilde{\mu}+\mu)$ instead of~\eqref{eq:CI bootstrap}.
The use of~\eqref{eq:CI bootstrap} avoids bias estimation.
However, since the effect of the bias is of the order~$n^{-2/5}$, the results are still satisfactory. 
In} order to further reduce the effect of the bias, we also investigated the possibility of constructing bootstrap confidence intervals with undersmoothing, i.e, 
we repeat the previous procedure with bandwidth $b=n^{-1/4}$. 
Results are shown in Table~\ref{tab:bootstrap2}. 
\begin{table}[h]
\begin{tabular}{ccccccc}
	\toprule
	&&    \multicolumn{2}{c}{SMLE}&&   \multicolumn{2}{c}{SG} \\
	\cline{3-4} \cline{6-7}
	\\[-5pt]
$n$     && AL    & CP   && AL &CP  \\
100   && 1.901 & 0.954  && 1.415 & 0.900 \\
500   && 0.749 & 0.951 && 0.672 & 0.918 \\
1000  &&0.540& 0.950  && 0.501 & 0.924 \\
{5000} && 0.262 & 0.965 && 0.252 & 0.952\\
\bottomrule\\
\end{tabular}
\caption{The average length (AL) and the coverage probabilities (CP) for the $95\%$ bootstrap confidence intervals
of the baseline hazard rate at the point $x_0=0.5$, {using the tri-weight kernel and $b=n^{-1/4}$.}}
\label{tab:bootstrap2}
\end{table}
We notice that the length of the confidence interval increases slightly and the coverage probabilities 
{improve significantly.
To summarize,} also the bootstrap confidence intervals are affected by the choice of the bandwidth, 
but the results are more satisfactory in comparison with the ones in Table~\ref{tab:1}.

\section{Discussion}
\label{sec:discussion}
In this paper we considered smooth estimation under monotonicity constraints of the baseline hazard rate in the Cox model. 
We {investigated} the asymptotic behavior of two estimators, which are the kernel smoothed versions of the monotone MLE and Grenander-type estimator. 
The main result is that they are asymptotically equivalent with a normal limit distribution at rate $n^{-m/(2m+1)}$, 
where $m$ is the degree of smoothness assumed for the baseline hazard. 
{Two other methods that} combine smoothing and isotonization for estimation of the baseline hazard in the Cox model 
were considered in~\cite{LopuhaaMusta2016}. 
As shown in Theorems~3.6 and~4.4 in~\cite{LopuhaaMusta2016}, the smoothed Grenander-type estimator, 
the smoothed maximum likelihood estimator, and the isotonized kernel estimator, are {all} asymptotically equivalent, 
while the maximum smoothed likelihood estimator exhibits a different asymptotic bias (which might be smaller or larger than the one of the previous three estimators). 
This means that, from the theoretical point of view, there is no reason to prefer one estimator with respect to the other (apart from the fact that the kernel smoothed estimators are differentiable while the other two are usually only continuous). 

{The method used to establish asymptotic normality for the estimators in this paper is quite different from the ones in~\cite{LopuhaaMusta2016}.
In the latter paper, the isotonization step was performed after a smoothing step.
As a consequence, the resulting estimators are asymptotically equivalent to corresponding naive estimators that 
are combinations of ordinary kernel type estimators, to which standard techniques apply.
This approach does not apply to the smoothed isotonic estimators in this paper.
Alternatively, we followed the approach from~\cite{GJ14} based on $L_2$-bounds for the isotonic estimator.
The approach had to be adapted at several points leading to $L_2$-bounds that are suboptimal, 
but sufficient for our purposes.}

Furthermore, we investigated also the finite sample performance of these estimators by constructing pointwise confidence intervals. First, making use of the theoretical results, we construct pointwise confidence intervals based on the limit distribution with undersmoothing to avoid bias estimation. Results confirm the comparable behavior of the four methods and favor the use of the smoothed isotonic estimators instead of the unsmoothed Grenander-type estimator or the non-isotonic kernel estimator. However, coverage probabilities are far from the nominal level and strongly depend on the choice of the bandwidth and the accuracy in the estimation of the regression coefficient $\beta_0$. Since most of the well-known methods to overcome these problems {do not seem} to work in our setting, a thorough investigation is still needed for improving the performance of the confidence intervals. 
Instead, we choose to exploit pointwise confidence intervals based on bootstrap procedures. 
{As it turns out,} the simple percentile bootstrap works better than the studentized one. 
{Such} a phenomenon was also observed in \cite{Burr94}. 
The four estimators, 
the smoothed  maximum likelihood estimator (SMLE), 
the smoothed Grenander-type estimator, 
the maximum smoothed likelihood estimator (MSLE) and 
the isotonized smoothed Breslow estimator, again exhibit comparable behavior but the SMLE and the MSLE have slightly better coverage probabilities. 
The performance is satisfactory, but still further investigation is required for bandwidth selection and correcting the asymptotic bias, which might improve the results.

\section{Proofs}
\label{sec:proofs}
\begin{proof}[Proof of Lemma~\ref{le:2-1}]
Define
$D_n^{(1)}(x;\beta)=\partial \Phi_n(x;\beta)/\partial \beta$ 
and let $D_{nj}^{(1)}(x;\beta)$ be the $j$th component of $D_n^{(1)}(x;\beta)$, for $j=1,\dots,p$. 
Then according to the proof of Lemma~3(iv) in~\cite{LopuhaaNane2013},  for any sequence~$\beta^*_n$, such that $\beta^*_n\to \beta_0$ almost surely,
it holds
\[
\limsup_{n\to\infty}\sup_{x\in\R}|D_n^{(1)}(x;\beta^*_n)|<\infty.
\]
In fact, from its proof, it can be seen that
\[
\sup_{x\in\R}|D_{nj}^{(1)}(x;\beta^*_n)|\leq \sum_{I_k\subseteq I}\left[\frac{1}{n}\sum_{i=1}^n|Z_i|\,\mathrm{e}^{\gamma'_kZ_i} \right]\to \sum_{I_k\subseteq I}\E\left[|Z|\mathrm{e}^{\gamma'_kZ} \right]<2^p\sup_{|\beta-\beta_0|\leq\epsilon}\E\left[|Z|\mathrm{e}^{\beta'Z}\right]<\infty
\]
with probability $1$, where the summations are over all subsets $I_k=\{i_1,\ldots,i_k\}$ of $I=\{1,\ldots,p\}$,
and $\gamma_k$ is the vector consisting of coordinates
$\gamma_{kj}=\beta_{0j}+\epsilon/(2\sqrt{p})$, for $j\in I_k$, and
$\gamma_{kj}=\beta_{0j}-\epsilon/(2\sqrt{p})$, for $j\in I\setminus I_k$.
Therefore,
\[
\sup_{x\in\R}|D_{n}^{(1)}(x;\beta^*_n)|
\leq
\sqrt{p}\sum_{I_k\subseteq I}\left(\frac{1}{n}\sum_{i=1}^n|Z_i|\,\mathrm{e}^{\gamma'_kZ_i} \right)
\to
\sqrt{p}\sum_{I_k\subseteq I}\E\left[|Z|\mathrm{e}^{\gamma'_kZ} \right]
\]
with probability one.
Hence, if for some $\xi_1>0$,
\begin{equation}
\label{eqn:E1}
E_{n,1}=\left\{\left|\sqrt{p}\sum_{I_k\subseteq I}
\left(\frac{1}{n}\sum_{i=1}^n|Z_i|\,\mathrm{e}^{\gamma'_kZ_i} \right)
-
\sqrt{p}\sum_{I_k\subseteq I}\E\left[|Z|\mathrm{e}^{\gamma'_kZ} \right]\right|\leq\xi_1 \right\},
\end{equation}
then $\1_{E_{n,1}}\to 1$ in probability.
Moreover, on this event, we have
\begin{equation}
\label{eqn:dn2}
\sup_{x\in\R}|D_{n}^{(1)}(x;\beta^*_n)|\leq \sqrt{p}\sum_{I_k\subseteq I}\E\left[|Z|\mathrm{e}^{\gamma'_kZ} \right]+\xi_1,
\end{equation}
i.e., $\sup_{x\in\R}|D_{n}^{(1)}(x;\beta^*_n)|$ is bounded uniformly in $n$.
For $\xi_2,\xi_3,\xi_4>0$ and $0<M<\tau_H$ define
\begin{equation}
\label{eqn:ef}
\begin{split}
E_{n,2}
&=
\left\{n^{2/3}|\hat{\beta}_n-\beta_0|^2<\xi_2\right\},\qquad
\displaystyle E_{n,3}=\left\{\sup_{x\in[0,M]}|\tilde{\Lambda}_n(x)-\Lambda_0(x)|<\xi_3 \right\},\\
\displaystyle E_{n,4}
&=
\left\{n^{1/3}\sup_{x\in \mathbb{R}}\left|\Phi_n(x;\beta_0)-\Phi(x;\beta_0)\right|\leq \xi_4\right\},
\qquad
E_{n,5}
=
\left\{T_{(n)}>M\right\}
\end{split}
\end{equation}
where $T_{(n)}$ denotes the last observed time.
Because $\sqrt{n}(\hat\beta_n-\beta_0)=O_p(1)$ (see Theorem~3.2 in~\cite{Tsiatis81}), together with~\eqref{eqn:Breslow}
and Lemma~4 in~\cite{LopuhaaNane2013}, it follows that
$\1_{E_n}\to 1$ in probability, for the event $E_n=E_{n,1}\cap E_{n,2}\cap E_{n,3}\cap E_{n,4}\cap E_{n,5}$.

From the definitions of $a_{n,x}$, $\theta_{n,x}$ and $H^{uc}$,
in~\eqref{eqn:an}, \eqref{eqn:thetan}, and~\eqref{eq:def Huc}, respectively,
we have
\[
\int
\theta_{n,x}(u,\delta,z)\,\mathrm{d}\p(u,\delta,z)
=
\1_{E_n}\left\{
\int a_{n,x}(u)\,\mathrm{d}H^{uc}(u)-\int \mathrm{e}^{\hat{\beta}'_n\,z}\,\int_{v=0}^u a_{n,x}(v)\,\mathrm{d}\tilde{\Lambda}_n(v)\,\mathrm{d}\p(u,\delta,z)
\right\}.
\]
Then, by applying Fubini's theorem, together with~\eqref{eqn:lambda0}, we obtain
\[
\begin{split}
&\int\theta_{n,x}(u,\delta,z)\,\mathrm{d}\p(u,\delta,z)\\
&=
\1_{E_n}\left\{\int a_{n,x}(u)\,\mathrm{d}H^{uc}(u)-\int a_{n,x}(v)\int_{u=v}^\infty \mathrm{e}^{\hat{\beta}'_n\,z}\,\mathrm{d}\p(u,\delta,z)\,\mathrm{d}\tilde{\Lambda}_n(v)\right\}\\
&=
\1_{E_n}\left\{\int a_{n,x}(u)\,\mathrm{d}H^{uc}(u)-\int a_{n,x}(v)\,\Phi(v;\hat{\beta}_n)\,\mathrm{d}\tilde{\Lambda}_n(v)\right\}\\
&=
\1_{E_n}\left\{\int\frac{k_b(x-u)}{\Phi(u;\beta_0)}\,\mathrm{d}H^{uc}(u)-\int k_b(x-u)\,\frac{\Phi(u;\hat{\beta}_n)}{\Phi(u;\beta_0)}\,\mathrm{d}\tilde{\Lambda}_n(u)\right\}\\
&=
\1_{E_n}\left\{-\int k_b(x-u)\,\mathrm{d}(\tilde{\Lambda}_n-\Lambda_0)(u)+\int k_b(x-u)\,\left(1-\frac{\Phi(u;\hat{\beta}_n)}{\Phi(u;\beta_0)}\right)\,\mathrm{d}\tilde{\Lambda}_n(u)\right\}.
\end{split}
\]
The mean value theorem yields
\[
\begin{split}
\int k_b(x-u)\,\left|1-\frac{\Phi(u;\hat{\beta}_n)}{\Phi(u;\beta_0)}\right|\,\mathrm{d}\tilde{\Lambda}_n(u)
&=
\int k_b(x-u)\,\frac{|\Phi(u;\beta_0)-\Phi(u;\hat{\beta}_n)|}{\Phi(u;\beta_0)}\,\mathrm{d}\tilde{\Lambda}_n(u)\\
&\leq
|\hat{\beta}_n-\beta_0|
\sup_{y\in\R}\left|\frac{\partial\Phi(y;\beta^*)}{\partial \beta}\right|\,\frac{\tilde{\lambda}_n^{SG}(x)}{\Phi(x+b;\beta_0)},
\end{split}
\]
with $|\beta^*-\beta_0|\leq|\hat{\beta}_n-\beta_0|$.
According to Lemma~3(iii) in~\cite{LopuhaaNane2013}, for $\epsilon>0$ from (A2),
\[
\sup_{y\in\R}\left|\frac{\partial\Phi(y;\beta^*)}{\partial \beta}\right|
<
\sup_{y\in\R}\sup_{|\beta-\beta_0|<\epsilon}
\left|\frac{\partial\Phi(y;\beta)}{\partial \beta}\right|
<\infty.
\]
Furthermore, there exists $M<\tau_H$, such that for sufficiently large $n$ we have $x+b\leq M$.
This yields the following bound $\Phi(x+b;\beta_0)\geq\Phi(M;\beta_0)>0$.
Moreover, according to~\eqref{eq:unif conv SG}, $\tilde{\lambda}_n^{SG}(x)\to\lambda_0(x)$ with probability one.
Since $|\hat{\beta}_n-\beta_0|=O_p(n^{-1/2})$ (see Theorem~3.1 in~\cite{Tsiatis81}),
it follows that
\[
\1_{E_n}\int k_b(x-u)\,\left|1-\frac{\Phi(u;\hat{\beta}_n)}{\Phi(u;\beta_0)}\right|\,\mathrm{d}\tilde{\Lambda}_n(u)=O_p(n^{-1/2}),
\]
which finishes the proof.
\end{proof}

\begin{proof}[Proof of Lemma~\ref{le:2-2}]
By means of Fubini's theorem
\[
\begin{split}
&
\int\overline{\theta}_{n,x}(u,\delta,z)\,\mathrm{d}\p_n(u,\delta,z)\\
&=
\int \delta\,\overline{\Psi}_{n,x}(u)\,\mathrm{d}\p_n(u,\delta,z)
-
\int \mathrm{e}^{\hat{\beta}'_n\,z}\int_{v=0}^u \overline{\Psi}_{n,x}(v)\,\mathrm{d}\tilde{\Lambda}_n(v)\,\mathrm{d}\p_n(u,\delta,z)\\
&=
\int \delta\,\overline{\Psi}_{n,x}(u)\,\mathrm{d}\p_n(u,\delta,z)
-
\int \overline{\Psi}_{n,x}(v)\int \1_{\{u\geq v\}} \mathrm{e}^{\hat{\beta}'_n\,z}\,\mathrm{d}\p_n(u,\delta,z) \,\mathrm{d}\tilde{\Lambda}_n(v)\\
&=
\1_{E_n}\left\{
\int \1_{[x-b,x+b]}(u) \delta\,\frac{\overline{a}_{n,x}\overline{\Phi}_n(u;\hat{\beta}_n)}{\Phi_n(u;\hat{\beta}_n)}\,\mathrm{d}\p_n(u,\delta,z)-\int_{x-b}^{x+b}\overline{a}_{n,x}\overline{\Phi}_n(v;\hat{\beta}_n) \,\mathrm{d}\tilde{\Lambda}_n(v)\right\}\\
&=
\1_{E_n}
\sum_{i=0}^{m}
\overline{a}_{n,x}\overline{\Phi}_n(\tau_{i+1};\hat{\beta}_n)
\left\{
\int \frac{\1_{(\tau_i,\tau_{i+1}]}(u)\delta}{\Phi_n(u;\hat{\beta}_n)}\,\mathrm{d}\p_n(u,\delta,z)
-
\bigg(\tilde{\Lambda}_n(\tau_{i+1})-\tilde{\Lambda}_n(\tau_i)\bigg)
\right\}
\end{split}
\]
and~\eqref{eqn:int} follows from the characterization of the Breslow estimator in~\eqref{eq:Breslow}.
\end{proof}
To obtain suitable bounds for~\eqref{eq:E sup}, we will establish bounds on the tail probabilities of~$\tilde U_n(a)$ defined in~\eqref{eq:def Un}.
To this end we consider a suitable martingale that will approximate the process $\Lambda_n-\Lambda_0$.
For $i=1,2,\ldots,n$, let $N_i(t)=\1_{\{X_i\leq t\}}\Delta_i$ be the right continuous counting process for the number of observed failures
on $(0,t]$ and $Y_i(t)=\1_{\{T_i\geq t\}}$ be the at-risk process.
Then, for each $i=1,2,\ldots,n$,
$M_i(t)=N_i(t)-A_i(t)$,
with $A_i(t)=\int_0^t Y_i(s)\mathrm{e}^{\beta'_0Z_i}\,\mathrm{d}\Lambda_0(s)$,
is a mean zero martingale with respect to the filtration
\[
\F^n_t=
\sigma
\left\{
\1_{\{X_i\leq s\}}\Delta_i,\1_{\{T_i\geq s\}}, Z_i:1\leq i\leq n,0\leq s\leq t
\right\}.
\]
(e.g., see~\cite{kalbfleischprentice2002}).
Furthermore, it is square integrable, since
\[
\E\left[M_i(t)^2\right]
\leq
2+2\int_0^t\E\left[\1_{\{T_i\geq s\}}\mathrm{e}^{2\beta'_0Z_i}\right]\lambda^2_0(s)\,\mathrm{d}s
\leq
2+2\tau_H\,\lambda^2_0(\tau_H)\Phi(0;2\beta_0)<\infty.
\]
Finally, it has predictable variation process
$\langle M_i\rangle=A_i(t)$ (e.g., see~\cite{gill84} or Theorem~2 of Appendix B in~\cite{SW86}).
For each $n\geq 1$, define
\begin{equation}
\label{eq:def Mn}
\mathbb{N}_n(t)
=
\sum_{i=1}^n N_i(t),
\qquad
\mathbb{A}_n(t)
=
\sum_{i=1}^n A_i(t),
\qquad
\mathbb{M}_n(t)
=\mathbb{N}_n(t)-\mathbb{A}_n(t).
\end{equation}
Then $\m_n(t)$ is a mean zero square integrable martingale with predictable variation process
\[
\langle\m_n\rangle(t)
=
\sum_{i=1}^n
\langle M_i\rangle(t)
=
\sum_{i=1}^n\int_0^t\1_{\{T_i\geq s\}}\,\mathrm{e}^{\beta'_0Z_i}\,\mathrm{d}\Lambda_0(s)
=
\int_0^tn\Phi_n(s;\beta_0)\,\mathrm{d}\Lambda_0(s),
\]
where $\Phi_n$ is defined in~\eqref{eq:def Phin}.

\begin{lemma}
\label{lem:mart2}
Suppose that (A1)--(A2) hold.
Let $0<M<\tau_H$ and let $\Phi$ be defined in~\eqref{eq:def Phi}.
Then, the process
\begin{equation}
\label{eqn:b_n}
\w_n(t)=\int_0^{t\wedge M}\frac{1}{n\Phi(s;\beta_0)}\,\mathrm{d}\m_n(s)
\end{equation}
is a mean zero, square integrable martingale with respect to the filtration $\F^n_t$,
Moreover, $\w_n$ has  predictable variation process
\[
\langle\w_n\rangle(t)
=
\int_0^{t\wedge M}
\frac{\lambda_0(s)\Phi_n(s;\beta_0)}{n\Phi^2(s;\beta_0)}\,\mathrm{d}s.
\]
\end{lemma}
\begin{proof}
Write
\[
\w_n(t)
=
\int_0^t Y_n(s)\,\mathrm{d}\m_n(s),
\quad
\text{where }
Y_n(s)=\frac{\1_{\{s\leq M\}}}{n\Phi(s;\beta_0)},
\]
and $\m_n=\mathbb{N}_n-\mathbb{A}_n$.
We apply Theorem B.3.1c in~\cite{SW86} with $Y$, $H$, $M$, $N$, and $A$,
replaced by $\w_n$, $Y_n$, $\m_n$, $\mathbb{N}_n$, and $\mathbb{A}_n$, respectively.
In order to check the conditions of this theorem,
note that $Y_n$ is a predictable process satisfying $|Y_n(t)|<\infty$, almost surely,
for all~$t\geq 0$, and that
\[
\begin{split}
\int_0^tY_n(s)\,\mathrm{d}A_n(s)
&=
\sum_{i=1}^n\int_0^t
\frac{\1_{\{s\leq M\}}}{n\Phi(s;\beta_0)}
\1_{\{T_i\geq s\}}\mathrm{e}^{\beta'_0Z_i}\,\mathrm{d}\Lambda_0(s)\\
&=
\int_0^t\frac{\1_{\{s\leq M\}}}{\Phi(s;\beta_0)}\Phi_n(s;\beta_0)\,\mathrm{d}\Lambda_0(s)<\infty,
\quad\text{a.s.}.
\end{split}
\]
Moreover,
since for $s\leq M$ we have $\Phi(s;\beta_0)\geq \Phi(M;\beta_0)>0$, it follows that
\[
\begin{split}
\E\left[\int_0^\infty Y^2_n(s)\,\mathrm{d}\langle\m_n\rangle(s)\right]
&=
\E\left[\int_0^\infty
\frac{\1_{\{s\leq M\}}}{n\Phi^2(s;\beta_0)}\Phi_n(s;\beta_0)\,\mathrm{d}\Lambda_0(s)\right]\\
&\leq
\frac{\lambda_0(\tau_H)M}{n^2\Phi^2(M;\beta_0)}
\sum_{i=1}^n \E\left[ \mathrm{e}^{\beta'_0Z_i}\right]<\infty,
\end{split}
\]
because of the assumption~(A2).
It follows from Theorem B.3.1c in~\cite{SW86}, that $\w_n$ is a square integrable martingale with mean zero and predictable variation process
\[
\langle\w_n\rangle(t)
=
\int_0^tY_n^2(s)\,\mathrm{d}\langle\m_n\rangle(s)
=
\int_0^t \frac{\1_{\{s\leq M\}}}{n\Phi^2(s;\beta_0)}\Phi_n(s;\beta_0)\,\mathrm{d}\Lambda_0(s),
\]
where $\Phi$ and $\Phi_n$ are defined in~\eqref{eq:def Phi} and~\eqref{eq:def Phin}, respectively.
\end{proof}
It is straightforward to verify that for $t\in[0,M]$ and $M<T_{(n)}$,
\begin{equation}
\label{eq:martingale relation}
\Lambda_n(t)-\Lambda_0(t)
=
\w_n(t)+R_n(t),
\end{equation}
where
\begin{equation}
\label{eqn:r}
R_n(t)
=
\int_0^t\frac{\Phi_n(s;\beta_0)}{\Phi(s;\beta_0)}\,\mathrm{d}\Lambda_0(s)-\Lambda_0(t)
+
\int_0^t
\left(
\frac{1}{\Phi_n(s;\hat\beta_n)}-\frac{1}{\Phi(s;\beta_0)}
\right)
\,\mathrm{d}H_n^{uc}(s),
\end{equation}
with
\begin{equation}
\label{def:Hnuc}
H_n^{uc}(x)
=
\int \delta \1_{\{t\leq x\}}\,\mathrm{d}\p_n(t,\delta,z).
\end{equation}
For establishing suitable bounds on the tail probabilities of $\tilde U_n(a)$, we need the following result
for the process $\w_n$, which is comparable to condition (A2) in~\cite{durot2007}.
\begin{lemma}
\label{le:mart}
Suppose that (A1)--(A2) hold.
Let $0<M<\tau_H$ and let $\w_n$ be defined as in~\eqref{eqn:b_n}.
Then, there exists a constant $C>0$ such that, for all $x>0$ and $t\in[0,M]$,
\[
\E\left[\sup_{u\in[0,M],|t-u|\leq x}
\big(\w_n(u)-\w_n(t)\big)^2\right]\leq \frac{C\,x}{n}.
\]
\end{lemma}
\begin{proof}
The proof is similar to that of Theorem~3 in~\cite{durot2007}.
First consider the case $t\leq u\leq t+x$.
According to Lemma~\ref{lem:mart2}, $\w_n$ is a martingale.
Hence,
by Doob's inequality, we have
\begin{equation}
\label{eqn:eq1}
\begin{split}
\E\left[\sup_{u\in[0,M],\,t\leq u\leq t+x}(\w_n(u)-\w_n(t))^2\right]
&\leq
4\E\left[\big(\w_n\big((t+x)\wedge M\big)-\w_n(t)\big)^2\right]\\
&=
4\E\left[\w_n\big((t+x)\wedge M\big)^2-\w_n(t)^2\right]\\
&=
4\E\left[\int_{t}^{(t+x)\wedge M}
\frac{\Phi_n(s;\beta_0)\lambda_0(s)}{n\Phi^2(s;\beta_0)}\,\mathrm{d}s \right]\\
&\leq
\frac{4\lambda(M)x}{n\Phi^2(M;\beta_0)}
\E\left[\Phi_n(0;\beta_0)\right],
\end{split}
\end{equation}
where according to~(A2),
\[
\E\left[\Phi_n(0;\beta_0)\right]
=
\frac{1}{n}\sum_{i=1}^n\E\left[\mathrm{e}^{\beta'_0\,Z_i}\right]\leq C,
\]
for some $C>0$.
This proves the lemma for the case $t\leq u\leq t+x$.

For the case $t-x\leq u\leq t$, we can write
\[
\begin{split}
&
\E\left[\sup_{u\in[0,M],t-x\leq u\leq t}
\big(\w_n(u)-\w_n(t)\big)^2\right]
=
\E\left[\sup_{0\vee(t-x)\leq u\leq t}
\big(\w_n(u)-\w_n(t)\big)^2\right]\\
&\quad\leq
2\E\left[\big(\w_n(t)-\w_n(0\vee(t-x))\big)^2\right]
+
2\E\left[\sup_{0\vee(t-x)\leq u<t}
\big(\w_n(u)-\w_n(0\vee(t-x))\big)^2\right].
\end{split}
\]
Then similar to~\eqref{eqn:eq1}, the right hand side is bounded by
\[
\begin{split}
&
2\E\left[\big(\w_n(t)-\w_n(0\vee(t-x))\big)^2\right]
+
8\E\left[\big(\w_n(t)-\w_n(0\vee(t-x))\big)^2\right]\\
&=
10\E\left[\w_n(t)^2-\w_n(0\vee(t-x))^2\right]
=
10\E\left[\int_{0\vee(t-x)}^{t}\frac{\Phi_n(s;\beta_0)\,\lambda_0(s)}{n\Phi^2(s;\beta_0)}\,\mathrm{d}s  \right]\\
&\leq
\frac{10\,\lambda(M)\,x}{n\Phi^2(M;\beta_0)}\,\E\left[\Phi_n(0;\beta_0) \right]
\leq
\frac{Cx}{n},
\end{split}
\]
for some $C>0$.
This concludes the proof.
\end{proof}

In what follows, let $0<M<\tau_H$.
Moreover, let $U$ be the inverse of $\lambda_0$ on $[\lambda_0(0),\lambda_0(M)]$, i.e.,
\begin{equation}
\label{def:U SG}
U(a)
=
\begin{cases}
0 & a< \lambda_0(0);\\
\lambda_0^{-1}(a) & a\in[\lambda_0(0),\lambda_0(M)];\\
M & a>\lambda_0(M).
\end{cases}
\end{equation}
Note that $U$ is continuous and differentiable on $(\lambda_0(0),\lambda_0(M)),$ but it is different from the inverse of $\lambda_0$ on the entire interval $[\lambda_0(0),\lambda_0(\tau_H)]$.
\begin{lemma}
\label{le:inv}
Suppose that (A1)--(A2) hold.
Let $0<M<\tau_H$ and let $\tilde U_n$ and $U$ be defined in~\eqref{eq:def Un} and~\eqref{def:U SG}, respectively.
Suppose that $H^{uc}$, defined in~\eqref{eq:def Huc}, has a bounded derivative $h^{uc}$ on $[0,M]$
and that $\lambda_0'$ is bounded below by a strictly positive constant.
Then, there exists an event~$E_n$, such that $\1_{E_n}\to 1$ in probability, and
a constant $K$ such that, for every $a\geq 0$ and $x> 0$,
\begin{equation}
\label{eqn:inv}
\p\left(
\left\{|\tilde{U}_n(a)-U(a)|\geq x\right\}
\cap E_n\cap
\left\{\tilde{U}_n(a)\leq M\right\}
\right)
\leq
K
\max\left\{\frac{1}{nx^3},\frac{1}{n^3x^5}\right\},
\end{equation}
for $n$ sufficiently large.
\end{lemma}
Note that Lemma~\ref{le:mart} and Lemma~\ref{le:inv} correspond to Theorem~3(i) and Lemma~2 in~\cite{durot2007}.
It is useful to spend some words on the restriction to the event $E_n\cap\{\tilde{U}_n(a)\leq M\}$.
The event $\{\tilde{U}_n(a)\leq M\}$ is implicit in~\cite{durot2007},
because there the Grenander-type estimator is defined by only considering $\Lambda_n$ on a compact interval not containing the end point of the support.
The event $E_n$ is needed in our setup because of the presence of the covariates, which lead to more complicated processes, and because
we require~\eqref{eq:sup E} for $p=2$.

\begin{proof}[Proof of Lemma~\ref{le:inv}]
First, we note that, from the definition of $U$ and the fact that $\tilde{U}_n$ is increasing, it follows that
$|\tilde{U}_n(a)-U(a)|
\leq
|\tilde{U}_n(\lambda_0(0))-U(\lambda_0(0))|$,
if $a\leq\lambda_0(0)$,
and
\[
\1_{\{\tilde{U}_n(a)\leq M\}}|\tilde{U}_n(a)-U(a)|
\leq
\1_{\{\tilde{U}_n(a)\leq M\}}|\tilde{U}_n(\lambda_0(M))-U(\lambda_0(M))|,\qquad \text{if } a\geq\lambda_0(M),
\]
Hence, it suffices to prove~\eqref{eqn:inv} only for $a\in[ \lambda_0(0),\lambda_0(M)]$.
Let $E_n$ be the event from Lemma~\ref{le:2-1}.
We start by writing
\begin{equation}
\label{eqn:main}
\begin{split}
&
\p\left(
\left\{|\tilde{U}_n(a)-U(a)|\geq x\right\}\cap E_n\cap \left\{\tilde{U}_n(a)\leq M\right\}\right)\\
&\quad=
\p\left(\left\{
U(a)+x\leq \tilde{U}_n(a)\leq M \right\}\cap E_n\right)
+
\p\left(\left\{\tilde{U}_n(a)\leq U(a)-x\right\}\cap E_n\right).
\end{split}
\end{equation}
First consider the first probability on the right hand side of~\eqref{eqn:main}.
It is zero, if $U(a)+x>M$.
Otherwise, if $U(a)+x\leq M$, then $x\leq M$ and
\[
\begin{split}
&
\p\left(\left\{
U(a)+x\leq \tilde{U}_n(a)\leq M \right\}\cap E_n\right)\\
&\leq
\p\Big(\left\{
\Lambda_n(y)-ay\leq \Lambda_n(U(a))-aU(a),\text{for some }y\in[U(a)+x,M]
\right\}\cap E_n
\Big)\\
&\leq
\p\left(\left\{
\inf_{y\in[U(a)+x,M]}
\Big(
\Lambda_n(y)-ay-\Lambda_n(U(a))+a\,U(a)
\Big)
\leq 0\right\}\cap E_n\right).
\end{split}
\]
From Taylor's expansion, we obtain
$\Lambda_0(y)-\Lambda_0(U(a))\geq \big(y-U(a)\big)a+c \big(y-U(a)\big)^2$,
where $c=\inf_{t\in[0,\tau_F)}\lambda'_0(t)/2>0$,
so that with~\eqref{eq:martingale relation},
the probability on the right hand side is bounded by
\[
\p\left(\left\{
\inf_{y\in[U(a)+x,M]}
\Big(\w_n(y)-\w_n(U(a))+R_n(y)-R_n(U(a))+c(y-U(a))^2\Big)\leq 0\right\}\cap E_n\right).
\]%
Let $i\geq 0$ be such that $M-U(a)\in[x2^i,x2^{i+1})$ and note that, on the event $E_n$ one has $T_{(n)}\geq M$.
Therefore, if $U(a)<y\leq M$, then $y\leq T_{(n)}$ and $U(a)<T_{(n)}$.
It follows that the previous probability can be bounded by
\[
\sum_{k=0}^i
\p\left(\left\{
\sup_{y\in I_k}
\Big(
\big|\w_n(y)-\w_n(U(a))\big|+|R_n(y)-R_n(U(a))|
\Big)
\geq c\,x^2\,2^{2k}\right\}\cap E_n\right).
\]
where the supremum is taken over $y\in[0,M]$, such that $y-U(a)\in[x2^k,x2^{k+1})$.
Using that $\p(X+Y\geq \epsilon)\leq \p(X\geq \epsilon/2)+\p(Y\geq \epsilon/2)$, together with the Markov inequality, we can bound this probability by
\begin{equation}
\label{eqn:1}
\begin{split}
&
4\sum_{k=0}^i
\left(c^2x^42^{4k}\right)^{-1}
\E\left[\sup_{y\leq M,\,y-U(a)\in[x2^k,x2^{k+1})}\big|\w_n(y)-\w_n(U(a))\big|^2\right]\\
&\quad+
8\sum_{k=0}^i
\left(
c^3x^62^{6k}
\right)^{-1}
\E\left[\sup_{y<M,\,y-U(a)\in[x2^k,x2^{k+1})}\1_{E_n}\,\big|R_n(y)-R_n(U(a))\big|^3\right].
\end{split}
\end{equation}
We have
\begin{equation}
\label{eqn:bound r}
\begin{split}
&
\E\left[\sup_{y<M,\,y-U(a)\in[x2^k,x2^{k+1})}\1_{E_n}\,\big|R_n(y)-R_n(U(a))\big|^3\right]\\
&\leq
4\E\left[\sup_{y<M,\,y-U(a)\in[x2^k,x2^{k+1})}\1_{E_n}\,\left| \int_{U(a)}^y \left(\frac{\Phi_n(s;\beta_0)}{\Phi(s;\beta_0)}-1 \right)\lambda_0(s)\mathrm{d}s\right|^3\right]\\
&\quad+
4\E\left[\sup_{y<M,\,y-U(a)\in[x2^k,x2^{k+1})}\1_{E_n}\,\left| \int_{U(a)}^y \left(\frac{1}{\Phi_n(s;\hat{\beta}_n)}-\frac{1}{\Phi(s;\beta_0)} \right)\mathrm{d}H^{uc}_n(s)\right|^3\right]
\end{split}
\end{equation}
For the first term in the right hand side of~\eqref{eqn:bound r} we have
\[
\begin{split}
&
\E\left[\sup_{y<M,\,y-U(a)\in[x2^k,x2^{k+1})}\1_{E_n}
\left|
\int_{U(a)}^y \left(\frac{\Phi_n(s;\beta_0)}{\Phi(s;\beta_0)}-1 \right)\lambda_0(s)\mathrm{d}s
\right|^3\right]\\
&\quad\leq
\E\left[\1_{E_n}\left(\int_{U(a)}^{(U(a)+x2^{k+1})\wedge M}\frac{\left|\Phi_n(s;\beta_0)-\Phi(s;\beta_0)\right|}{\Phi(s;\beta_0)}\lambda_0(s)\,\mathrm{d}s\right)^3\right]\\
&\quad\leq
\frac{x^32^{3(k+1)}\lambda^3_0(M)}{\Phi(M;\beta_0)^3}
\E\left[\1_{E_n}\,\sup_{s\in[0,M]}\left|\Phi_n(s;\beta_0)-\Phi(s;\beta_0)\right|^3\right]
\leq
\frac{x^32^{3(k+1)}\,\lambda^3_0(M)\xi_4}{n\Phi(M;\beta_0)^3},
\end{split}
\]
where we have used~\eqref{eqn:ef}.
In order to bound the second term on the right hand side of~\eqref{eqn:bound r}, note that on the event $E_n$,
\begin{equation}
\label{eqn:E_n5}
\begin{split}
\sup_{x\in\R}|\Phi_n(x;\hat{\beta}_n)-\Phi(x;\beta_0)|
&\leq
\sup_{x\in\R}|\Phi_n(x;\hat{\beta}_n)-\Phi_n(x;\beta_0)|+\sup_{x\in\R}|\Phi_n(x;\beta_0)-\Phi(x;\beta_0)|\\
&\leq
|\hat{\beta}_n-\beta_0|\sup_{x\in\R}|D_n^{(1)}(x;\beta^*)|+ \frac{\xi_4}{n^{1/3}}
\leq
\frac{\sqrt{\xi_2}L+\xi_4}{n^{1/3}}.
\end{split}
\end{equation}
In particular, for sufficiently large $n$ we have
$\sup_{x\in\R}\left|\Phi_n(x;\hat{\beta}_n)-\Phi(x;\beta_0)\right|\leq \Phi(M;\beta_0)/2$,
which yields that, for $x\in[0,M]$,
\begin{equation}
\label{eq:bound Phihat}
\Phi_n(x;\hat{\beta}_n)
\geq
\Phi(x;\beta_0)- \frac{1}{2}\Phi(M;\beta_0)
\geq \frac{1}{2}\Phi(M;\beta_0).
\end{equation}
Using~\eqref{eqn:E_n5}, on the event $E_n$, for $n$ sufficiently large, we can write
\[
\begin{split}
\sup_{s\in[0,M]}\left|\frac{1}{\Phi_n(s;\hat{\beta}_n)}-\frac{1}{\Phi(s;\beta_0)} \right|
&\leq
\sup_{s\in[0,M]}\frac{\left|\Phi_n(s;\hat{\beta}_n)-\Phi(s;\beta_0)\right|}{\Phi_n(s;\hat{\beta}_n)\,\Phi(s;\beta_0)}\\
&\leq
\frac{2}{\Phi^2(M;\beta_0)}\sup_{s\in[0,M]}\left|\Phi_n(s;\hat{\beta}_n)-\Phi(s;\beta_0)\right|\\
&\leq Cn^{-1/3},
\end{split}
\]
for some $C>0$.
Consequently, for the second term in the right hand side of~\eqref{eqn:bound r} we obtain
\[
\begin{split}
&
\E\left[\sup_{y<M,\,y-U(a)\in[x2^k,x2^{k+1})}\1_{E_n}
\left| \int_{U(a)}^y \left(\frac{1}{\Phi_n(s;\hat{\beta}_n)}-\frac{1}{\Phi(s;\beta_0)} \right)\mathrm{d}H^{uc}_n(s)\right|^3\right]\\
&\quad\leq
\frac{C^3}{n}\E\left[\1_{E_n}\left(
\frac{1}{n}\sum_{i=1}^n \Delta_i \1_{\{T_i\in[U(a),(U(a)+x2^{k+1})\wedge M) ]\}}
\right)^3\right]
\leq
\frac{C^3}{n^4}
\E[N^3],
\end{split}
\]
where $N$ is a binomial distribution with probability of success
\[
\gamma
=
H^{uc}((U(a)+x2^{k+1})\wedge M))-H^{uc}(U(a))
\leq
\sup_{s\in[0,M]}|h^{uc}(s)|
x2^{k+1}.
\]
Furthermore,
\[
\mathbb{E}[N^3]=n\gamma(1-3\gamma+3n\gamma+2\gamma^2-3n\gamma^2+n^2\gamma^2)
\leq
\begin{cases}
7n\gamma & \text{, if }n\gamma\leq 1;\\
7n^3\gamma^3 & \text{, if } n\gamma>1.
\end{cases}
\]
Using Lemma~\ref{le:mart} and the bound in~\eqref{eqn:1}, for the first probability on the right hand side of~\eqref{eqn:main},
it follows that there exist $K_1,K_2>0$, such that for all $a\geq 0$, $n\geq 1$ and $x>0$,
\begin{equation}
\label{eqn:inv1}
\begin{split}
&
\p\left(\left\{U(a)+x\leq \tilde{U}_n(a)\leq M \right\}\cap E_n\right)\\
&\quad\leq
K_1
\sum_{k=0}^i\frac{x2^{k+1}}{nx^42^{4k}}
+
K_2\sum_{k=0}^i
\max\left\{
\frac{x2^{k+1}}{n^3x^62^{6k}},\frac{x^32^{3(k+1)}}{nx^62^{6k}}
\right\}\\
&\quad\leq
\frac{2K_1}{nx^3}\sum_{k=0}^\infty 2^{-3k}
+
\max
\left\{
\frac{2K_2}{n^3x^5}\sum_{k=0}^\infty 2^{-5k},
\frac{8K_2}{nx^3}\sum_{k=0}^\infty 2^{-3k}
\right\}
\leq
K
\max\left\{
\frac{1}{nx^3},\frac{1}{n^3x^5}
\right\}.
\end{split}
\end{equation}
We proceed with the second probability on the right hand side of~\eqref{eqn:main}.
We can assume $x\leq U(a)$, because otherwise $\p\big(\tilde{U}_n(a)\leq U(a)-x\big)=0$.
We have
\[
\begin{split}
&
\p\left(\left\{\tilde{U}_n(a)\leq U(a)-x\right\}\cap E_n\right)\\
&\leq
\p\left(\left\{\inf_{y\in[0,U(a)-x]}\big[\Lambda_n(y)-ay-\Lambda_n(U(a))+a\,U(a)\big]\leq 0\right\}\cap E_n\right).
\end{split}
\]
Let $i\geq 0$ be such that $U(a)\in[x2^{i},x2^{i+1})$.
By a similar argument used to obtain the bound~\eqref{eqn:1},
this probability is bounded by
\begin{equation}
\label{eqn:2}
\begin{split}
&
4\sum_{k=0}^i
\left(c^2x^42^{4k}\right)^{-1}
\E\left[
\sup_{y\leq U(a),U(a)-y\in[x2^{k},x2^{k+1})}\big|\w_n(y)-\w_n(U(a))\big|^2\right]\\
&\quad
+8\sum_{k=0}^i
\left(
c^3x^62^{6k}
\right)^{-1}
\E\left[\sup_{y\leq U(a),U(a)-y\in[x2^{k},x2^{k+1})}
\1_{E_n}
\big|R_n(y)-R_n(U(a))\big|^3\right].
\end{split}
\end{equation}
In the same way as in the first case, we also have
\[
\E\left[\sup_{y\leq U(a),U(a)-y\in[x2^{k},x2^{k+1})}
\1_{E_n}
\big|R_n(y)-R_n(U(a))\big|^3\right]
\leq
K_2\max\left\{\frac{x2^{k+1}}{n^3},\frac{x^32^{3(k+1)}}{n}\right\}.
\]
Exactly as in~\eqref{eqn:inv1}, Lemma~\ref{le:mart} and~\eqref{eqn:2} imply that
\[
\p
\left(\left\{
\tilde{U}_n(a)\leq U(a)-x
\right\}\cap E_n\right)\leq
K
\max\left\{
\frac{1}{nx^3},\frac{1}{n^3x^5}
\right\},
\]
for some positive constant $K$.
Together with~\eqref{eqn:main} and~\eqref{eqn:inv1}, this finishes the proof.
\end{proof}

\begin{lemma}
\label{le:lambda}
Suppose that (A1)--(A2) hold.
Let $0<\epsilon<M'<M<\tau_H$ and suppose that $H^{uc}$, defined in~\eqref{eq:def Huc}, has a bounded derivative $h^{uc}$ on $[0,M]$.
Let $\tilde{\lambda}_n$ be the Grenander-type estimator of a nondecreasing baseline hazard rate $\lambda_0$, which is differentiable with $\lambda_0'$ bounded above and below by strictly positive constants.
Let $E_n$ be the event from Lemma~\ref{le:2-1} and take $\xi_3$ in~\eqref{eqn:ef} such that
\begin{equation}
\label{eqn:xi}
0<\xi_3<
\frac{1}{8}
\min\left\{(M-M')^2,\epsilon^2\right\}
\inf_{x\in[0,\tau_H]}\lambda'_0(x).
\end{equation}
Then, there exists a constant $C$  such that, for $n$ sufficiently large,
\[
\sup_{t\in[\epsilon,M'] }
\E\left[n^{2/3}\1_{E_n}\left(\lambda_0(t)-\tilde{\lambda}_n(t)\right)^2\right]\leq C.
\]
\end{lemma}

\begin{proof}
It is sufficient to prove that there exist some constants $C_1,C_2>0$,
such that for each $n\in\N$ and each $t\in(\epsilon,M']$, we have
\begin{gather}
\label{eqn:exp1}
\E\left[n^{2/3}
\1_{E_n}
\left\{(\tilde{\lambda}_n(t)-\lambda_0(t))_+\right\}^2\right]\leq C_1,\\
\label{eqn:exp2}
\E\left[n^{2/3}\1_{E_n}
\left\{(\lambda_0(t)-\tilde{\lambda}_n(t))_+\right\}^2\right]\leq C_2.
\end{gather}
Lets first consider~\eqref{eqn:exp1}.
We will make use of the following result
\[
\begin{split}
\E\left[n^{2/3}
\1_{E_n}
\left\{(\tilde{\lambda}_n(t)-\lambda_0(t))_+\right\}^2\right]
&=
2\int_0^\infty
\p\left(
n^{1/3}\1_{E_n}(\tilde{\lambda}_n(t)-\lambda_0(t))\geq x
\right)x\,\mathrm{d}x\\
&=
2\int_0^{2\eta}
\p\left(n^{1/3}\1_{E_n}(\tilde{\lambda}_n(t)-\lambda_0(t))\geq x\right)x\,\mathrm{d}x\\
&\quad+
2\int_{2\eta}^\infty
\p\left(n^{1/3}\1_{E_n}(\tilde{\lambda}_n(t)-\lambda_0(t))\geq x\right)\,x\,\mathrm{d}x\\
&\leq
4\eta^2+2\int_{2\eta}^\infty
\p\left(n^{1/3}\1_{E_n}(\tilde{\lambda}_n(t)-\lambda_0(t))> x/2\right)x\,\mathrm{d}x\\
&\leq
4\eta^2+4\int_{\eta}^\infty
\p\left(n^{1/3}\1_{E_n}(\tilde{\lambda}_n(t)-\lambda_0(t))> x\right)x\,\mathrm{d}x\\
\end{split}
\]
for a fixed $\eta>0$.
We distinguish between the cases
$a+n^{-1/3}x\leq \lambda_0(M)$ and $a+n^{-1/3}x> \lambda_0(M)$,
where $a=\lambda_0(t)$.	
We prove that, in the first case, there exist a positive constant $C$ such that for all $t\in(\epsilon,M']$
and $n\in\N$,
$\p
\big(n^{1/3}\1_{E_n}(\tilde{\lambda}_n(t)-\lambda_0(t))> x\big)
\leq
C/x^3$,
for all $x\geq \eta$, and in the second case $\p\big(n^{1/3}\1_{E_n}(\tilde{\lambda}_n(t)-\lambda_0(t))> x\big)=0$.
Then~\eqref{eqn:exp1} follows immediately.	

First, assume $a+n^{-1/3}x\leq \lambda_0(M)$.
By the switching relation, we get
\[
\begin{split}
\p\left(
n^{1/3}\1_{E_n}(\tilde{\lambda}_n(t)-\lambda_0(t))> x\right)
&=
\p\left(\left\{
\tilde{\lambda}_n(t)> a+n^{-1/3}x\right\}\cap E_n\right)\\
&=
\p\left(\left\{\tilde{U}_n(a+n^{-1/3}x)< t\right\}\cap E_n\right).
\end{split}
\]
Because $a+n^{-1/3}x\leq \lambda_0(M)$, we have $U(a+n^{-1/3}x)\geq M>t$.
Furthermore, $\{\tilde{U}_n(a+n^{-1/3}x)< t\}\subset\{\tilde{U}_n(a+n^{-1/3}x)< M\}$.
Hence, together with Lemma~\ref{le:inv}, we can write
\begin{equation}
\label{eq:bound case 1}
\begin{split}
&
\p\left(\left\{\tilde{U}_n(a+n^{-1/3}x)< t\right\}\cap E_n\right)\\
&\quad\leq
\p\Bigg(
\Big\{
\left|U(a+n^{-1/3}x)-\tilde{U}_n(a+n^{-1/3}x)\right|>U(a+n^{-1/3}x)-t
\Big\}\\
&\qquad\qquad\qquad\qquad\qquad\qquad\qquad\qquad\qquad
\cap E_n\cap\left\{\tilde{U}_n(a+n^{-1/3}x)<M\right\}
\Bigg)\\
&\quad\leq
K
\max
\left\{
\frac{1}{n\big(U(a+n^{-1/3}x)-t\big)^3},\frac{1}{n^3\big(U(a+n^{-1/3}x)-t\big)^5}
\right\}
\leq
\frac{C}{x^3},
\end{split}
\end{equation}
because
$U(a+n^{-1/3}x)-t=U'(\xi_n)n^{-1/3}x$, for some $\xi_n\in(a,a+n^{-1/3}x)$,
where $U'(\xi_n)=\lambda'_0(\lambda_0^{-1}(\xi_n))^{-1}\geq 1/\sup_{t\in[0,\tau_H]}\lambda_0'(t)>0$.
	
Next, consider the case $a+n^{-1/3}x> \lambda_0(M)$.
Note that, we cannot argue as in the previous case,
because for $a+n^{-1/3}x> \lambda_0(M)$ we always have $U(a+n^{-1/3}x)=M$,	so that we loose the dependence on $x$.
However, if $n^{1/3}(\tilde{\lambda}_n(t)-\lambda_0(t))> x$,
then for each $y>t$, we have
\[
\tilde{\Lambda}_n(y)-\tilde{\Lambda}_n(t)
\geq
\tilde{\lambda}_n(t)\,(y-t)
>
(a+n^{-1/3}x)\,(y-t),	
\]
where $a=\lambda_0(t)$.
In particular for $y=\tilde{M}=M'+(M-M')/2$, we obtain
\begin{equation}
\label{eq:zero prob1}
\begin{split}
&
\p\left\{n^{1/3}\1_{E_n}(\tilde{\lambda}_n(t)-\lambda_0(t))> x\right\}\\
&\leq
\p\left(\left\{\tilde{\Lambda}_n(\tilde{M})-\tilde{\Lambda}_n(t)
>
\left(a+n^{-1/3}x\right)(\tilde{M}-t)\right\}\cap E_n\right)\\
&\leq
\p\Bigg(\Bigg\{\tilde{\Lambda}_n(\tilde{M})-\tilde{\Lambda}_n(t)-\big(\Lambda_0(\tilde{M})-\Lambda_0(t)\big)>\\
&\qquad\qquad\qquad\qquad\qquad
\left(a+n^{-1/3}x\right)(\tilde{M}-t)-\big(\Lambda_0(\tilde{M})-\Lambda_0(t)\big)\Bigg\}\cap E_n\Bigg)\\
&\leq
\p\left(\left\{
2\sup_{x\in[0,M]}|\tilde{\Lambda}_n(x)-\Lambda_0(x) |
>
\left(a+n^{-1/3}x-\lambda_0(\tilde{M})\right)(\tilde{M}-t)\right\}\cap E_n\right),
\end{split}
\end{equation}
also using that $\Lambda_0(\tilde{M})-\Lambda_0(t)\geq \lambda_0(\tilde M)(\tilde M-t)$.
Furthermore, since $a+n^{-1/3}x> \lambda_0(M)$, it follows from~\eqref{eqn:xi} that
\begin{equation}
\label{eq:prob zero}
\left(a+n^{-1/3}x-\lambda_0(\tilde{M})\right)(\tilde{M}-t)
\geq
\frac{1}{4}(M-M')^2\inf_{x\in[0,\tau_H]}\lambda'_0(x)
\geq
2\xi_3,
\end{equation}
so that, by the definition of $\xi_3$ in~\eqref{eqn:ef}, the probability on the right hand side~\eqref{eq:zero prob1} is zero.
This concludes the proof of~\eqref{eqn:exp1}.

Next, we have to deal with ~\eqref{eqn:exp2}.
Arguing as in the proof of~\eqref{eqn:exp1}, we obtain
\[
\E\left[n^{2/3}\1_{E_n}
\left\{(\lambda_0(t)-\tilde{\lambda}_n(t))_+\right\}^2\right]
\leq
\eta^2+2\int_{\eta}^\infty
\p\left(
n^{1/3}\1_{E_n}(\lambda_0(t)-\tilde{\lambda}_n(t))\geq x
\right)x\,\mathrm{d}x,
\]
for a fixed $\eta>0$,
where
\[
\p\left(
n^{1/3}\1_{E_n}(\lambda_0(t)-\tilde{\lambda}_n(t))\geq x
\right)
=
\p\left(\left\{\tilde{U}_n(a-n^{-1/3}x)\geq t\right\}\cap E_n\right),
\]
with $a=\lambda_0(t)$.
First of all, we can assume that $a-n^{-1/3}x\geq 0$, because otherwise $\p\{\tilde{\lambda}_n(t)\leq a-n^{-1/3}x\}=0$.
Since $t=U(a)$, as before, we write
\[
\begin{split}
&
\p\left(\left\{\tilde{U}_n(a-n^{-1/3}x)\geq t\right\}\cap E_n\right)\\
&\quad\leq
\p\left(\left\{
\left|\tilde{U}_n(a-n^{-1/3}x)-U(a-n^{-1/3}x)\right|\geq
t-U(a-n^{-1/3}x)\right\}\cap E_n\right).
\end{split}
\]
In order to apply Lemma~\ref{le:inv}, we intersect with the event $\tilde{U}_n(a-n^{-1/3}x)\leq M$.
Note that
\[
\p\left(\left\{\tilde{U}_n(a-n^{-1/3}x)> M\right\}\cap E_n\right)
\leq
\p\left(\left\{\tilde{\lambda}_n(M)\leq a-n^{-1/3}x\right\}\cap E_n\right)=0.
\]
This can be seen as follows.
If $\tilde{\lambda}_n(M)\leq a-n^{-1/3}x $, then for each $y<M$, we have
\[
\tilde{\Lambda}_n(M)-\tilde{\Lambda}_n(y)
\leq
\tilde{\lambda}_n(M)(M-y)
\leq
(a-n^{-1/3}x)(M-y).
\]
In particular for $y=\tilde{M}=M'+(M-M')/2$, similar to~\eqref{eq:zero prob1},
we obtain
\[
\begin{split}
&
\p\left(\left\{\tilde{\lambda}_n(M)\leq a-n^{-1/3}x\right\}\cap E_n\right)\\
&\quad\leq
\p\left(\left\{
2\sup_{x\in[0,M]}|\tilde{\Lambda}_n(x)-\Lambda_0(x)|
\geq
\left(-a+n^{-1/3}x+\lambda_0(\tilde{M})\right)
\left(M-\tilde{M}\right)\right\}\cap E_n\right).
\end{split}
\]
Because $a=\lambda_0(t)\leq \lambda_0(M')$,
we can argue as in~\eqref{eq:prob zero} and conclude that the probability
on the right hand side is zero.
It follows that
\[
\begin{split}
&
\p\left(\left\{\tilde{U}_n(a-n^{-1/3}x)\geq t\right\}\cap E_n\right)\\
&\quad\leq
\p\Bigg(
\left\{\left|\tilde{U}_n(a-n^{-1/3}x)-U(a-n^{-1/3}x)\right|\geq t-U(a-n^{-1/3}x)\right\}\\
&\qquad\qquad\qquad\qquad\qquad\qquad\qquad\qquad\qquad
\cap E_n\cap
\left\{
\tilde{U}_n(a-n^{-1/3}x)\leq M
\right\}
\Bigg)\\
&\quad\leq
K
\max
\left\{
\frac{1}{n\big(t-U(a-n^{-1/3}x)\big)^3},\frac{1}{n^3\big(t-U(a-n^{-1/3}x)\big)^5}
\right\}.
\end{split}
\]
To bound the right hand side, we have to distinguish between $a-n^{-1/3}x>\lambda_0(0)$
and $a-n^{-1/3}x\leq \lambda_0(0)$.
If $a-n^{-1/3}x>\lambda_0(0)$, then the right hand side is bounded by $K/x^3$,
because $t-U(a-n^{-1/3}x)=U'(\xi_n)n^{-1/3}x$, for some $\xi_n\in(a-n^{-1/3}x,a)$,
where $U'(\xi_n)=\lambda'_0(\lambda_0^{-1}(\xi_n))^{-1}\geq 1/\sup_{t\in[0,\tau_H]}\lambda_0'(t)>0$.
Otherwise, if $a-n^{-1/3}x\leq\lambda_0(0)$, then we are done because then
$\p\big(
n^{1/3}\1_{E_n}(\lambda_0(t)-\tilde{\lambda}_n(t))\geq x
\big)
=0$.
This can be seen as follows.
When $a-n^{-1/3}x\leq\lambda_0(0)$, then for each $y<t$, we have
\[
\tilde{\Lambda}_n(t)-\tilde{\Lambda}_n(y)
\leq
\tilde{\lambda}_n(t)(t-y)
\leq
(a-n^{-1/3}x)(t-y).	
\]
In particular, for $y=\epsilon'=\epsilon/2$, we obtain
\[
\begin{split}
&
\p\left(
n^{1/3}\1_{E_n}\left(\lambda_0(t)-\tilde{\lambda}_n(t)\right)\geq x
\right)\\
&\leq
\p\left(\left\{
\tilde{\Lambda}_n(t)-\tilde{\Lambda}_n(\epsilon')\leq \left(a-n^{-1/3}x\right)\left(t-\epsilon'\right)
\right\}\cap E_n\right)\\
&\leq
\p\Bigg(\Bigg\{
\tilde{\Lambda}_n(t)-\tilde{\Lambda}_n(\epsilon')-\left(\Lambda_0(t)-\Lambda_0(\epsilon')\right)\\
&\qquad\qquad\qquad\qquad\leq
\left(a-n^{-1/3}x\right)(t-\epsilon')-\left(\Lambda_0(t)-\Lambda_0(\epsilon')\right)
\Bigg\}\cap E_n\Bigg)\\
&\leq
\p\left(\left\{
2\sup_{x\in[0,M]}|\tilde{\Lambda}_n(x)-\Lambda_0(x)|
\geq
\left(-a+n^{-1/3}x\right)
\left(t-\epsilon'\right)+\lambda_0(\epsilon')\left(t-\epsilon'\right)\right\}\cap E_n\right).
\end{split}
\]
Because $a-n^{-1/3}x\leq\lambda_0(0)$,
we can argue as in~\eqref{eq:prob zero},
\begin{equation}
\label{eq:prob zero2}
\begin{split}
\left(-a+n^{-1/3}x+\lambda_0(\epsilon')\right)(t-\epsilon')
&\geq
\left(\lambda_0(\epsilon')-\lambda_0(0)\right)\left(\epsilon-\epsilon'\right)\\
&\geq
\frac{1}{4}\epsilon^2\inf_{x\in[0,\tau_H]}\lambda'_0(x)
\geq
2\xi_3.
\end{split}
\end{equation}
and conclude that the probability on the right hand side is zero.
This concludes the proof of~\eqref{eqn:exp2}.
\end{proof}

\begin{lemma}
\label{le:lambda2}
Suppose that (A1)--(A2) hold.
Fix $x\in(0,\tau_h)$.
Let $0<\epsilon<x<M'<M<\tau_H$ and suppose that $H^{uc}$, defined in~\eqref{eq:def Huc}, has a bounded derivative $h^{uc}$ on $[0,M]$.
Let $\tilde{\lambda}_n$ be the Grenander-type estimator of a nondecreasing baseline hazard rate $\lambda_0$,
which is differentiable with $\lambda'_0$ bounded above and below by strictly positive constants.
Let $E_n$ be the event from Lemma~\ref{le:2-1} and assume that $\xi_3$ satisfies~\eqref{eqn:xi}.
Then
\[
\1_{E_n}\int_{x-b}^{x+b}(\lambda_0(t)-\tilde{\lambda}_n(t))^2\mathrm{d}t=O_p(bn^{-2/3}).
\]
\end{lemma}
\begin{proof}
Markov's inequality and Fubini, yield
\[
\begin{split}
&
\p\left(b^{-1}n^{2/3}\1_{E_n}\int_{x-b}^{x+b}(\lambda_0(t)-\tilde{\lambda}_n(t))^2\,\mathrm{d}t>K\right)\\
&\quad\leq
\frac{1}{K}
\E\left[
b^{-1}n^{2/3}\1_{E_n}
\int_{x-b}^{x+b}(\lambda_0(t)-\tilde{\lambda}_n(t))^2\,\mathrm{d}t
\right]\\
&\quad\leq
\frac2K
\sup_{t\in[x-b,x+b]}
\E\left[
n^{2/3}\1_{E_n}(\lambda_0(t)-\tilde{\lambda}_n(t))^2
\right].
\end{split}
\]
For $n$ sufficiently large $[x-b,x+b]\subset [\epsilon, M']$, so that according to Lemma~\ref{le:lambda},
the right hand side is bounded by $2C/K$, for some constant $C>0$.
This proves the lemma.
\end{proof}
\begin{proof}[Proof of Lemma~\ref{le:result4}]
Take $x<M<\tau_H$ and $n$ sufficiently large such that $x+b\leq M$.
With~$\overline{a}_{n,x}\overline{\Phi}_n$ defined in~\eqref{def:piecewise constant}, we have
\begin{equation}
\label{eq:decomp thetabar}
\begin{split}
&
\int\left\{
\overline{\theta}_{n,x}(u,\delta,z)-\theta_{n,x}(u,\delta,z)
\right\}\,\mathrm{d}\p(u,\delta,z)\\
&\quad=
\1_{E_n}\int \1_{[x-b,x+b]}(u)\delta
\left(
\frac{\overline{a}_{n,x}\overline{\Phi}_n(u;\hat{\beta}_n)}{\Phi_n(u;\hat{\beta}_n)}-a_{n,x}(u)
\right)\,\mathrm{d}\p(u,\delta,z)\\
&\quad\qquad-
\1_{E_n}\int \mathrm{e}^{\hat{\beta}'_nz}
\int_0^u
\left(
\overline{\Psi}_{n,x}(v)-a_{n,x}(v)
\right)\mathrm{d}\,\tilde{\Lambda}_n(v)\,\mathrm{d}\p(u,\delta,z)\\
&\quad=
\1_{E_n}\int_{x-b}^{x+b}
\left(
\frac{\overline{a}_{n,x}\overline{\Phi}_n(u;\hat{\beta}_n)}{\Phi_n(u;\hat{\beta}_n)}-a_{n,x}(u)
\right)\,\mathrm{d}H^{uc}(u)\\
&\quad\qquad-
\1_{E_n}\int_{x-b}^{x+b}
\left(
\frac{\overline{a}_{n,x}\overline{\Phi}_n(v;\hat{\beta}_n)}{\Phi_n(v;\hat{\beta}_n)}-a_{n,x}(v)
\right)\Phi(v;\hat{\beta}_n)\,\mathrm{d}\tilde{\Lambda}_n(v)
\end{split}
\end{equation}
using Fubini, the definition~\eqref{eq:def Phi} of $\Phi$,
and the fact that $\overline{a}_{n,x}\overline{\Phi}_n$ and $a_{n,x}$ are zero outside $[x-b,x+b]$.
Write $\widehat{\Phi}_n(u)=\Phi_n(u;\hat\beta_n)$,
$\widehat{\Phi}(u)=\Phi(u;\hat\beta_n)$, and
$\Phi_0(u)=\Phi(u;\beta_0)$.
Then the right hand side can be written as
\[
\1_{E_n}\int_{x-b}^{x+b}
\frac{a_{n,x}(\hat{A}_n(u))\widehat{\Phi}_n(\hat{A}_n(u))-a_{n,x}(u)\widehat{\Phi}_n(u)}{\widehat{\Phi}_n(u)}
\left(
\Phi_0(u)\lambda_0(u)-\widehat{\Phi}(u)\tilde{\lambda}_n(u)
\right)\,\mathrm{d}u
\]
where $\hat{A}_n(u)$ is defined in~\eqref{eq:An}.
The Cauchy-Schwarz inequality then yields
\begin{equation}
\label{eqn:sch}
\begin{split}
&
\left|\int\left\{\overline{\theta}_{n,x}(u,\delta,z)-\theta_{n,x}(u,\delta,z)\right\}\,\mathrm{d}\p(u,\delta,z) \right|\\
&\leq
\1_{E_n}
\left\Vert
\frac{(a_{n,x}\circ\hat{A}_n)(\widehat{\Phi}_n\circ\hat{A}_n)-a_{n,x}\widehat{\Phi}_n}{\widehat{\Phi}_n}
\1_{[x-b,x+b]}
\right\Vert_{\mathcal{L}_2}
\left\Vert
\left(\Phi_0\lambda_0-\widehat{\Phi}\tilde{\lambda}_n\right)
\1_{[x-b,x+b]}
\right\Vert_{\mathcal{L}_2}.
\end{split}
\end{equation}
Furthermore,
\begin{equation}
\label{eqn:l2}
\begin{split}
&
\1_{E_n}\left\Vert
\frac{(a_{n,x}\circ\hat{A}_n)(\widehat{\Phi}_n\circ\hat{A}_n)-a_{n,x}\widehat{\Phi}_n}{\widehat{\Phi}_n}
\1_{[x-b,x+b]}
\right\Vert_{\mathcal{L}_2}\\
&\leq
2\1_{E_n}\int_{x-b}^{x+b}
\left(
\frac{k_b(x-\hat{A}_{n,x}(u))}{\Phi_0(\hat{A}_n(u))}-\frac{k_b(x-u)}{\Phi_0(u)}
\right)^2\,\mathrm{d}u\\
&\qquad+
2\1_{E_n}\int_{x-b}^{x+b}
\left(
\frac{k_b(x-\hat{A}_n(u))}{\Phi_0(\hat{A}_n(u))}
\right)^2
\frac{\left(\widehat{\Phi}_n(\hat{A}_n(u))-\widehat{\Phi}_n(u)\right)^2}{\widehat{\Phi}_n(u)^2}\,\mathrm{d}u\\
&\leq
2\1_{E_n}
\int_{x-b}^{x+b}
\left(
\frac{\mathrm{d}}{\mathrm{d}y}\frac{k_b(x-y)}{\Phi_0(y)}\bigg|_{y=\xi_u}
\right)^2
\left(\hat{A}_n(u)-u\right)^2\,\mathrm{d}u\\
&\qquad+
\1_{E_n}\frac{c_1}{b^2\widehat{\Phi}_n(M)^2}
\int_{x-b}^{x+b}
\left(
\widehat{\Phi}_n(\hat{A}_n(u))-\widehat{\Phi}_n(u)
\right)^2\,\mathrm{d}u\\
&\leq
\1_{E_n}\frac{c_2}{b^4}\int_{x-b}^{x+b}\left(\hat{A}_n(u)-u\right)^2\,\mathrm{d}u
+
\1_{E_n}\frac{c_1}{b^2\widehat{\Phi}_n(M)^2}
\int_{x-b}^{x+b}
\left(
\widehat{\Phi}_n(\hat{A}_n(u))-\widehat{\Phi}_n(u)
\right)^2\,\mathrm{d}u,
\end{split}
\end{equation}
for some constants $c_1,c_2>0$,
where we use the boundedness of $k'$, $\mathrm{d}\Phi(x;\beta_0)/\mathrm{d}x$ and $1/\Phi(x;\beta_0)$ on $[0,x+b]\subseteq[0,M]$.
Then, since
$\lambda_0(u)-\lambda_0(\hat{A}_n(u))=\lambda'_0(\xi)(u-\hat{A}_n(u))$
and $\lambda'_0$ is bounded and strictly positive on $[0,M]\supseteq[x-b,x+b]$, there exists a constant $K>0$ such that
\[
|u-\hat{A}_n(u)|\leq K\,|\lambda_0(u)-\lambda_0(\hat{A}_n(u))|.
\]
If $u\in(\tau_i,\tau_{i+1}]$ and $\hat{A}_n(u)>\tau_i$, then $\tilde{\lambda}_n(u)=\tilde{\lambda}_n(\hat{A}_n(u))$ and we obtain
\begin{equation}
\label{eqn:lambda}
\begin{split}
|u-\hat{A}_n(u)|	
&\leq
K|\lambda_0(u)-\tilde{\lambda}_n(u)|+K|\tilde{\lambda}_n(\hat{A}_n(u))-\lambda_0(\hat{A}_n(u))|\\
&\leq
2K|\lambda_0(u)-\tilde{\lambda}_n(u)|.
\end{split}
\end{equation}
This holds also in the case $\hat{A}_n(u)=\tau_i$,
simply because $|\lambda_0(u)-\lambda_0(\hat{A}_n(u))|\leq |\lambda_0(u)-\tilde{\lambda}_n(u)|$.
As a result, using Lemma~\ref{le:lambda2},
for the first term on the right hand side of~\eqref{eqn:l2},
we derive that
\begin{equation}
\label{eq:first term}
\1_{E_n}\frac{1}{b^4}\int_{x-b}^{x+b}
\left(\hat{A}_n(u)-u\right)^2\,\mathrm{d}u
\leq
\frac{C}{b^4}\1_{E_n}\int_{x-b}^{x+b}
\left(\lambda_0(u)-\tilde{\lambda}_n(u)\right)^2\,\mathrm{d}u
=
O_p(b^{-3}n^{-2/3}).
\end{equation}
For the second term on the right hand side of~\eqref{eqn:l2}, we find
\[
\begin{split}
\big|
\Phi_n(\hat{A}_n(u);\hat\beta_n)-\Phi_n(u;\hat{\beta}_n)\big|
&\leq
2\sup_{x\in\R}\big|\Phi_n(x;\hat{\beta}_n)-\Phi_n(x;\beta_0)\big|+\big|\Phi(\hat{A}_n(u);\beta_0)-\Phi(u;\beta_0)\big|\\
&\qquad+
2\sup_{x\in\R}|\Phi_n(x;\beta_0)-\Phi(x;\beta_0)|\\
&\leq
2|\hat{\beta}_n-\beta_0|\sup_{x\in\R}\big|D^{(1)}_n(x;\beta^*)\big|+|\Phi'(\xi;\beta_0)||\hat{A}_n(u)-u|\\
&\qquad+
2\sup_{x\in\R}|\Phi_n(x;\beta_0)-\Phi(x;\beta_0)|,
\end{split}
\]
which, using~\eqref{eqn:dn2},~\eqref{eqn:Phi} and $|\hat{\beta}_n-\beta_0|=O_p(n^{-1/2})$ (see Theorem 3.2 in~\cite{Tsiatis81}), leads to
\[
\begin{split}
&
b^{-2}\1_{E_n}\int_{x-b}^{x+b}
\left(
\Phi_n(\hat{A}_n(u);\hat\beta_n)-\Phi_n(u;\hat{\beta}_n)
\right)^2\,\mathrm{d}u\\
&\qquad\leq
8b^{-1}\1_{E_n}|\hat{\beta}_n-\beta_0|^2\sup_{x\in\R}\big|D^{(1)}_n(x;\beta^*)\big|^2\\
&\qquad\quad+
2b^{-2}\1_{E_n}\sup_{s\in[x-b,x+b]}\Phi'(s;\beta_0)
\int_{x-b}^{x+b}
\left(\hat{A}_n(u)-u\right)^2\,\mathrm{d}u\\
&\qquad\qquad+
8b^{-1}
\sup_{x\in\R}|\Phi_n(x;\beta_0)-\Phi(x;\beta_0)|^2\\
&\qquad=
O_p(b^{-1}n^{-1})+O_p(b^{-1}n^{-2/3})+O_p(b^{-1}n^{-1})
=
O_p(b^{-1}n^{-2/3}).
\end{split}
\]
Consequently, from~\eqref{eqn:l2} together with~\eqref{eq:bound Phihat},
for the first term on the right hand side of~\eqref{eqn:sch}, we obtain
\[
\begin{split}
\1_{E_n}\left\Vert
\frac{(a_{n,x}\circ\hat{A}_n)(\widehat{\Phi}_n\circ\hat{A}_n)-a_{n,x}\widehat{\Phi}_n}{\widehat{\Phi}_n}
\1_{[x-b,x+b]}
\right\Vert_{\mathcal{L}_2}^2
&=
O_p(b^{-3}n^{-2/3})+O_p(b^{-1}n^{-2/3})\\
&=O_p(b^{-3}n^{-2/3}).
\end{split}
\]
For the second term on the right hand side of~\eqref{eqn:sch}, we first write
\begin{equation}
\label{eq:bound Phi0lambda0}
\begin{split}
&
\left\Vert
\left(\Phi_0\lambda_0-\widehat{\Phi}\tilde{\lambda}_n\right)
\1_{[x-b,x+b]}
\right\Vert_{\mathcal{L}_2}\\
&\leq
\left\Vert
\left(\Phi_0-\widehat{\Phi}\right)\lambda_0
\1_{[x-b,x+b]}
\right\Vert_{\mathcal{L}_2}
+
\left\Vert
\left(\lambda_0-\tilde{\lambda}_n\right)\widehat{\Phi}
\1_{[x-b,x+b]}
\right\Vert_{\mathcal{L}_2}
\end{split}
\end{equation}
On the event $E_n$, we find
\[
\begin{split}
\1_{E_n}
\left\Vert
\left(\Phi_0-\widehat{\Phi}\right)\lambda_0
\1_{[x-b,x+b]}
\right\Vert_{\mathcal{L}_2}^2
&\leq
2b\1_{E_n}|\hat{\beta}_n-\beta_0|^2\sup_{x\in\R}\big|D^{(1)}_n(x;\beta^*)\big|^2\sup_{u\in[0,M]}\lambda_0(u)\\
&=
O_p(bn^{-1}),
\end{split}
\]
and
\[
\1_{E_n}
\left\Vert
\left(\lambda_0-\tilde{\lambda}_n\right)\widehat{\Phi}
\1_{[x-b,x+b]}
\right\Vert_{\mathcal{L}_2}^2
\leq
\Phi(0,\hat\beta_n)
\1_{E_n}
\int_{x-b}^{x+b}
\left(\lambda_0(u)-\tilde{\lambda}_n(u)\right)^2
\,\mathrm{d}u
=
O_p(bn^{-2/3}),
\]
due to Lemma~\ref{le:lambda2}.
It follows that
\[
\1_{E_n}
\left\Vert
\left(\Phi_0\lambda_0-\widehat{\Phi}\tilde{\lambda}_n\right)
\1_{[x-b,x+b]}
\right\Vert_{\mathcal{L}_2}
=
O_p(b^{1/2}n^{-1/3}).
\]
Together with~\eqref{eqn:sch}, this concludes the proof.
\end{proof}

To establish Lemma~\ref{le:result5}, we need a slightly stronger version of Lemma~\ref{le:lambda}.
Note that, in order to have the uniform result in~\eqref{eqn:lemma},
we loose a factor $n^{2/9}$ with respect to the bound in Lemma~\ref{le:lambda}.
This might not be optimal, but it is sufficient for our purposes.
\begin{lemma}
\label{le:sup}
Suppose that (A1)--(A2) hold.
Let $0<\epsilon<M'<M<\tau_H$ and suppose that $H^{uc}$, defined in~\eqref{eq:def Huc}, has a bounded derivative $h^{uc}$ on $[0,M]$.
Let $\tilde{\lambda}_n$ be the Grenander-type estimator of a nondecreasing baseline hazard rate $\lambda_0$,
which is differentiable with $\lambda'_0$ bounded above and below by strictly positive constants.
Let $E_n$ be the event from Lemma~\ref{le:2-1} and assume that $\xi_3$ satisfies~\eqref{eqn:xi}.
Then, there exists a constant $C>0$, such that, for each $n\in\N$,
\begin{equation}
\label{eqn:lemma}
\E\left[
n^{4/9}\1_{E_n}\sup_{t\in(\epsilon,M'] }\left(\lambda_0(t)-\tilde{\lambda}_n(t)\right)^2
\right]
\leq C.
\end{equation}
\end{lemma}
\begin{proof}
We decompose $(\epsilon,M']$ in $m$ intervals $(c_k,c_{k+1}]$, where
$c_k=
\epsilon+k(M'-\epsilon)/m$,for $k=0,1,\dots,m$,
and $m=(M'-\epsilon)n^{2/9}$.
Then, we have
\[
\sup_{t\in(\epsilon,M']}
\left(
\lambda_0(t)-\tilde{\lambda}_n(t)
\right)^2
=
\max_{0\leq k\leq m-1}
\sup_{t\in(c_k,c_{k+1}]}
\left(\lambda_0(t)-\tilde{\lambda}_n(t)\right)^2.
\]
On the other hand, the fact that $\lambda_0$ is differentiable with bounded derivative implies that
\[
\begin{split}
&
\sup_{t\in(c_k,c_{k+1}]}
\big(\lambda_0(t)-\tilde{\lambda}_n(t)\big)^2\\
&\quad\leq
2\sup_{t\in(c_k,c_{k+1}]}\big(\lambda_0(t)-\lambda_0(c_{k+1})\big)^2
+
2\sup_{t\in(c_k,c_{k+1}]}\big(\lambda_0(c_{k+1})-\tilde{\lambda}_n(t)\big)^2\\
&\quad\leq
2\left(\sup_{u\in[0,M']}\lambda'_0(u)\right)^2\big(c_k-c_{k+1}\big)^2\\
&\qquad\qquad+
2\max\left\{
\big(\lambda_0(c_{k+1})-\tilde{\lambda}_n(c_k)\big)^2,\big(\lambda_0(c_{k+1})-\tilde{\lambda}_n(c_{k+1})\big)^2
\right\}\\
&\quad\leq
2\left(
\sup_{u\in[0,M']}\lambda'_0(u)
\right)^2\frac{(M'-\epsilon)^2}{m^2}\\
&\qquad\qquad+
\max_{0\leq k\leq m}\big(\lambda_0(c_{k})-\tilde{\lambda}_n(c_{k})\big)^2
+
\max_{0\leq k\leq m-1}\big(\lambda_0(c_{k+1})-\lambda_0(c_k)\big)^2.
\end{split}
\]
Here we used that $\tilde{\lambda}_n$ is non-decreasing,
and therefore $\sup_{t\in(c_k,c_{k+1}]}\big(\lambda_0(c_{k+1})-\tilde{\lambda}_n(t)\big)^2 $ is achieved either at $t=c_k$ or $t=c_{k+1}$, for $k=0,1,\dots,m-1$.
Hence,
\[
\begin{split}
\sup_{t\in(\epsilon,M']}
\big(\lambda_0(t)-\tilde{\lambda}_n(t)\big)^2
&\leq
4\max_{0\leq k\leq m}
\big(\lambda_0(c_{k})-\tilde{\lambda}_n(c_{k})\big)^2
+
6
\left(
\sup_{u\in[0,M']}\lambda'_0(u)
\right)^2\frac{(M'-\epsilon)^2}{m^2}\\
&\leq
4\max_{0\leq k\leq m}
\big(\lambda_0(c_{k})-\tilde{\lambda}_n(c_{k})\big)^2
+
C_1n^{-4/9},
\end{split}
\]
where $C_1=6\big(\sup_{u\in[\epsilon,M']}\lambda'_0(u)\big)^2$.
Consequently, using Lemma~\ref{le:lambda}, we derive
\[
\begin{split}
\E\left[
n^{4/9}\1_{E_n}\sup_{t\in(\epsilon,M']}\big(\lambda_0(t)-\tilde{\lambda}_n(t)\big)^2
\right]
&\leq
4\E\left[
n^{4/9}\1_{E_n}\max_{0\leq k\leq m}\big(\lambda_0(c_{k})-\tilde{\lambda}_n(c_{k})\big)^2
\right]+C_1\\
&\leq
4n^{-2/9}\sum_{k=0}^m
\E\left[n^{2/3}\1_{E_n}\big(\lambda_0(c_{k})-\tilde{\lambda}_n(c_{k})\big)^2\right]+C_1\\
&\leq
4\left(M'-\epsilon+1\right)C+C_1.
\end{split}
\]
This concludes the proof of~\eqref{eqn:lemma}.
\end{proof}

\begin{proof}[Proof of Lemma~\ref{le:result5}]
Let $n$ be sufficiently large such that $x+b\leq M'<M<\tau_H$.
Denote by~$R_n$ the left hand side of~\eqref{eqn:statement}.
Write $\mathbb{G}_n=\sqrt{n}(\mathbb{P}_n-\mathbb{P})$ and decompose $R_n=R_{n1}+R_{n2}$,
where
\[
\begin{split}
R_{n1}
&=
n^{-1/2}
\1_{E_n}\int\frac{\delta\1_{[x-b,x+b]}(u)}{\Phi_n(u;\hat{\beta}_n) }
\left(
\bar{a}_{n,x}\bar{\Phi}_n(u;\hat{\beta}_n)-a_{n,x}(u)\Phi_n(u;\hat{\beta}_n)
\right)\,\mathrm{d}\mathbb{G}_n(u,\delta,z),\\
R_{n2}
&=
n^{-1/2}
\1_{E_n}
\int \1_{\{u>x-b\}}
\Bigg[e^{\hat{\beta}'_nz}
\int_{x-b}^{u\wedge (x+b)}
\frac{\bar{a}_{n,x}\bar{\Phi}_n(v;\hat{\beta}_n)}{\Phi_n(v;\hat{\beta}_n)}\,\mathrm{d}\tilde{\Lambda}_n(v)\\
&\qquad\qquad\qquad\qquad\qquad\qquad\qquad\qquad-
\mathrm{e}^{\beta'_0z}
\int_{x-b}^{u\wedge (x+b)}
 a_{n,x}(v)\,\mathrm{d}\Lambda_0(v)
\Bigg]\,\mathrm{d}\mathbb{G}_n(u,\delta,z).
\end{split}
\]
Choose $\eta>0$.
We prove separately that there exists $\nu>0$, such that
\begin{equation}
\label{eq:bigO Rn}
\begin{split}
\limsup_{n\to\infty}
\p\left(b^{3/2}n^{13/18}|R_{n1}|>\nu\right)&\leq \eta\\
\limsup_{n\to\infty}
\p\left(n^{1/2}|R_{n2}|>\nu\right)&\leq \eta.
\end{split}
\end{equation}
We consider the following events.
\begin{equation}
\label{def:events}
\begin{split}
\A_{n1}&=\left\{\tilde{\lambda}_n(M)>K_1\right\},\\
\A_{n2}&=\left\{\sup_{t\in[\epsilon,M']}\left|\lambda_0(t)-\tilde{\lambda}_n(t)\right|>K_2n^{-2/9}\right\},
\end{split}
\end{equation}
where $K_1,K_2>0$, and let $\A_n=\A_{n1}\cup\A_{n2}$.
Since $\tilde{\lambda}_n(M)=O_p(1)$ we can choose $K_1>0$, such that $\p(\A_{n1})\leq\eta/3$ and
from Lemma~\ref{le:sup} we find that we can choose $K_2>0$, such that $\p(\A_{n2})\leq\eta/3$,
so that $\p(\A_n)\leq 2\eta/3$.
First consider $R_{n1}$.
Since
\begin{equation}
\label{eq:Rn1 Markov}
\begin{split}
\p\left(b^{3/2}n^{13/18}|R_{n1}|>\nu\right)
&\leq
\p\left(\A_n\right)
+
\p\left(\left\{b^{3/2}n^{13/18}|R_{n1}|>\nu\right\}\cap\A_n^c\right)\\
&\leq
2\eta/3+b^{3/2}n^{13/18}\nu^{-1}
\E\left[ |R_{n1}|\1_{\A_n^c}\right].
\end{split}
\end{equation}
It suffices to show that there exists $\nu>0$, such that
$b^{3/2}n^{13/18}\nu^{-1}\E\big[|R_{n1}|\1_{\A_n^c}\big]\leq\eta/3$,
for all~$n$ sufficiently large.
Write
\[
\begin{split}
w(u)
&=
\frac{1}{\Phi_n(u;\hat{\beta}_n)}
\left(
a_{n,x}(\hat{A}_n(u))\Phi_n(\hat{A}_n(u);\hat{\beta}_n)-a_{n,x}(u)\Phi_n(u;\hat{\beta}_n)
\right)\\
&=
\frac{a_{n,x}(\hat{A}_n(u))}{\Phi_n(u;\hat{\beta}_n)}
\left(
\Phi_n(\hat{A}_n(u);\hat{\beta}_n)-\Phi(\hat{A}_n(u);\beta_0)
\right)
+
\frac{a_{n,x}(u)}{\Phi_n(u;\hat{\beta}_n)}
\left(
\Phi(u;\beta_0)-\Phi_n(u;\hat{\beta}_n)
\right)\\
&\qquad+
\frac{1}{\Phi_n(u;\hat{\beta}_n)}
\left(
a_{n,x}(\hat{A}_n(u))\Phi(\hat{A}_n(u);\beta_0)-a_{n,x}(u)\Phi(u;\beta_0)
\right).
\end{split}
\]
We will argue that the function $W_n=b^{2}n^{2/9}w$ is uniformly bounded and of bounded variation.
Because of~\eqref{eqn:E_n5},
$n^{1/3}(\Phi_n(\hat{A}_n(u);\hat{\beta}_n)-\Phi(\hat{A}_n(u);\beta_0))$
and $n^{1/3}(\Phi(u;\beta_0)-\Phi_n(u;\hat{\beta}_n))$
are uniformly bounded.
Moreover, they are of bounded variation, as being the difference of two monotone functions.
Similarly, $1/\Phi_n(u;\hat{\beta}_n)$ is of bounded variation and on the event $E_n$ it is also uniformly bounded.
Furthermore, by the definition of $a_{n,x}$, we have
\[
a_{n,x}(\hat{A}_n(u))\Phi(\hat{A}_n(u);\beta_0)-a_{n,x}(u)\Phi(u;\beta_0)
=
k_b\left(x-\hat{A}_n(u)\right)-k_b(x-u).
\]
This is a function of bounded variation, such that multiplied by $b^{2}n^{2/9}$ it is uniformly bounded on the event $\A_n^c$,
because using~\eqref{eqn:lambda}, we obtain
\begin{equation}
\label{eq:bound Rn1}
\begin{split}
\left|k_b\left(x-\hat{A}_n(u)\right)-k_b(x-u)\right|
&\leq
b^{-2}
|\hat{A}_n(u)-u|\sup_{x\in[-1,1]}|k'(x)|\\
&\leq
2Kb^{-2}|\tilde{\lambda}_n(u)-\lambda_0(u)|\sup_{x\in[-1,1]}|k'(x)|\\
&\leq
b^{-2}n^{-2/9}2KK_2\sup_{x\in[-1,1]}|k'(x)|.
\end{split}
\end{equation}
Finally, $ba_{n,x}(u)=bk_b(x-u)/\Phi(u;\beta_0)$
is also a function of bounded variation, as being the product of a function of bounded variation $bk_b(x-u)$ with
the monotone function $1/\Phi(u;\beta_0)$, and it is uniformly bounded.
Then, since $ba_{n,x}(\hat{A}_n(u))$ is the composition of an
increasing function with a function of bounded variation that is uniformly bounded,
it is also a function of bounded variation and uniformly bounded.
As a result, being the sum and product of functions of bounded variation that are uniformly bounded,
$W_n=b^{2}n^{2/9}w$ belongs to the class~$\B_{\tilde{K}}$ of functions of bounded variation, uniformly bounded by some constant $\tilde{K}$.
Consequently, it holds
\[
\begin{split}
R_{n1}
&=
n^{-1/2}
\1_{E_n}\int \delta\1_{[x-b,x+b]}(u)w(u)\,\mathrm{d}\sqrt{n}(\mathbb{P}_n-\mathbb{P})(u,\delta,z)\\
&=
b^{-2}n^{-13/18}
\1_{E_n}\int \delta\1_{[x-b,x+b]}(u)W_n(u)\,\mathrm{d}\sqrt{n}(\mathbb{P}_n-\mathbb{P})(u,\delta,z).
\end{split}
\]
Let $\B_{\tilde K}$ be the class of functions of bounded variation on $[0,M]$, that are uniformly bounded by $\tilde{K}>0$,
and let $\G_n=\{\zeta_{B,n}:B\in\B_{\tilde K}\}$, where $\zeta_{B,n}(u,\delta)=\delta\1_{[x-b,x+b]}(u)B(u)$.
Then, $\delta\1_{[x-b,x+b]}W_n$ is a member of the class $\G_n$, which has envelope $F_n(u,\delta)=\tilde{K}\delta\1_{[x-b,x+b]}(u)$.
Furthermore, if $J(\delta,\G_n)$ is the corresponding entropy-integral (see Section~2.14 in~\cite{VW96}),
then according to Lemma~\ref{lem:covering number} in~\cite{LopuhaaMustaSM2016},
$J(\delta,\G_n)\leq \int_0^\delta \sqrt{1+C/\epsilon}\,\mathrm{d}\epsilon$, for some $C>0$.
Consequently, together with Theorem 2.14.1 in~\cite{VW96}, we obtain that
\[
\begin{split}
\E\left[
|R_{n1}|\1_{\A_n^c}
\right]
&\leq
b^{-2}n^{-13/18}
\E\sup_{\zeta\in\G_n}\left|\int \zeta_{B,n}(u,\delta,z)\,\mathrm{d}\sqrt{n}(\p_n-\p)(u,\delta,z)\right|	\\
&\leq
KJ(1,\G_n)\Vert F_n\Vert_{L_2(\p)}b^{-2}n^{-13/18}
\\
&\leq
K'\left(H^{uc}(x+b)-H^{uc}(x-b)\right)^{1/2}b^{-2}n^{-13/18}
\leq
K''
b^{-3/2}n^{-13/18},
\end{split}
\]
because $H^{uc}$ is absolutely continuous.
As a result,
\[
b^{3/2}n^{13/18}
\nu^{-1}
\E\left[|R_{n1}|\1_{\A_n^c}\right]
\leq
\frac{K''}{\nu}\leq \eta/3
\]
for sufficiently large $\nu$.
This proves the first part of~\eqref{eq:bigO Rn}.

We proceed with $R_{n2}$.
Similar to~\eqref{eq:Rn1 Markov},
\begin{equation}
\label{eq:Rn2 Markov}
\p\left(n^{1/2}|R_{n2}|>\nu\right)
\leq
2\eta/3+n^{1/2}\nu^{-1}
\E\left[ |R_{n2}|\1_{\A_n^c}\right].
\end{equation}
and
it suffices to show that there exists $\nu>0$, such that
$n^{1/2}\nu^{-1}\E\big[|R_{n2}|\1_{\A_n^c}\big]\leq\eta/3$,
for all~$n$ sufficiently large.
We write
\[
n^{1/2}R_{n2}
=
\1_{E_n}
\int
\left(
\mathrm{e}^{\hat{\beta}'_nz}r_{1,n}(u)
-
\mathrm{e}^{\beta'_0z}r_{2,n}(u)
\right)\,\mathrm{d}\sqrt{n}(\p_n-\p)(u,\delta,z),
\]
where
\[
\begin{split}
r_{1,n}(u)
&=
\1_{\{u>x-b\}}\int_{x-b}^{u\wedge(x+b)}\frac{a_{n,x}(\hat{A}_n(v))}{\Phi_n(v;\hat{\beta}_n)}\Phi_n(\hat{A}_n(v);\hat{\beta}_n)
\tilde{\lambda}_n(v)\,\mathrm{d}v,\\
r_{2,n}(u)
&=
\1_{\{u>x-b\}}\int_{x-b}^{u\wedge(x+b)}a_{n,x}(v)\,\lambda_0(v)\,\mathrm{d}v,
\end{split}
\]
are both monotone functions, uniformly bounded by some constant $C$ on the event $\A_n^c\cap E_n$.
Let $\M_C$ be the class of monotone functions bounded uniformly by $C>0$
and let $\G_n=\{\zeta_{r,\beta}(u,z):r\in\M_C,\beta\in\R^P,|\beta-\beta_0|\leq \xi_2\}$,
where $\xi_2$ is chosen as in~\eqref{eqn:ef} and $\zeta_{r,\beta}(u,z)=r(u)\mathrm{e}^{\beta'z}$.
Then $\mathrm{e}^{\hat{\beta}'_nz}r_{1,n}(u)$ is a member of the class $\G_n$,
which has envelope
\[
F_n(u,z)=C\exp\left\{\sum_{j=1}^p (\beta_{0,j}-\sigma_n)z_j\vee(\beta_{0,j}+\sigma_n)z_j\right\},
\]
with $\sigma_n=\sqrt{\xi_2n^{-2/3}}$ is the envelope of $\G_n$.
If $J_{[\,]}(\delta,\G_n,L_2(\mathbb{P}))$ is the bracketing integral
(see Section~2.14 in~\cite{VW96}),
then according to Lemma~\ref{lem:bracketing covering number} in~\cite{LopuhaaMustaSM2016},
$J_{[\,]}(\delta,\G_n,L_2(\mathbb{P}))\leq \int_0^\delta \sqrt{1+C/\epsilon}\,\mathrm{d}\epsilon$,
for some $C>0$.
Hence, together with Theorem 2.14.2 in~\cite{VW96}, we obtain
\[
\begin{split}
&
\E
\left[
\left|
\1_{A_n^c\cap E_n}
\int \mathrm{e}^{\hat{\beta}'_nz}r_{1,n}(u)\,\mathrm{d}\sqrt{n}(\p_n-\p)(u,\delta,z)
\right|\right]\\
&\quad\leq
\E\sup_{\zeta\in\G_n}
\left|\int \zeta_{r,\beta}(u,z)\,\mathrm{d}\sqrt{n}(\p_n-\p)(u,\delta,z)\right|	\\
&\quad\leq
K
J_{[\,]}(1,\G_n,L_2(\p))
\Vert F_n\Vert_{L_2(\p)}
\leq
K',
\end{split}
\]
for some $K'>0$.
We conclude that,
\[
\nu^{-1}
\E\left[
\left|
\1_{A_n^c\cap E_n}
\int
\mathrm{e}^{\hat{\beta}'_nz}r_{1,n}(u)\,\mathrm{d}\sqrt{n}(\p_n-\p)(u,\delta,z)
\right|
\right]
\leq
\frac{K'}{\nu}
\leq
\eta/6,
\]
for sufficiently large $\nu$.
In the same way, it can also be proved
\[
\nu^{-1}
\E\left[
\left|
\1_{A_n^c\cap E_n}
\int \mathrm{e}^{\beta_0z}r_{2,n}(u)\,\mathrm{d}\sqrt{n}(\p_n-\p)(u,\delta,z)
\right|\right]
\leq
\frac{K}{\nu}\leq \eta/6
\]
for sufficiently large $\nu$, concluding the proof of~\eqref{eq:Rn2 Markov}
and therefore the second part of~\eqref{eq:bigO Rn}.
\end{proof}

\section*{Supplementary Material}
Supplement to ''Smoothed isotonic estimators of a monotone baseline hazard in the Cox model''.
\begin{itemize}
\item Supplement~\ref{app:SG}: Entropy bounds for the smoothed Grenander-type estimator.
\item Supplement~\ref{sec:SMLE}: Smooth maximum likelihood estimator.
\item Supplement~\ref{sec:bootstrap}: Consistency of the bootstrap.
\end{itemize}

\bibliographystyle{imsart-nameyear}
\bibliography{shapeconstrained-estimation}

\setcounter{equation}{0}
\renewcommand{\theequation}{S\arabic{equation}}
\newpage

\medskip

\newpage

\centerline{\LARGE\bf Smoothed isotonic estimators of a monotone baseline}
\smallskip
\centerline{\LARGE\bf hazard in the Cox model}
\bigskip
\centerline{\LARGE Supplementary Material}
\bigskip
\centerline{\Large Hendrik Paul Lopuha\"a$^\dagger$ and Eni Musta$^\dagger$}
\medskip
\centerline{\it Delft University of Technology$^\dagger$}
\bigskip
\centerline{\date{\today}}

\appendix

\section{Entropy bounds for the smoothed Grenander-type estimator}
\label{app:SG}

\begin{lemma}
\label{lem:covering number}
Let $0<x<M<\tau_H$ and let $\B_{\tilde K}$ be the class of functions of bounded variation on $[0,M]$, that are uniformly bounded by $\tilde{K}>0$.
Let $\G_n=\{\zeta_{B,n}:B\in\B_{\tilde K}\}$, where $\zeta_{B,n}(u,\delta)=\delta\1_{[x-b,x+b]}(u)B(u)$ and $n^{1/5}b\to c>0$.
For $\delta>0$, let
\[
J(\delta,\G_n)
=
\sup_{Q}\int_0^\delta
\sqrt{1+\log N(\epsilon\Vert F_n\Vert_{L_2(Q)},\G_n,L_2(Q))}\,\mathrm{d}\epsilon
\]
where $N(\epsilon,\G_n,L_2(Q))$ is the minimal number of $L_2(Q)$-balls
of radius $\epsilon$ needed to cover $\G_n$,
$F_n=\tilde K\delta\1_{[x-b,x+b]}$ is the envelope of $\G_n$,
and the supremum is taken over measures $Q$ on $[0,M]\times\{0,1\}\times\R^p$,
for which $Q([x-b,x+b]\times\{1\}\times\R^p)>0$.
Then $J(\delta,\G_n)\leq \int_0^\delta \sqrt{1+C/\epsilon}\,\mathrm{d}\epsilon$,
for some $C>0$.
\end{lemma}
\begin{proof}
We first bound the entropy of $\G_n$ with respect to any probability measure $Q$ on $[0,M]\times\{0,1\}\times\R^p$, such that
\begin{equation}
\label{eq:envelope}
\Vert F_n\Vert_{L_2(Q)}^2
=
\tilde K
Q\left([x-b,x+b]\times\{1\}\times\R^p\right)>0.
\end{equation}
Fix such a probability measure  $Q$ and let $\epsilon>0$.
Let $Q'$ be the probability measure on $[0,M]$ defined by
\[
Q'(S)=\frac{Q(S\cap [x-b,x+b]\times\{1\}\times\R^p)}{Q([x-b,x+b]\times\{1\}\times\R^p)},
\qquad
S\in[0,M],
\]
For a given $\epsilon'>0$, select a minimal $\epsilon'$-net $B_1,\dots,B_{N}$ in $\B_{\tilde{K}}$ with respect to $L_2(Q')$,
where $N=N(\epsilon',\B_{\tilde{K}},L_2(Q'))$.
Then, from (2.6) in~\cite{Geer}
\begin{equation}
\label{eqn:entr}
\log N(\epsilon',\B_{\tilde{K}},L_2(Q'))\leq\frac{K}{\epsilon'},
\end{equation}
for some constant $K>0$.
Then, consider the functions $\zeta_{B_1,n},\ldots,\zeta_{B_N,n}$ corresponding to $B_1,\ldots,B_N$.
For any $\zeta_{B,n}\in\G_n$, there exists a $B_i$ in the $\epsilon'$-net
that is closest function to~$B$, i.e.,
$\Vert
B-B_i
\Vert_{L_2(Q')}\leq\epsilon'$.
Then, we find
\begin{equation}
\label{eqn:res5}
\begin{split}
\Vert
\zeta_{B,n}-\zeta_{B_i,n}
\Vert_{L_2(Q)}
&=
Q([x-b,x+b]\times\{1\}\times\R^p)^{1/2}
\Vert\,B-B_i\Vert_{L_2(Q')}\\
&\leq
Q([x-b,x+b]\times\{1\}\times\R^p)^{1/2}\epsilon'.
\end{split}
\end{equation}
Hence, if we take $\epsilon'=\epsilon/Q([x-b,x+b]\times\{1\}\times\R^p)^{1/2}$,
it follows that $\zeta_{B_1,n},\ldots,\zeta_{B_N,n}$ forms an $\epsilon$-net in $\G_n$
with respect to $L_2(Q)$,
and $N(\epsilon,\G_n,L_2(Q))\leq N(\epsilon',\B_{\tilde{K}},L_2(Q'))$.
Using~\eqref{eqn:entr}, this implies that
\[
\log N(\epsilon,\G_n,L_2(Q))
\leq
\frac{K}{\epsilon'}
=
\frac{K\cdot Q([x-b,x+b]\times\{1\}\times\R^p)^{1/2}}{\epsilon}
\]
where $K$ does not depend on $Q$, and according to~\eqref{eq:envelope},
\[
J(\delta,\G_n)
=
\sup_{Q}\int_0^\delta
\sqrt{1+\log N(\epsilon\Vert F_n\Vert_{L_2(Q)},\G_n,L_2(Q))}\,\mathrm{d}\epsilon
\leq
\int_0^\delta
\sqrt{1+\frac{C}{\epsilon\tilde K^{1/2}}}\,\mathrm{d}\epsilon,
\]
for some $C>0$.
\end{proof}

\begin{lemma}
\label{lem:bracketing covering number}
Let $\M_C$ be the class of monotone functions bounded uniformly by $C>0$.
Let $\G_n=\{\zeta_{r,\beta}(u,z):r\in\M_C,\beta\in\R^P,|\beta-\beta_0|\leq \xi_2\}$,
where $\xi_2$ is chosen as in~\eqref{eqn:ef} and $\zeta_{r,\beta}(u,z)=r(u)\mathrm{e}^{\beta'z}$.
For $\delta>0$, let
\[
J_{[\,]}(\delta,\G_n,L_2(\mathbb{P}))
=
\int_0^\delta
\sqrt{1+\log N_{[\,]}(\epsilon\Vert F_n\Vert_{L_2(\mathbb{P})},\G_n,L_2(\mathbb{P}))}\,\mathrm{d}\epsilon.
\]
where $N_{[\,]}(\epsilon,\G_n,L_2(\mathbb{P}))$ is the bracketing number
and
\[
F_n(u,z)=C\exp\left\{\sum_{j=1}^p (\beta_{0,j}-\sigma_n)z_j\vee(\beta_{0,j}+\sigma_n)z_j\right\},
\]
with $\sigma_n=\sqrt{\xi_2n^{-2/3}}$ is the envelope of $\G_n$.
Then $J_{[\,]}(\delta,\G_n,L_2(\mathbb{P}))\leq \int_0^\delta \sqrt{1+C/\epsilon}\,\mathrm{d}\epsilon$,
for some $C>0$.
\end{lemma}
\begin{proof}
The entropy with bracketing for the class of bounded monotone functions on $\R$ satisfies
\begin{equation}
\label{eq:entropy MC}
\log N_{[\,]}(\gamma,\M_C,L_p(Q))\leq \frac{K}{\gamma},
\end{equation}
for every probability measure $Q$, $\gamma>0$ and $p\geq 1$ (e.g., see Theorem 2.7.5 in~\cite{VW96}).
Define the probability measure $Q$ on $\R$ by
\[
Q(S)=\int \1_S(u)\,\rm{d}\p(u,\delta,z),
\qquad
\text{for all }S\subseteq\R,
\]
and fix $\epsilon>0$.
For a given $\gamma>0$, take a net of $\gamma$-brackets $\{(l_1,L_1)\,\dots,(l_N,L_N)\}$ in $\M_C$
with respect to~$L_4(Q)$,
where $N=N_{[\,]}(\gamma,\M_C,L_4(Q))$.
For every $j=1,\ldots,p$,
divide the interval $[\beta_{0,j}-\sigma_n,\beta_{0,j}+\sigma_n]$ in subintervals of length $\gamma$,
i.e.,
\[
\left[a_{k_j},b_{k_j}\right]
=
\left[\beta_{0,j}-\sigma_n+k_j\gamma,\beta_{0,j}-\sigma_n+(k_j+1)\gamma\right],
\]
for $k_j=0,\dots,\bar{N}-1$,
where $\bar{N}=2\sigma_n/\gamma$.
Then, construct brackets $(l_{i,k_1,\ldots,k_p},L_{i,k_1,\ldots,k_p})$,
for $i=1,\dots,N$ and $k_j=0,\dots,\bar{N}-1$, for $j=1,\ldots,p$, in the following way:
\[
\begin{split}
l_{i,k_1,\ldots,k_p}(u,z)
&=
l_i(u)
\exp\left\{\sum_{j=1}^p a_{k_j}z_j\wedge b_{k_j}z_j\right\},\\
L_{i,k_1,\ldots,k_p}(u,z)
&=
L_i(u)
\exp\left\{\sum_{j=1}^p a_{k_j}z_j\vee b_{k_j}z_j\right\}.
\end{split}
\]
By construction, using~\eqref{eq:entropy MC}, the number of these brackets is
\begin{equation}
\label{eq:covering bracketing number}
N\cdot\bar{N}^p
\leq
\exp\{K/\gamma\}\cdot\left(\frac{2\sigma_n}{\gamma}\right)^p
\leq
\exp\{K/\gamma\}\cdot\exp\left\{2p\sqrt{\xi_2}/\gamma\right\}
=
\mathrm{e}^{K_1/\gamma},
\end{equation}
for some $K_1>0$ independent of $n\geq 1$.
For any $\zeta_{r,\beta}\in\G_n$, there exist a $\gamma$-bracket $(l_i,L_i)$, such that $r\in[l_i,L_i]$,
and intervals $[a_{k_j},b_{k_j}]$, such that $\beta_j\in[a_{k_j},b_{k_j}]$,
for all $j=1,\ldots,p$.
It follows that there exists a bracket $(l_{i,k_1,\ldots,k_p},L_{i,k_1,\ldots,k_p})$, such that
$\zeta_{r,\beta}\in[l_{i,k_1,\ldots,k_p},L_{i,k_1,\ldots,k_p}]$.
Moreover, for each $i=1,\dots,n$ and $k_j=0,\ldots,\bar N-1$, for $j=1,\ldots,p$,
\[
\begin{split}
&
\Vert l_{i,k_1,\ldots,k_p}-L_{i,k_1,\ldots,k_p}\Vert_{L_2(\p)}\\
&=
\left(\int \left|l_{i,k_1,\ldots,k_p}(u,z)-L_{i,k_1,\ldots,k_p}(u,z)\right|^2\,\mathrm{d}\p(u,\delta,z) \right)^{1/2}\\
&\leq
\left(\int
\exp\left\{
2\sum_{j=1}^p a_{k_j}z_j\wedge b_{k_j}z_j
\right\}
\left|l_i(u)-L_i(u)\right|^2\,\mathrm{d}\p(u,\delta,z) \right)^{1/2}\\
&\qquad+
C\left(\int
\left|
\exp\left\{
\sum_{j=1}^p a_{k_j}z_j\wedge b_{k_j}z_j
\right\}
-
\exp\left\{
\sum_{j=1}^p a_{k_j}z_j\vee b_{k_j}z_j
\right\}
\right|^2\,\mathrm{d}\p(u,\delta,z)
\right)^{1/2}\\
&\leq
\left(\int
\exp\left\{
4\sum_{j=1}^p a_{k_j}z_j\wedge b_{k_j}z_j
\right\}
\,\mathrm{d}\p(u,\delta,z) \right)^{1/4}
\left(\int
|l_i(u)-L_i(u)|^4\,\mathrm{d}\p(u,\delta,z) \right)^{1/4}\\
&\qquad+
C\left(\int
\left|
\exp\left\{
\sum_{j=1}^p a_{k_j}z_j\wedge b_{k_j}z_j
\right\}
-
\exp\left\{
\sum_{j=1}^p a_{k_j}z_j\vee b_{k_j}z_j
\right\}
\right|^2\,\mathrm{d}\p(u,\delta,z)
\right)^{1/2}.
\end{split}
\]
For the first term on the right hand side, we use that for all $j=1,\ldots,p$,
\begin{equation}
\label{eq:bound aj bj}
\left|a_{k_j}-\beta_{0,j}\right|\leq \sigma_n
\quad\text{and}\quad
\left|b_{k_j}-\beta_{0,j}\right|\leq \sigma_n,
\end{equation}
so that, according to~(A3), there exist a $K_1>0$, such that for $n$ sufficiently large
\[
\int
\exp\left\{
4\sum_{j=1}^p a_{k_j}z_j\wedge b_{k_j}z_j
\right\}
\,\mathrm{d}\p(u,\delta,z)
\leq
\sup_{|\beta-\beta_0|\leq p\sigma_n}
\E
\left[
\mathrm{e}^{4\beta'Z}
\right]
\leq
K_1.
\]
For the second term we have
\[
\left(\int
|l_i(u)-L_i(u)|^4\,\mathrm{d}\p(u,\delta,z) \right)^{1/4}
=
\|l_i-L_i\|_{L_4(Q)}\leq \gamma,
\]
by construction,
and by the mean value theorem, the third term on the right hand side
can be bounded by
\[
C\left(\int \left(\gamma\sum_{i=1}^p |z_i|\right)^2
\mathrm{e}^{2\theta_z}\mathrm{d}\p(u,\delta,z)\right)^{1/2}
\leq
C'\gamma
\left(
\int
|z|^2
\mathrm{e}^{2\theta_z}\mathrm{d}\p(u,\delta,z)\right)^{1/2},
\]
for some $\sum_{j=1}^p a_{k_j}z_j\wedge b_{k_j}z_j
\leq
\theta_z
\leq
\sum_{j=1}^p a_{k_j}z_j\vee b_{k_j}z_j$.
Consequently, in view of~\eqref{eq:bound aj bj},
we find
\[
\Vert l_{i,k_1,\ldots,k_p}-L_{i,k_1,\ldots,k_p}\Vert_{L_2(\p)}
\leq
\gamma
K_1^{1/4}
+
C'\gamma
\left(\sup_{|\beta-\beta_0|\leq p\sigma_n}
\E\left[|Z|^2\mathrm{e}^{2\beta'Z}\right]
\right)^{1/2}
\leq
K'\gamma,
\]
for some $K'>0$, using~(A2).
Hence, if we take $\gamma=\epsilon/K'$, then
$\{(l_{i,k_1,\ldots,k_p},L_{i,k_1,\ldots,k_p})\}$,
for $i=1,\ldots,n$ and $k_j=0,\ldots,\bar{N}-1$, for $j=1,\ldots,p$, forms
a net of $\epsilon$-brackets, and according to~\eqref{eq:covering bracketing number},
there exists a $K>0$, such that
\[
\log N_{[\,]}(\epsilon,\G_n,L_2(\p))\leq \frac{K}{\epsilon}.
\]
As a result,
\[
J_{[\,]}(1,\G_n,L_2(\p))
=
\int_0^1\sqrt{1+\log N_{[\,]}(\epsilon\,\Vert F\Vert_{L_2(\p)},\G_n,L_2(\p))}\,\mathrm{d}\epsilon
=
\int_0^1\sqrt{1+\frac{K}{\epsilon}}\,\mathrm{d}\epsilon,
\]
for some $K>0$.
\end{proof}

\section{Smooth maximum likelihood estimator}
\label{sec:SMLE}

To derive the pointwise asymptotic distribution of $\hat{\lambda}_n^{SM}$, we follow the same approach as the one used for $\tilde\lambda_n^{SG}$.
We will go through the same line of reasoning as used to obtain Theorem~\ref{theo:distr}.
However, large parts of the proof are very similar, if not the same.
We briefly sketch the main steps.

Instead of $\tilde{\Lambda}_n$, we use
\[
\hat{\Lambda}_n(x)=\int_0^x\hat{\lambda}_n(u)\,\mathrm{d}u,
\]
where $\hat{\lambda}_n$ is the MLE for $\lambda_0$.
First, the analogue of Lemma~\ref{le:2-1} still holds.
\begin{lemma}
\label{le:2-1 SM}
Suppose that (A1)--(A2) hold.
Let $a_{n,x}$ be defined by~\eqref{eqn:an} and let $\hat\beta_n$ be the partial maximum likelihood estimator for $\beta_0$.
Define
\begin{equation}
\label{eqn:thetan SM}
\theta_{n,x}(u,\delta,z)
=
\1_{E_n}
\left\{
\delta
a_{n,x}(u)-\mathrm{e}^{\hat{\beta}'_n\,z}\int_0^u a_{n,x}(v)\,\mathrm{d}\hat{\Lambda}_n(v)
\right\},
\end{equation}
Then, there exists an event $E_n$, with $\1_{E_n}\to1$ in probability, such that
\[
\int\theta_{n,x}(u,\delta,z)\,\mathrm{d}\p(u,\delta,z)=-\1_{E_n}\int k_b(x-u)\,\mathrm{d}(\hat{\Lambda}_n-\Lambda_0)(u)+O_p(n^{-1/2}).
\]
\end{lemma}
\begin{proof}
We modify the definition of the event $E_n$ from Lemma~\ref{le:2-1} as follows.
The events~$E_{n,1}$, $E_{n,2}$, $E_{n,4}$,  and $E_{n,5}$, from~\eqref{eqn:E1} and~\eqref{eqn:ef} remain the same.
Replace $E_{n,3}$ in~\eqref{eqn:ef} by
\begin{equation}
\label{eq:def En3}
E_{n,3}
=
\left\{
\sup_{x\in[T_{(1)},T_{(n)}]}|V_n(x)-H^{uc}(x)|<\xi_3
\right\},
\end{equation}
and let
\begin{equation}
\label{eq:def En6}
E_{n,6}
=
\left\{
T_{(1)}\leq\epsilon
\right\},
\end{equation}
for some $\epsilon>0$ and $\xi_3>0$.
Then $\p(E_{n,6})\to 1$, and also $\p(E_{n,3})\to 1$ according to Lemma~5 in~\cite{LopuhaaNane2013}.
As before $E_n=\bigcap_{i=1}^6E_{n,i}$.
Similar to the proof of Lemma~\ref{le:2-1}, we obtain
\[
\begin{split}
&
\int\theta_{n,x}(u,\delta,z)\,\mathrm{d}\p(u,\delta,z)\\
&=
\1_{E_n}
\left\{
\int a_{n,x}(u)\,\mathrm{d}H^{uc}(u)-\int \mathrm{e}^{\hat{\beta}'_nz}
\int_{v=0}^u a_{n,x}(v)\,\mathrm{d}\hat{\Lambda}_n(v)\,\mathrm{d}\p(u,\delta,z)
\right\}\\
&=
\1_{E_n}
\left\{
-\int k_b(x-u)\,\mathrm{d}(\hat{\Lambda}_n-\Lambda_0)(u)+\int k_b(x-u)\,\left(1-\frac{\Phi(u;\hat{\beta}_n)}{\Phi(u;\beta_0)}\right)\,\mathrm{d}\hat{\Lambda}_n(u)
\right\},
\end{split}
\]
and
\[
\1_{E_n}\int k_b(x-u)\,\left|1-\frac{\Phi(u;\hat{\beta}_n)}{\Phi(u;\beta_0)}\right|\,\mathrm{d}\hat{\Lambda}_n(u)=O_p(n^{-1/2}).
\]
which proves the lemma.
\end{proof}
Next, we slightly change the definition of $\overline{\Psi}_{n,x}$ from~\eqref{def:Psibar}, i.e.,
\begin{equation}
\label{def:Psibar SM}
\overline{\Psi}_{n,x}(u)
=
\overline{a}_{n,x}(u)\1_{E_n}=a_{n,x}(\hat{A}_n(u))\1_{E_n}
\end{equation}
where $E_n$ is the event from Lemma~\ref{le:2-1 SM} and $\hat{A}_n$, as defined in~\eqref{eq:An}, is now taken constant on $[\tau_i,\tau_{i+1})$ and consider
\begin{equation}
\label{eq:def theta bar SM}
\bar{\theta}_{n,x}(u,\delta,z)
=
\delta
\overline{\Psi}_{n,x}(u)-\mathrm{e}^{\hat{\beta}'_nz}
\int_0^u \overline{\Psi}_{n,x}(v)\,\mathrm{d}\hat{\Lambda}_n(v).
\end{equation}
With this $\bar{\theta}_{n,x}$, we have the same property as in Lemma~\ref{le:2-2}.
\begin{lemma}
\label{le:2-2 SM}
Let $\bar{\theta}_{n,x}$ be defined in~\eqref{eq:def theta bar SM}.
Then
\[
\int\bar{\theta}_{n,x}(u,\delta,z)\,\mathrm{d}\p_n(u,\delta,z)=0.
\]
\end{lemma}
\begin{proof}
Similar to the proof of Lemma~\ref{le:2-2}, we have
\[
\begin{split}
&
\int\bar{\theta}_{n,x}(u,\delta,z)\,\mathrm{d}\p_n(u,\delta,z)\\
&=
\1_{E_n}
\sum_{i=0}^{m}
\overline{a}_{n,x}(\tau_i)
\left\{
\int
\1_{[\tau_i,\tau_{i+1})}(u)\delta\,\mathrm{d}\p_n(u,\delta,z)
-
\int_{\tau_i}^{\tau_{i+1}}\Phi_n(v;\hat{\beta}_n)\,\mathrm{d}\hat{\Lambda}_n(v)
\right\}\\
&=
\1_{E_n}\sum_{i=0}^{m} \overline{a}_{n,x}(\tau_i)
\left\{
V_n(\tau_{i+1})-V_n(\tau_i)-\hat{\lambda}_n(\tau_i)
\left(\hat{W}_n(\tau_{i+1})-\hat{W}_n(\tau_i)\right)
\right\}\\
&=0,
\end{split}
\]
The last equality follows from the characterization of the maximum likelihood estimator.
\end{proof}
Furthermore, for  $\bar{\theta}_{n,x}$ defined in~\eqref{eq:def theta bar SM}, we also have
\begin{equation}
\label{eq:Lemma33 SM}
\int
\left\{
\bar{\theta}_{n,x}(u,\delta,z)-\theta_{n,x}(u,\delta,z)
\right\}\,\mathrm{d}\p(u,\delta,z)
=
O_p(b^{-1}n^{-2/3}),
\end{equation}
see Lemma~\ref{le:result4 SM},
and
\begin{equation}
\label{eq:Lemma34 SM}
\int
\left\{
\bar{\theta}_{n,x}(u,\delta,z)-\eta_{n,x}(u,\delta,z)
\right\}
\,\mathrm{d}(\p_n-\p)(u,\delta,z)
=
O_p(b^{-3/2}n^{-13/18})+O_p(n^{-1/2}),
\end{equation}
see Lemma~\ref{le:result5 SM},
where $\eta_{n,x}$ is defined similar to~\eqref{eqn:eta}, but with $E_n$ taken from Lemma~\ref{le:2-1 SM}.
Similar to the proof of Lemma~\ref{le:result4},
the proof of Lemma~\ref{le:result4 SM} is quite technical and involves
bounds on the tail probabilities of the inverse process corresponding to $\hat\lambda_n$
(see Lemma~\ref{le:invSM}), used to obtain the analogue of~\eqref{eq:E sup}
(see Lemma~\ref{le:lambdaSM}).
It is defined by
\[
\hat{U}_n(a)
=
\argmin_{x\in[T_{(1)},T_{(n)}]}
\left\{V_n(x)-a\hat{W}_n(x)\right\}
\]
where
\[
\begin{split}
\hat{W}_n(x)
&=
\int
\left(
\mathrm{e}^{\hat{\beta}'_nz}
\int_{T_{(1)}}^x\1_{\{u\geq s\}}\,\mathrm{d}s \right)\,\mathrm{d}\p_n(u,\delta,z),\quad x\geq T_{(1)},\\
V_n(x)
&=
\int\delta\1_{\{u<x\}}\,\mathrm{d}\p_n(u,\delta,z),
\end{split}
\]
with $\hat\beta_n$ being the partial maximum likelihood estimator (see Lemma~1 in~\cite{LopuhaaNane2013})
and it satisfies the switching relation
$\hat{\lambda}_n(x)\leq a$ if and only if $\hat{U}_n(a)\geq x$.
Let $U$ be the inverse of $\lambda_0$ on $[\lambda_0(\epsilon),\lambda_0(M)]$, for some $0<\epsilon<M<\tau_H$, i.e.,
\begin{equation}
\label{def:U SM}
U(a)
=
\begin{cases}
\epsilon & a< \lambda_0(\epsilon);\\
\lambda_0^{-1}(a) & a\in[\lambda_0(\epsilon),\lambda_0(M)];\\
M & a>\lambda_0(M).
\end{cases}
\end{equation}
Similar to the bounding the tail probabilities of $\tilde U_n(a)$, in order to bound
the tail probabilities of $\hat U_n(a)$ we first introduce a suitable martingale that will approximate the process $V_n(t)-a\hat W_n(t)$.
\begin{lemma}
\label{lem:mart1}
Suppose that (A1)--(A2) hold.
Define
\begin{equation}
\label{eq:def Bn bar}
\bar{\w}_n(t)=V_n(t)-\int_{0}^t \Phi_n(s;\beta_0)\lambda_0(s)\,\mathrm{d}s.
\end{equation}
The process $\{(\bar{\w}_n(t),\F^n_t):0\leq t <\tau_H\}$ is a square integrable martingale with mean zero and predictable variation process
\[
\langle\bar{\w}_n\rangle(t)=\frac{1}{n}\int_0^t\Phi_n(s;\beta_0)\,\mathrm{d}\Lambda_0(s).
\]
\end{lemma}
\begin{proof}
Note that
\[
\bar{\w}_n(t)
=
\frac{1}{n}\m_n(t)-\frac{1}{n}\sum_{i=1}^n\1_{\{T_i=t\}}\Delta_i
\]
where $\m_n$ is defined in~\eqref{eq:def Mn}.
Since $H^{uc}$ is absolutely continuous, $\1_{\{T_i=t\}}\Delta_i=0$ a.s., which means that $\bar{\w}_n=\m_n$ a.s..
Hence $\bar{\w}_n$ is a mean zero martingale
and has the same predictable variation as $n^{-1}\m_n$.
\end{proof}

\begin{lemma}
\label{le:martSM}
Suppose that (A1)--(A2) hold.
There exists a constant $C>0$ such that, for all $x>0$ and $t\in[0,\tau_H]$,
\[
\E\left[\sup_{u\in[0,\tau_H],|t-u|\leq x}
\left(\bar{\w}_n(u)-\bar{\w}_n(t)\right)^2\right]\leq \frac{Cx}{n}.
\]
\end{lemma}
\begin{proof}
First, consider the case $t\leq u\leq t+x$.
Then, by Doob's inequality, we have
\[
\begin{split}
\E\left[\sup_{u\in[0,\tau_H],t\leq u\leq t+x}(\bar{\w}_n(u)-\bar{\w}_n(t))^2\right]
&\leq
4\E\left[\left(\bar{\w}_n\big((t+x)\wedge \tau_H\big)-\bar{\w}_n(t)\right)^2\right]\\
&=
4\E\left[\left(\bar{\w}_n\big((t+x)\wedge \tau_H\big)\big)^2-\big(\bar{\w}_n(t)\right)^2\right]\\
&=
\frac{4}{n}\E\left[\int_{t}^{(t+x)\wedge \tau_H}\Phi_n(s;\beta_0)\lambda_0(s)\,\mathrm{d}s \right]\\
&\leq
\frac{4\lambda_0(\tau_H)x}{n}\E\left[\Phi_n(0;\beta_0)\right]\\
&=
\frac{4\lambda_0(\tau_H)x}{n}\frac{1}{n}\sum_{i=1}^n\E\left[\mathrm{e}^{\beta'_0Z_i}\right]
\leq \frac{Kx}{n}
\end{split}
\]
for some $K>0$, using~(A2).
For the case $t-x\leq u\leq t$, we can write
\[
\begin{split}
&
\E\left[\sup_{u\in[0,\tau_H],t-x\leq u\leq t}
\big(\bar\w_n(u)-\bar\w_n(t)\big)^2\right]
=
\E\left[\sup_{0\vee(t-x)\leq u\leq t}
\big(\bar\w_n(u)-\bar\w_n(t)\big)^2\right]\\
&\quad\leq
2\E\left[\big(\bar\w_n(t)-\bar\w_n(0\vee(t-x))\big)^2\right]
+
2\E\left[\sup_{0\vee(t-x)\leq u<t}
\big(\bar\w_n(u)-\bar\w_n(0\vee(t-x))\big)^2\right].
\end{split}
\]
Then similar, the right hand side is bounded by
\[
\begin{split}
&
2\E\left[\big(\bar\w_n(t)-\bar\w_n(0\vee(t-x))\big)^2\right]
+
8\E\left[\big(\bar\w_n(t)-\bar\w_n(0\vee(t-x))\big)^2\right]\\
&=
10\E\left[\bar\w_n(t)^2-\bar\w_n(0\vee(t-x))^2\right]
=
\frac{10}n
\E\left[\int_{0\vee(t-x)}^{t}
\Phi_n(s;\beta_0)\lambda_0(s)\,\mathrm{d}s  \right]\\
&\leq
\frac{10\lambda_0(\tau_H)x}{n}
\E\left[\Phi_n(0;\beta_0) \right]
\leq
\frac{Cx}{n},
\end{split}
\]
for some $C>0$.
This concludes the proof.
\end{proof}

\begin{lemma}
\label{le:invSM}
Suppose that (A1)--(A2) hold.
Let $0<\epsilon<M<\tau_H$ and let $\hat U_n(a)$ and $U$ be defined in~\eqref{eq:def Un SM} and~\eqref{def:U SM}.
Suppose that $\lambda_0'$ is uniformly bounded below by a strictly positive constant.
Then, there exists an event~$E_n$, such that $\1_{E_n}\to 1$ in probability, and
a constant $K$ such that, for every $a\geq 0$ and $x> 0$,
\begin{equation}
\label{eqn:invSM}
\p\left(
\left\{|\hat{U}_n(a)-U(a)|\geq x\right\}
\cap
E_n
\cap
\left\{\epsilon\leq\hat{U}_n(a)\leq M\right\}
\right)
\leq
\frac{K}{nx^3}.
\end{equation}
\end{lemma}
\begin{proof}
Similar to the proof of Lemma~\ref{le:inv}, we start by writing
\begin{equation}
\label{eqn:mainSM}
\begin{split}
&
\p\left(
\left\{|\hat{U}_n(a)-U(a)|\geq x\right\}
\cap
E_n
\cap
\left\{\epsilon\leq\hat{U}_n(a)\leq M\right\}
\right)\\
&=
\p\left(\left\{U(a)+x\leq \hat{U}_n(a)\leq M\right\}\cap E_n\right)
+
\p\left(\left\{\epsilon\leq \hat{U}_n(a)\leq U(a)-x\right\}\cap E_n\right).
\end{split}
\end{equation}
The first probability is zero if $U(a)+x>M$.
Otherwise, if $U(a)+x\leq M$, then $x\leq M$ and we get
\[
\begin{split}
&
\p\left(\left\{U(a)+x\leq \hat{U}_n(a)\leq M \right\}\cap E_n\right)\\
&\leq
\p\left(\left\{
\inf_{y\in[U(a)+x,M]}
\left(
V_n(y)-a\hat{W}_n(y)-V_n(U(a))+a\hat{W}_n(U(a))
\right)\leq 0\right\}\cap E_n\right).
\end{split}
\]
Define
\begin{equation}
\label{eq:def Rn bar}
\bar{R}_n(t)=a\int_{0}^{t}\left(\Phi_n(s;\beta_0)-\Phi_n(s;\hat{\beta}_n)\right)\,\mathrm{d}s.
\end{equation}
Then, for $T_{(1)}<U(a)<y$,
\[
\begin{split}
&
\bar{\w}_n(y)-\bar{\w}_n(U(a))+\bar{R}_n(y)-\bar{R}_n(U(a))\\
&=
V_n(y)-V_n(U(a))-\int_{U(a)}^y \Phi_n(s;\beta_0)\lambda_0(s)\,\mathrm{d}s
+
a\int_{U(a)}^{y}\left(\Phi_n(s;\beta_0)-\Phi_n(s;\hat{\beta}_n)\right)\,\mathrm{d}s\\
&=
V_n(y)-V_n(U(a))-a\int_{U(a)}^{y}\Phi_n(s;\hat{\beta}_n)\,\mathrm{d}s
-\int_{U(a)}^y \Phi_n(s;\beta_0)\lambda_0(s)\,\mathrm{d}s
+
a\int_{U(a)}^{y}\Phi_n(s;\beta_0)\,\mathrm{d}s\\
&=
V_n(y)-a\hat{W}_n(y)-V_n(U(a))+a\hat{W}_n(U(a))
-
\int_{U(a)}^{y}\Phi_n(s;\beta_0)(\lambda_0(s)-a)\,\mathrm{d}s.
\end{split}
\]
On the event $E_n$, by Taylor expansion we find that
\[
\begin{split}
\int_{U(a)}^{y}\Phi_n(s;\beta_0)(\lambda_0(s)-a)\,\mathrm{d}s
&=
\int_{U(a)}^{y}\Phi_n(s;\beta_0)(\lambda_0(s)-\lambda_0(U(a)))\,\mathrm{d}s\\
&=
\int_{U(a)}^{y}\Phi_n(s;\beta_0)
\left(
\lambda_0'(\xi_s)(s-U(a))
\right)\,\mathrm{d}s\\
&\geq
\inf_{t\in[0,\tau_H)}\lambda_0'(t)
\left(\Phi(M;\beta_0)-\xi_4n^{-1/3}\right)
\frac12(y-U(a))^2\\
&\geq
c(y-U(a))^2
\end{split}
\]
for some $c>0$.
Similar to the proof of Lemma~\ref{le:inv}, it follows that
\[
\begin{split}
&
\p\left(\left\{
\inf_{y\in[U(a)+x,M]}
\left(
V_n(y)-a\hat{W}_n(y)-V_n(U(a))+a\hat{W}_n(U(a))
\right)\leq 0\right\}\cap E_n\right)\\
&\leq
\p\left(\left\{
\inf_{y\in[U(a)+x,M]}
\left(
\bar{\w}_n(y)-\bar{\w}_n(U(a))+\bar{R}_n(y)-\bar{R}_n(U(a))+c(y-U(a))^2
\right)\leq 0\right\}\cap E_n\right).
\end{split}
\]
Hence, as before
\[
\begin{split}
&
\p\left(\left\{U(a)+x\leq \hat{U}_n(a)\leq M \right\}\cap E_n\right)\\
&\leq
\sum_{k=0}^i
\p\left(\left\{
\sup_{I_k}
\Bigg(
\left|
\bar{\w}_n(y)-\bar{\w}_n(U(a))
\right|+\left|\bar{R}_n(y)-\bar{R}_n(U(a))\right|
\Bigg)
\geq
cx^22^{2k}
\right\}\cap E_n\right)\\
\end{split}
\]
where the supremum runs over $y\leq M$, such that $y-U(a)\in[x2^k,x2^{k+1})$.
With Markov, we can bound this probability by
\begin{equation}
\label{eqn:1SM}
\begin{split}
&
4\sum_{k=0}^i
\left(c^2x^42^{4k}\right)^{-1}
\E\left[\sup_{y\leq M,y-U(a)\in[x2^k,x2^{k+1})}\big|\bar\w_n(y)-\bar\w_n(U(a))\big|^2\right]\\
&\quad+
8\sum_{k=0}^i
\left(
c^3x^62^{6k}
\right)^{-1}
\E\left[\sup_{y<M,y-U(a)\in[x2^k,x2^{k+1})}\1_{E_n}
\big|\bar{R}_n(y)-\bar{R}_n(U(a))\big|^3\right].
\end{split}
\end{equation}
As in the proof of Lemma~\ref{le:inv}, we will bound both expectations separately.
We have
\[
\begin{split}
&
\E\left[\sup_{y<M,\,y-U(a)\in[x2^k,x2^{k+1})}\1_{E_n}\,\big|\bar{R}_n(y)-\bar{R}_n(U(a))\big|^3\right]\\
&\leq
\E\left[\1_{E_n}
\left(
\int_{U(a)}^{(U(a)+x2^{k+1})\wedge M}
a\left|\Phi_n(s;\beta_0)-\Phi_n(s;\hat{\beta}_n)\right|\,\mathrm{d}s
\right)^3
\right]\\
&\leq
x^32^{3(k+1)}\lambda_0(M)^3
\E\left[\1_{E_n}\,\sup_{s\in[0,M]}\left|\Phi_n(s;\beta_0)-\Phi_n(s;\hat{\beta}_n)\right|^3\right]\\
&\leq
x^32^{3(k+1)}\lambda_0(M)^3
\E\left[\1_{E_n}\,|\hat{\beta}_n-\beta_0|^3\sup_{x\in\R}\left|D^{(1)}_n(\beta^*;x)\right |^3\right]\\
&\leq
x^32^{3(k+1)}\lambda_0(M)^3
\frac{L^3\xi^{3/2}_2}{n}
\leq
\frac{Cx^32^{3(k+1)}}{n},
\end{split}
\]
for some $C>0$.

To bound the first expectation in~\eqref{eqn:1SM}, we use Lemma~\ref{le:martSM} and we can argue as in the proof of Lemma~\ref{le:inv} to obtain
\[
\p\left(\left\{U(a)+x\leq \hat{U}_n(a)\leq M \right\}\cap E_n\right)\\
\leq
\frac{K}{nx^3}.
\]
We can deal in the same way as in the proof of Lemma~\ref{le:inv} with the second probability on the right hand side of~\eqref{eqn:mainSM}, using the properties of $\bar{\w}_n$ and $\bar{R}_n$.
\end{proof}
Note that on the event $E_n$ from Lemma~\ref{le:2-1 SM},
similar to~\eqref{eqn:E_n5}, we have
\begin{equation}
\label{eqn:E_n5SM}
\sup_{x\in\R}
|\Phi_n(x;\hat{\beta}_n)-\Phi(x;\beta_0)|
\leq \frac{C_\phi}{n^{1/3}},
\end{equation}
where $C_\phi=\sqrt{\xi_2}L+\xi_4$, with $L$ the upper bound from~\eqref{eqn:dn2}.

\begin{lemma}
\label{le:lambdaSM}
Suppose that (A1)--(A2) hold.
Take $0<\epsilon<\epsilon'<M'<M<\tau_H$.
Let $\hat{\lambda}_n$ be the maximum likelihood estimator of a nondecreasing baseline hazard rate $\lambda_0$,
which is differentiable with $\lambda_0'$ uniformly bounded above and below by strictly positive constants.
Let $E_n$ be the event from Lemma~\ref{le:2-1 SM}.
Take $\xi_2>0$ and $\xi_4>0$ in~\eqref{eqn:ef} sufficiently small, such that
\begin{equation}
\label{eqn:C SM}
C_\phi<\frac{\Phi(M;\beta_0)}{2\lambda_0(M)}
\min\left\{\epsilon'-\epsilon, M-M'\right\}
\inf_{x\in[0,\tau_H]}\lambda'_0(x).
\end{equation}
and take $\xi_3$ in~\eqref{eq:def En3} sufficiently small, such that
\begin{equation}
\label{eqn:xiSM}
\xi_3\leq \frac{1}{4}
\left\{
\frac{(M-M')}{2}\inf_{x\in[0,\tau_H]}\lambda'_0(x)-\frac{C_\phi}{\Phi(M;\beta_0)}\lambda_0(M)
\right\}
(M-M')\Phi(M;\beta_0).
\end{equation}
Then, there exists a constant $C$  such that, for each $n\in\N$,
\[
\sup_{t\in[\epsilon',M']}
\E\left[
n^{2/3}\1_{E_n}\big(\lambda_0(t)-\hat{\lambda}_n(t)\big)^2
\right]\leq C.
\]
\end{lemma}
\begin{proof}
It is sufficient to prove that there exist some constants $C_1,$ $C_2$  such that for each $n\in\N$ and each $t\in[\epsilon',M']$, we have
\begin{gather}
\label{eqn:exp1SM}
\E\left[n^{2/3}
\1_{E_n}
\left\{(\hat{\lambda}_n(t)-\lambda_0(t))_+\right\}^2\right]\leq C_1,\\
\label{eqn:exp2SM}
\E\left[n^{2/3}\1_{E_n}
\left\{(\lambda_0(t)-\hat{\lambda}_n(t))_+\right\}^2\right]\leq C_2.
\end{gather}
Lets first consider~\eqref{eqn:exp1SM}.
Then completely similar to the proof of Lemma~\ref{le:lambda} we have
\[
\E\left[n^{2/3}
\1_{E_n}
\left\{(\hat{\lambda}_n(t)-\lambda_0(t))_+\right\}^2\right]
\leq
4\eta^2+4\int_{\eta}^\infty
\p
\left(
n^{1/3}\1_{E_n}(\hat{\lambda}_n(t)-\lambda_0(t))> x
\right)x\,\mathrm{d}x,
\]
for a fixed $\eta>0$,
where
\[	
\p
\left(
n^{1/3}\,\1_{E_n}\,(\hat{\lambda}_n(t)-\lambda_0(t))> x
\right)
=
\p\left(\{\hat{U}_n(a+n^{-1/3}x)< t\}\cap E_n\right).
\]
We distinguish between the cases
\[
a+n^{-1/3}x\leq \lambda_0(M)\quad\text{ and }\quad a+n^{-1/3}x> \lambda_0(M),
\]
where $a=\lambda_0(t)$.	
We prove that, in the first case, there exist a positive constant $C$ such that for all $t\in(\epsilon,M']$, and $n\in\N$,
\[
\p\{n^{1/3}\1_{E_n}(\hat{\lambda}_n(t)-\lambda_0(t))> x\}\leq \frac{C}{x^3},
\]
for all $x\geq \eta$, and in the second case $\p\left(n^{1/3}\1_{E_n}(\hat{\lambda}_n(t)-\lambda_0(t))> x\right)=0$.
Then~\eqref{eqn:exp1SM} follows immediately. 	

First assume $a+n^{-1/3}x>\lambda_0(M)$.
Note that, if $\hat{\lambda}_n(t)>a+n^{-1/3}x $, then for each $y>t$, we have
\[	
V_n(y)-V_n(t)
\geq
\hat{\lambda}_n(t)
\left(\hat{W}_n(y)-\hat{W}_n(t)\right)
>
\left(a+n^{-1/3}x\right)
\left(\hat{W}_n(y)-\hat{W}_n(t)\right).
\]
In particular for $y=\tilde{M}=M'+(M-M')/2$, we obtain
\begin{equation}
\label{eq:tail1}
\begin{split}
&
\p
\left(n^{1/3}\1_{E_n}(\hat{\lambda}_n(t)-\lambda_0(t))> x\right)\\
&\leq
\p
\Big(\Big\{V_n(\tilde{M})-V_n(t)-\big(H^{uc}(\tilde{M})-H^{uc}(t)\big)\\
&\qquad\qquad>
\left(a+n^{-1/3}x\right)\left(\hat{W}_n(\tilde{M})-\hat{W}_n(t)\right)
-
\left(H^{uc}(\tilde{M})-H^{uc}(t)\right)\Big\}\cap E_n\Big)\\
&\leq
\p\Bigg(\Bigg\{
2\sup_{x\in[T_{(1)},T_{(n)}]}|V_n(x)-H^{uc}(x)|\\
&\qquad\qquad>
\int_t^{\tilde{M}}
\left\{
\left(a+n^{-1/3}x\right)\Phi_n(s;\hat{\beta}_n)-\lambda_0(s)\Phi(s;\beta_0)
\right\}\,\mathrm{d}s
\Bigg\}\cap E_n\Bigg).
\end{split}
\end{equation}
Note that according to~\eqref{eqn:E_n5SM}, $\Phi_n(s;\hat{\beta}_n)-\Phi(s,\beta_0)\geq -C_\phi$,
and that $a+n^{-1/3}x> \lambda_0(M)>\lambda_0(\tilde M)\geq \lambda_0(s)$.
Therefore, since $C_\phi\leq \Phi(M;\beta_0)$,
from~\eqref{eqn:xiSM}, we have
\begin{equation}
\label{eq:prob zero3}
\begin{split}
&
\int_t^{\tilde{M}}
\left\{
\left(a+n^{-1/3}x\right)\Phi_n(s;\hat{\beta}_n)-\lambda_0(s)\Phi(s;\beta_0)
\right\}\,\mathrm{d}s\\
&\geq
\Phi(M;\beta_0)\int_t^{\tilde{M}}
\left\{
\left(a+n^{-1/3}x\right)\left(1-\frac{C_\phi}{\Phi(M;\beta_0)}\right)-\lambda_0(s)
\right\}\,\mathrm{d}s\\
&>
\Phi(M;\beta_0)(\tilde{M}-M')\left(\lambda_0(M)\left(1-\frac{C_\phi}{\Phi(M;\beta_0)}\right)-\lambda_0(\tilde{M})\right)\\
&\geq
\frac{1}{2}
\left\{
\frac{M-M'}{2}\inf_{x\in[0,\tau_H]}\lambda'_0(x)-\frac{C_\phi}{\Phi(M;\beta_0)}\lambda_0(M)
\right\}
(M-M')\Phi(M;\beta_0)\geq2\xi_3.
\end{split}
\end{equation}
Hence, similar to~\eqref{eq:prob zero} and~\eqref{eq:prob zero2},
we conclude that the probability
on the right hand side of~\eqref{eq:tail1} is zero.

Then, consider the case $a+n^{-1/3}x\leq \lambda_0(M)$.
Similar to~\eqref{eq:bound case 1}, from Lemma~\ref{le:invSM}, we have
\[
\p\left(
\left\{\epsilon\leq \hat{U}_n(a+n^{-1/3}x)< t\right\}\cap E_n
\right)
\leq
\frac{C}{x^3}
\]
for some $C>0$.
Moreover, for the case $\hat{U}_n(a+n^{-1/3}x)<\epsilon$, we find
\[
\p\left(\left\{\hat{U}_n(a+n^{-1/3}x)<\epsilon\right\}\cap E_n\right)
=
\p\left(\left\{\hat{\lambda}_n(\epsilon)> a+n^{-1/3}x\right\}\cap E_n\right).
\]
Note that, if $\hat{\lambda}_n(\epsilon)>a+n^{-1/3}x $, then for each $y>\epsilon$, we have
\[	
V_n(y)-V_n(\epsilon)
\geq
\hat{\lambda}_n(\epsilon)
\left(\hat{W}_n(y)-\hat{W}_n(\epsilon)\right)
>
\left(a+n^{-1/3}x\right)\left(\hat{W}_n(y)-\hat{W}_n(\epsilon)\right).	
\]
In particular for $y=\tilde{\epsilon}=\epsilon+(\epsilon'-\epsilon)/2$,
similar to~\eqref{eq:tail1},
we obtain
\begin{equation}
\label{eq:tail2}
\begin{split}
&
\p\left(\left\{\hat{\lambda}_n(\epsilon)> a+n^{-1/3}x\right\}\cap E_n\right)\\
&\leq
\p\left(\left\{V_n(\tilde{\epsilon})-V_n(\epsilon)
>
\left(a+n^{-1/3}x\right)\left(\hat{W}_n(y)-\hat{W}_n(\epsilon)\right)\right\}\cap E_n\right)\\
&\leq
\p\Bigg(\Bigg\{
2\sup_{x\in[T_{(1)},T_{(n)}]}|V_n(x)-H^{uc}(x)|\\
&\qquad\qquad>
\int_\epsilon^{\tilde{\epsilon}}
\left\{
\left(a+n^{-1/3}x\right)\Phi_n(s;\hat{\beta}_n)-\lambda_0(s)\Phi(s;\beta_0)
\right\}\,\mathrm{d}s
\Bigg\}\cap E_n\Bigg).
\end{split}
\end{equation}
Then, similar to~\eqref{eq:prob zero3} and using $\tilde\epsilon>\epsilon'$,
from~\eqref{eqn:xiSM} we obtain
\[
\begin{split}
&
\int_\epsilon^{\tilde{\epsilon}}
\left\{
\left(a+n^{-1/3}x\right)\Phi_n(s;\hat{\beta}_n)-\lambda_0(s)\Phi(s;\beta_0)
\right\}\,\mathrm{d}s\\
&\geq
\frac{1}{2}
\left\{
\frac{\epsilon'-\epsilon}{2}\inf_{x\in[0,\tau_H]}\lambda'_0(x)-\frac{C_\phi}{\Phi(M,\beta_0)}\lambda_0(\epsilon')
\right\}
\left(\epsilon'-\epsilon\right)\Phi(M;\beta_0)
\geq
2\xi_3,
\end{split}
\]
and we conclude that the probability on the right hand side of~\eqref{eq:tail2} is zero.
This concludes the proof of~\eqref{eqn:exp1SM}.

We proceed with~\eqref{eqn:exp2SM}.
Arguing as in the proof of~\eqref{eqn:exp1SM}, we obtain
\[
\E\left[n^{2/3}\1_{E_n}
\left\{(\lambda_0(t)-\hat{\lambda}_n(t))_+\right\}^2\right]
\leq
\eta^2+2\int_{\eta}^\infty
\p\left(n^{1/3}\1_{E_n}\left(\lambda_0(t)-\hat{\lambda}_n(t)\right)\geq x\right)x\,\mathrm{d}x,
\]
for a fixed $\eta>0$, where
\[
\p\left(n^{1/3}\1_{E_n}\left(\lambda_0(t)-\hat{\lambda}_n(t)\right)\geq x\right)
=
\p\left(
\left\{\hat{U}_n(a-n^{-1/3}x)\geq t\right\}\cap E_n\right),
\]
where $a=\lambda_0(t)$.
First consider the case $0<a-n^{-1/3}x\leq\lambda_0(\epsilon)$.
For each $y<t$, we have
\[
V_n(t)-V_n(y)
\leq
\hat{\lambda}_n(t)
\left(\hat{W}_n(t)-\hat{W}_n(y)\right)
\leq
\left(a-n^{-1/3}x\right)
\left(\hat{W}_n(t)-\hat{W}_n(y)\right).	
\]
In particular, for $y=\tilde{\epsilon}=\epsilon+(\epsilon'-\epsilon)/2$,
similar to~\eqref{eq:tail2}, we obtain
\begin{equation}
\label{eq:tail3}
\begin{split}
&
\p\left(n^{1/3}\1_{E_n}\left(\lambda_0(t)-\hat{\lambda}_n(t)\right)\geq x\right)\\
&\leq
\p\Bigg(\Bigg\{
2\sup_{x\in[T_{(1)},T_{(n)}]}|V_n(x)-H^{uc}(x)|\\
&\qquad\qquad\geq
\int_{\tilde{\epsilon}}^t
\left\{\left(-a+n^{-1/3}x\right)\Phi_n(s;\hat{\beta}_n)+\lambda_0(s)\Phi(s;\beta_0)\right\}\,\mathrm{d}s
\Bigg\}\cap E_n\Bigg).
\end{split}
\end{equation}
As before, using that $-a+n^{-1/3}x+\lambda_0(s)>0$
and $t-\tilde\epsilon\geq \epsilon'-\tilde\epsilon=\frac12(\epsilon'-\epsilon)$,
similar to~\eqref{eq:prob zero3}, from~\eqref{eqn:xiSM} we have
\begin{equation}
\label{eq:prob zero4}
\begin{split}
&
\int_{\tilde{\epsilon}}^t
\left\{\left(-a+n^{-1/3}x\right)\Phi_n(s;\hat{\beta}_n)+\lambda_0(s)\Phi(s;\beta_0)\right\}\,\mathrm{d}s\\
&\geq
\frac{1}{2}
\left\{\frac{\epsilon'-\epsilon}{2}\inf_{x\in[0,\tau_H]}\lambda'_0(x)
-
\frac{C_\phi}{\Phi(M;\beta_0)}\lambda_0(\epsilon)\right\}(\epsilon'-\epsilon)\Phi(M;\beta_0)
\geq
2\xi_3.
\end{split}
\end{equation}
and we conclude that the probability on the right hand side of~\eqref{eq:tail3} is zero.

Next, suppose that $a-n^{-1/3}x>\lambda_0(\epsilon)$ and consider
\[
\p\left(\left\{\hat{U}_n(a-n^{-1/3}x)\geq t\right\}\cap E_n\right).
\]
In order to use Lemma~\ref{le:invSM}, we must intersect with the event $\{\epsilon\leq \hat U_n(a-n^{-1/3}x)\leq M\}$.
Note that since $t\in[\epsilon',M']$, we have that $\hat{U}_n(a-n^{-1/3}x)\geq t$ implies $\hat{U}_n(a-n^{-1/3}x)\geq \epsilon$.
Using Lemma~\ref{le:invSM} and the mean value theorem, we obtain
\[
\begin{split}
&
\p\left(\left\{t\leq \hat{U}_n(a-n^{-1/3}x)\leq M\right\}\cap E_n\right)\\
&=
\p\Bigg(
\left\{
\hat{U}_n(a-n^{-1/3}x)-U(a-n^{-1/3}x)\geq t-U(a-n^{-1/3}x)
\right\}\\
&\qquad\qquad\qquad\qquad\qquad\qquad\cap
\left\{\epsilon\leq\hat{U}_n(a-n^{-1/3}x)\leq M\right\}\cap E_n
\Bigg)\\
&\leq
\p\Bigg(
\left\{
|\hat{U}_n(a-n^{-1/3}x)-U(a-n^{-1/3}x)|\geq t-U(a-n^{-1/3}x)
\right\}\\
&\qquad\qquad\qquad\qquad\qquad\qquad\cap
\left\{\epsilon\leq\hat{U}_n(a-n^{-1/3}x)\leq M \right\}\cap E_n
\Bigg)\\
&\leq
\frac{K}{n\left\{t-U(a-n^{-1/3}x)\right\}^3}
\leq
\frac{K}{(U'(\xi_n))^3x^3}
\leq \frac{C}{x^3}
\end{split}
\]
where we use that $t=U(a)$, $\xi_n\in(a-n^{-1/3}x,a)$, and $U'(\xi_n)=1/\lambda'_0(\lambda_0^{-1}(\xi_n))$ is bounded.
Finally, note that
\[
\p\left(\left\{\hat{U}_n(a-n^{-1/3}x)> M\right\}\cap E_n\right)
\leq
\p\left(\left\{\hat{\lambda}_n(M)\leq a-n^{-1/3}x\right\}\cap E_n\right).
\]
If $\hat{\lambda}_n(M)\leq a-n^{-1/3}x $, then for each $y<M$, we have
\[
V_n(M)-V_n(y)
\leq
\hat{\lambda}_n(M)\left(\hat{W}_n(M)-\hat{W}_n(y)\right)
\leq
\left(a-n^{-1/3}x\right)
\left(\hat{W}_n(M)-\hat{W}_n(y)\right).	
\]
In particular for $y=\tilde{M}=M'+(M-M')/2$,
similar to~\eqref{eq:tail3}, we obtain
\begin{equation}
\label{eq:tail4}
\begin{split}
&
\p\left(\left\{\hat{\lambda}_n(M)\leq a-n^{-1/3}x\right\}\cap E_n\right)\\
&\leq
\p\Bigg(\Bigg\{
2\sup_{x\in[T_{(1)},T_{(n)}]}|V_n(x)-H^{uc}(x)|\\
&\qquad\qquad\geq
\int_{\tilde{M}}^M
\left\{
\left(a-n^{-1/3}x\right)\Phi_n(s;\hat{\beta}_n)-\lambda_0(s)\Phi(s;\beta_0)
\right\}\,\mathrm{d}s
\Bigg\}\cap E_n\Bigg)
\end{split}
\end{equation}
As before, similar to~\eqref{eq:prob zero4}, from~\eqref{eqn:xiSM} we have
\[
\begin{split}
&
\int_{\tilde{M}}^M
\left\{
\left(a-n^{-1/3}x\right)\Phi_n(s;\hat{\beta}_n)-\lambda_0(s)\Phi(s;\beta_0)
\right\}\,\mathrm{d}s\\
&\geq
\frac{1}{2}
\left\{
\frac{M-M'}2\inf_{x\in[0,\tau_H]}\lambda'_0(x)-\frac{C_\phi}{\Phi(M;\beta_0)}\lambda_0(M')
\right\}
(M-M')\Phi(M;\beta_0)
\geq2\xi_3.
\end{split}
\]
and we conclude that the probability on the right hand side of~\eqref{eq:tail4} is zero.
This finishes the proof.
\end{proof}

\begin{lemma}
\label{le:lambda2SM}
Suppose that (A1)--(A2) hold.
Fix $x\in(0,\tau_h)$ and take $0<\epsilon<\epsilon'<M'<M<\tau_H$.
Let $\hat{\lambda}_n$ be the maximum likelihood estimator of a nondecreasing baseline hazard rate~$\lambda_0$ which is differentiable with $\lambda_0'$
uniformly bounded above and below by strictly positive constants.
Let $E_n$ be the event from Lemma~\ref{le:2-1 SM} and choose $C_\phi$ and $\xi_3$
such that they satisfy~\eqref{eqn:C SM} and~\eqref{eqn:xiSM}, respectively.
Then
	\[
	\1_{E_n}\,\int_{x-b}^{x+b}(\lambda_0(t)-\hat{\lambda}_n(t))^2\,\mathrm{d}t=O_p(bn^{-2/3}).
\]
\end{lemma}
\begin{proof}
The proof is exactly the same as the proof of Lemma~\ref{le:lambda2} (replacing $\tilde{\lambda}_n$ with $\hat{\lambda}_n$).
\end{proof}
We are now in the position to establish the analogue~\eqref{eq:Lemma33 SM} of Lemma~\ref{le:result4}.

\begin{lemma}
\label{le:result4 SM}
Suppose that (A1)--(A2) hold.
Fix $x\in(0,\tau_h)$ and let $\theta_{n,x}$ and $\bar\theta_{n,x}$ be defined by~\eqref{eqn:thetan SM} and~\eqref{eq:def theta bar SM}, respectively.
Assume that $\lambda_0$ is differentiable, such that $\lambda'_0$ is uniformly bounded above and below by strictly positive constants
and let $k$ satisfy~\eqref{def:kernel}.
Then, it holds
\[
\int
\left\{
\bar{\theta}_{n,x}(u,\delta,z)-\theta_{n,x}(u,\delta,z)
\right\}
\,\mathrm{d}\p(u,\delta,z)=
O_p(b^{-1}n^{-2/3}).
\]
 \end{lemma}
\begin{proof}
Take $0<\epsilon<\epsilon'<x<M'<M<\tau_H$ and consider $n$ sufficiently large such that $[x-b,x+b]\subset[\epsilon',M']$.
Similar to~\eqref{eq:decomp thetabar}, we have
\[
\begin{split}
&
\int
\left\{
\bar{\theta}_{n,x}(u,\delta,z)-\theta_{n,x}(u,\delta,z)
\right\}
\,\mathrm{d}\p(u,\delta)\\
&=
\1_{E_n}\int_{x-b}^{x+b}
\left(
a_{n,x}(\hat{A}_n(u))-a_{n,x}(u)
\right)
\left(\Phi(u;\beta_0)\lambda_0(u)-\Phi(u;\hat{\beta}_n)\hat{\lambda}_n(u)\right)\,\mathrm{d}u
\end{split}
\]
so that by Cauchy-Schwarz inequality
\begin{equation}
\label{eqn:sch SM}
\begin{split}
&
\left|\int
\left\{
\bar{\theta}_{n,x}(u,\delta,z)-\theta_{n,x}(u,\delta,z)
\right\}\,\mathrm{d}\p(u,\delta,z) \right| \\
&\leq
\1_{E_n}
\left\Vert
\left(
a_{n,x}\circ\hat{A}_n-a_{n,x}
\right)\1_{[x-b,x+b]}
\right\Vert_{\mathcal{L}_2}
\left\Vert
\left(
\Phi_0\lambda_0-\hat\Phi_n\hat{\lambda}_n
\right)
\1_{[x-b,x+b]}
\right\Vert_{\mathcal{L}_2},
\end{split}
\end{equation}
where $\Phi_0(u)=\Phi(u;\beta_0)$ and $\hat\Phi_n(u)=\Phi_n(u;\hat\beta_n)$.
Similar to~\eqref{eqn:l2},
\[
\1_{E_n}\left\Vert (a_{n,x}(\hat{A}_n)-a_{n,x})\1_{[x-b,x+b]}\right\Vert_{\mathcal{L}_2}^2
\leq
\frac{c}{b^4}\1_{E_n}\int_{x-b}^{x+b}(\hat{A}_n(u)-u)^2\,\mathrm{d}u,
\]
for some constant $c$,
and by the same reasoning as in~\eqref{eqn:lambda},
for $u\in[\tau_i,\tau_{i+1})$ and $\hat{A}_n(u)<\tau_{i+1}$, we obtain
\[
|u-\hat{A}_n(u)|
\leq
2K|\lambda_0(u)-\hat{\lambda}_n(u)|,
\]
which also holds in the case $\hat{A}_n(u)=\tau_{i+1}$ simply because $|\lambda_0(u)-\lambda_0(\hat{A}_n(u))|\leq |\lambda_0(u)-\hat{\lambda}_n(u)|$.
As a result, using Lemma~\ref{le:lambda2SM}, we derive that
\[
\1_{E_n}\frac{1}{b^4}
\int_{x-b}^{x+b}
\left(\hat{A}_n(u)-u\right)^2\,\mathrm{d}u
\leq
\frac{C}{b^4}\1_{E_n}
\int_{x-b}^{x+b}
\left(\lambda_0(u)-\hat{\lambda}_n(u)\right)^2\,\mathrm{d}u
=
O_p(b^{-3}n^{-2/3}).
\]
The argument for second factor in~\eqref{eqn:sch SM} is the same as for~\eqref{eq:bound Phi0lambda0},
and yields
\[
\1_{E_n}
\left\Vert
\left(\Phi_0\lambda_0-\widehat{\Phi}\tilde{\lambda}_n\right)
\1_{[x-b,x+b]}
\right\Vert_{\mathcal{L}_2}
=
O_p(b^{1/2}n^{-1/3}).
\]
Together with~\eqref{eqn:sch SM}, this concludes the proof.
\end{proof}

To establish the analogue of Lemma~\ref{le:result5} for $\hat\lambda_n$,
similar to Lemma~\ref{le:sup}, we need a stronger version of Lemma~\ref{le:lambdaSM}.
As before, we loose a factor $n^{2/9}$ with respect to the bound in Lemma~\ref{le:lambdaSM},
which might not be optimal, but suffices for our purposes.
\begin{lemma}
\label{le:sup SM}
Suppose that (A1)--(A2) hold.
Take $0<\epsilon<\epsilon'<M'<M<\tau_H$.
Let~$\hat{\lambda}_n$ be the maximum likelihood estimator of a nondecreasing baseline hazard rate $\lambda_0$,
which is differentiable with $\lambda_0'$ uniformly bounded above and below by strictly positive constants.
Let $E_n$ be the event from Lemma~\ref{le:2-1 SM} and choose $C_\phi$ and $\xi_3$
such that they satisfy~\eqref{eqn:C SM} and~\eqref{eqn:xiSM}, respectively.
Then, there exists a constant $C>0$  such that, for each $n\in\N$,
\begin{equation}
\label{eqn:lemmaSM}
\E
\left[
n^{4/9}\1_{E_n}\,\sup_{t\in[\epsilon',M'] }
\left(\lambda_0(t)-\hat{\lambda}_n(t)\right)^2
\right]
\leq C.
\end{equation}
\end{lemma}
\begin{proof}
The proof is exactly the same as the proof of Lemma~\ref{le:sup} (replacing $\tilde{\lambda}_n$ with $\hat{\lambda}_n$).
\end{proof}

\begin{lemma}
\label{le:result5 SM}
Suppose that (A1)--(A2) hold.
Fix $x\in(0,\tau_h)$ and take $0<\epsilon<\epsilon'<x<M'<M<\tau_H$.
Assume that $\lambda_0$ is differentiable, and such that $\lambda'_0$ is uniformly
bounded above and below by strictly positive constants.
Assume that $x\mapsto \Phi(x;\beta_0)$ is differentiable with a bounded derivative in a neighborhood of $x$.
Let $\bar\theta_{n,x}$ be defined in~\eqref{eq:def theta bar SM} and let $\eta_{n,x}$ be defined by~\eqref{eqn:eta},
where $E_n$ is the event from Lemma~\ref{le:2-1 SM}.
Let $k$ satisfy~\eqref{def:kernel}.
Then, it holds
\begin{equation}
\label{eqn:statement SM}
\int
\left\{
\bar{\theta}_{n,x}(u,\delta,z)-\eta_{n,x}(u,\delta,z)
\right\}
\,\mathrm{d}(\p_n-\p)(u,\delta,z)
=
O_p(b^{-3/2}n^{-13/18})+O_p(n^{-1/2}).
\end{equation}
\end{lemma}
\begin{proof}
Let $n$ be sufficiently large, such that $\epsilon'<x-b<x+b<M'$.
Denote by $R_n$ the left hand side of~\eqref{eqn:statement SM} and write $R_n=R_{n1}+R_{n2}$,
where
\[
\begin{split}
R_{n1}
&=
n^{-1/2}
\1_{E_n}\1_{[x-b,x+b]}(u)
\int\delta
\left\{
\bar{a}_{n,x}(u)-a_{n,x}(u)
\right\}\,\mathrm{d}\sqrt{n}(\p_n-\p)(u,\delta,z),\\
R_{n2}
&=
n^{-1/2}
\1_{E_n}
\int \1_{\{u>x-b\}}
\Bigg\{
\mathrm{e}^{\hat{\beta}'_nz}\int_{x-b}^{u\wedge (x+b)}\bar{a}_{n,x}n(u)\,\mathrm{d}\hat{\Lambda}_n(v)\\
&\qquad\qquad\qquad\qquad\qquad\qquad\qquad-
\mathrm{e}^{\beta'_0z}\int_{x-b}^{u\wedge (x+b)} a_{n,x}(v)\,\mathrm{d}\Lambda_0(v)\Bigg\}\,\mathrm{d}\sqrt{n}(\p_n-\p)(u,\delta,z).
\end{split}
\]
Choose $\eta>0$.
Consider similar two events as in~\eqref{def:events}:
\begin{equation}
\label{def:eventsSM}
\begin{split}
\A_{n1}&=\left\{\hat{\lambda}_n(M)>K_1\right\},\\
\A_{n2}&=\left\{\sup_{t\in[\epsilon',M']}\left|\lambda_0(t)-\hat{\lambda}_n(t)\right|>K_2n^{-2/9}\right\},
\end{split}
\end{equation}
where $K_1,K_2>0$, and let $\A_n=\A_{n1}\cup\A_{n2}$.
From Lemma~\ref{le:sup SM} and the fact that  $\hat{\lambda}_n(M)=O_p(1)$, it follows
that we can choose $K_1,K_2>0$ such that $\p(\A_n)\leq 2\eta/3$.
As in the proof of Lemma~\ref{le:result5}, it suffices
to show that there exists $\nu>0$, such that
$b^{3/2}n^{13/18}\nu^{-1}\E\big[ |R_{n1}|\1_{\A_n^c}\big]\leq\eta/3$
and
$n^{1/2}\nu^{-1}\E\big[ |R_{n2}|\1_{\A_n^c}\big]\leq\eta/3$,
for all $n$ sufficiently large.

Let us first consider $R_{n1}$.
We have
\begin{equation}
\label{eq:bound Rn1 SM}
\begin{split}
a_{n,x}(\hat{A}_n(u))-a_{n,x}(u)
&=
\frac{k_b(x-\hat{A}_n(u))-k_b(x-u)}{\Phi(\hat{A}_n(u);\beta_0)}\\
&\qquad
+k_b(x-u)\frac{\Phi(u;\beta_0)-\Phi(\hat{A}_n(u);\beta_0)}{\Phi(\hat{A}_n(u);\beta_0)\Phi(u;\beta_0)}.
\end{split}
\end{equation}
Similar to~\eqref{eq:bound Rn1},
\[
\left|
k_b(x-\hat{A}_n(u))-k_b(x-u)
\right|
\leq
b^{-2}n^{-2/9}
K_2
\sup_{x\in[-1,1]}|k'(x)|,
\]
for some $K_2>0$,
and similarly, using that
$x\mapsto\Phi(x;\beta_0)$ is differentiable with bounded derivative in a neighborhood of $x$,
\[
\begin{split}
b^{-1}|\Phi(u;\beta_0)-\Phi(\hat{A}_n(u);\beta_0)|
\leq
Kb^{-1}\,|\hat{A}_n(u)-u|
\leq
b^{-1}n^{-2/9}
KK_3,
\end{split}
\]
Consequently, $R_{n1}$ can be written as
\[
R_{n1}
=
\1_{E_n}
b^{-2}n^{-13/18}
\int \1_{[x-b,x+b]}(u)
\delta
W_n(u)\,d\sqrt{n}(\p_n-\p)(u,\delta,z),
\]
where $W_n$ is a function of bounded variation, uniformly bounded.
Completely similar to the proof of Lemma~\ref{le:result5},
together with Lemma~\ref{lem:covering number} we find that
\[
b^{2}n^{13/18}
\nu^{-1}
\E
\left[
|R_{n1}|\1_{\A_n^c}
\right]
\leq
\frac{K''}{\nu n^{2/90}}\leq
\eta/3,
\]
for sufficiently large $\nu$.
For $R_{n2}$ we write
\[
n^{1/2}R_{n2}
=
\1_{E_n}
\int
\left(
\mathrm{e}^{\hat{\beta}'_nz}r_{1,n}(u)-\mathrm{e}^{\beta'_0z}r_{2,n}(u)
\right)
\,\mathrm{d}\sqrt{n}(\p_n-\p)(u,\delta,z),
\]
where
\[
\begin{split}
r_{1,n}(u)
&=
\1_{\{u>x-b\}}
\int_{x-b}^{u\wedge(x+b)}a_{n,x}(\hat{A}_n(v))\,\hat{\lambda}_n(v)\,\mathrm{d}v,\\
r_{2,n}(u)
&=
\1_{\{u>x-b\}}
\int_{x-b}^{u\wedge(x+b)}a_{n,x}(v)\,\lambda_0(v)\,\mathrm{d}v,
\end{split}
\]
are both monotone functions, uniformly bounded by some constant $C$ on the event $\A_n^c$.
Once more, from here we follow exactly the same proof as the one for Lemma~\ref{lem:covering number}.
\end{proof}

\begin{theo}
\label{theo:distrSM}
Suppose that (A1)--(A2) hold.
Fix $x\in(0,\tau_h)$.
Assume that $\lambda_0$ is $m\geq2$ times differentiable in $x$,
such that $\lambda'_0$ is uniformly bounded above and below by strictly positive constants.
Moreover, that $t\mapsto \Phi(t;\beta_0)$ is differentiable with a bounded derivative in a neighborhood of $x$,
and let $k$ satisfy~\eqref{def:kernel}.
Let~$\hat{\lambda}_n^{SM}$ be defined in~\eqref{def:lambdaSM} and assume that $n^{1/(2m+1)}b\to c>0$.
Then, it holds
\[
n^{m/(2m+1)}
\left(\hat{\lambda}_n^{SM}(x)-\lambda_0(x)\right)\xrightarrow{d}N(\mu,\sigma^2),
\]
where
\[
\mu=
\frac{(-c)^m}{m!}\lambda^{(m)}_0(x)
\int
k(u)u^m\,\mathrm{d}u
\quad\text{ and }\quad
\sigma^2
=
\frac{\lambda_0(x)}{c\Phi(x;\beta_0)}
\int k^2(u)\,\mathrm{d}u.
\]
\end{theo}
\begin{proof}
The proof is completely analogous to that of Theorem~\ref{theo:distrSM} and is based on a similar decomposition as in~\eqref{eq:decomposition lambdaSG}.
After using Lemmas~\ref{le:2-1 SM}, \ref{le:2-2 SM}, \ref{le:result4 SM}, and~\ref{le:result5 SM}, it remains to obtain the limit of~\eqref{eq:decomp eta},
where $E_n$ is the event from Lemma~\ref{le:2-1 SM}.
This is completely similar to the argument in the proof of Theorem~\ref{theo:distr}.
\end{proof}

\section{Consistency of the bootstrap}
\label{sec:bootstrap}
Instead of $\p_n$, we consider $\p_n^*$, the empirical measure corresponding to the bootstrap sample 
	$(T_1^*,\Delta_1^*,Z_1)$,$\ldots$, $(T_n^*,\Delta_n^*,Z_n)$, and instead of $\p$, we consider $P_n^*$, the measure corresponding to the bootstrap distribution of
	$(T^*,\Delta^*,Z)=(\min(X^*,C^*),\1_{\{X^*\leq C^*\}},Z)$ conditional on the data $(T_1,\Delta_1,Z_1)$, $\ldots$ ,$(T_n,\Delta_n,Z_n)$, where
	$X^*$ conditional on $Z$ has distribution function $\hat F_n(x\mid Z)$ and $C^*$ has distribution function $\hat{G}_n$.
	To prove~\eqref{eqn:asymptotic_normality_bootstrap}, we mimic the proof of Theorem~\ref{theo:distr}, 
	which means that one needs to establish the bootstrap versions of Lemmas~\ref{le:2-1}-\ref{le:result5}.
	
	In view of Remark~\ref{rem:partial MLE}, let $\hat\beta_n$ be an estimate for $\beta_0$ satisfying~\eqref{eq:cond PMLE}. 
	Let $\hat{\beta}^*_n$ be the bootstrap version and suppose that 
	$\hat\beta_n^*-\hat{\beta}_n\to 0$,
	for almost all sequences $(T_i^*,\Delta_i^*,Z_i)$, $i=1,2,\ldots$, conditional on the sequence  $(T_i,\Delta_i,Z_i)$, $i=1,2,\ldots$,
	and that $\sqrt{n}(\hat\beta_n^*-\hat{\beta}_n)=O_p^*(1)$,
	meaning that for all $\epsilon>0$, there exists $M>0$ such that
	\[
	\limsup_{n\to\infty}
	P_n^*
	\left(
	\sqrt{n}
	|\hat\beta_n^*-\hat{\beta}_n|
	>M
	\right)
	<
	\epsilon,
	\qquad
	\p-\text{almost surely}.
	\]
	Then, similar to~\eqref{eqn:an} and\eqref{eq:def Phi}, define
	\[
	a_{n,x}^*(u)=\frac{k_b(x-u)}{\Phi^*(u;\hat\beta_n)}
	\quad\text{and}\quad
	\Phi^*(u;\hat\beta_n)
	=
	\int \1_{ \{t\geq u\}}\,\mathrm{e}^{\hat\beta_n'z}\,\mathrm{d}P_n^*(t,\delta,z).
	\]
	and let
	\[
	\theta_{n,x}^*(u,\delta,z)
	=
	\1_{E_n^*}
	\left\{
	\delta\,a_{n,x}^*(u)-\mathrm{e}^{(\hat\beta_n^*)'\,z}\int_0^u a_{n,x}^*(v)\,\mathrm{d}\tilde{\Lambda}^*_n(v)
	\right\}.
	\]
	Here $\tilde{\Lambda}_n^*$ is the greatest convex minorant of the bootstrap Breslow estimator
	\[
	\Lambda_n^*(x)
	=
	\int \frac{\delta\1_{\{ t\leq x\}}}{\Phi_n^*(t;\hat{\beta}^*_n)}\,\mathrm{d}\p^*_n(t,\delta,z),
	\]
	with
	\[
	\Phi_n^*(x;\beta)=\int \1_{\{t\geq x\}} \mathrm{e}^{\beta'z}\,\mathrm{d}\p_n^*(t,\delta,z),
	\]
	and $E_n^*$ is an event such that $\1_{E_n^*}=1+o_p^*(1)$,
	meaning that for all $\epsilon>0$, 
	\[
	\limsup_{n\to\infty}
	P_n^*(|\1_{E_n^*}-1|>\epsilon)=0,
	\qquad
	\p-\text{almost surely}.
	\]
	To obtain the bootstrap equivalent of Lemma~\ref{le:2-1}, we first show that
	\[
	\Lambda_0^*(x):=\int \frac{\delta\1_{\{t\leq x\}}}{\Phi^*(t;\hat\beta_n)}\,\mathrm{d}P_n^*(t,\delta,z)=\Lambda_n^s(x)
	\]
	and for constructing of the event $E_n^*$, we prove that
	$\sqrt{n}\sup_{x\in[0,M]}|\tilde\Lambda_n^*(x)-\Lambda_n^s(x)|$ and 
	$\sqrt{n}\sup_{x\in\R}|\Phi_n^*(x;\hat\beta_n)-\Phi^*(x;\hat\beta_n)|$
	are of the order $O_p^*(1)$.
	This yields the bootstrap equivalent of Lemma~\ref{le:2-1}:
	\[
	\int\theta_{n,x}^*(u,\delta,z)\,\mathrm{d}P_n^*(u,\delta,z)
	=
	-\1_{E^*_n}\int k_b(x-u)\,\mathrm{d}(\tilde{\Lambda}^*_n-\Lambda_n^s)(u)+O^*_p(n^{-1/2}).
	\]
	The proof of the bootstrap version of Lemma~\ref{le:2-2} is completely the same as that of Lemma~\ref{le:2-2}.
	The bootstrap versions of Lemmas~\ref{lem:mart2}-\ref{le:lambda2}, which are preparatory for the bootstrap version of Lemma~\ref{le:result4}, 
	can be obtained by means of similar arguments.
	This requires a suitable martingale that approximates the process $\Lambda_n^*-\Lambda_n^s$.
	To this end we define
	\[
	\mathbb{M}_n^*(t)
	=
	\sum_{i=1}^n 
	\left(
	\1_{\{X_i^*\leq t\}}\Delta_i^*-\int_0^t Y_i^*(u)\mathrm{e}^{\hat\beta'_nZ_i}\,\mathrm{d}\Lambda_n^s(u)
	\right)
	\]
	and
	\[
	\w_n^*(t)=\int_0^{t\wedge M}\frac{1}{n\Phi^*(s;\hat\beta_n)}\,\mathrm{d}\m_n^*(s),
	\]
	where the latter can be shown to be a mean zero square integrable martingale, that satisfies
	a bound similar to one in Lemma~\ref{le:mart}.
	Similar to Lemma~\ref{le:inv}, this leads to a suitable bound on the tail probabilities of the bootstrap inverse process, 
	defined for $a\in[\tilde\lambda_n^s(\epsilon),\tilde\lambda_n^s(M)]$,
	for $0<\epsilon<M<\tau_H$, by
	\begin{equation}
	\label{def:U SG*}
	U^*_n(a)
	=
	\begin{cases}
	\epsilon & a< \tilde\lambda_n^s(\epsilon);\\
	(\lambda_n^s)^{-1}(a) & a\in[\tilde\lambda_n^s(\epsilon),\tilde\lambda_n^s(M)];\\
	M & a>\tilde\lambda_n^s(M).
	\end{cases}
	\end{equation}
	This enables us to obtain $L_2$-bounds similar to Lemmas~\ref{lem:mart2} and~\ref{le:mart},
	\[
	\begin{split}
	\sup_{t\in[\epsilon',M'] }
	\E^*\left[n^{2/3}\1_{E_n^*}\left(\tilde\lambda_n^s(t)-\tilde{\lambda}_n^*(t)\right)^2\right]&\leq C;\\
	\1_{E_n^*}\int_{x-b}^{x+b}(\tilde\lambda_n^s(t)-\tilde{\lambda}_n^*(t))^2\mathrm{d}t&=O_p^*(bn^{-2/3}),
	\end{split}
	\]
	for $0<\epsilon<\epsilon'<M'<M<\tau_H$, where $\E^*$ denotes the expectation with respect to~$P_n^*$.
	Moreover, since the proof of Lemma~\ref{le:result4} makes use of the derivative of $k_b(x-y)/\Phi(y;\beta_0)$,
	differentiation of its bootstrap counterpart $k_b(x-y)/\Phi^*(y;\hat\beta_n)$ has to be circumvented.
	This is done by a suitable differentiable approximation of $\Phi^*(y;\hat\beta_n)$,
	and we then obtain the bootstrap version of Lemma~\ref{le:result4}:
	\[
	\int
	\left\{
	\overline{\theta}_{n,x}^*(u,\delta,z)-\theta_{n,x}^*(u,\delta,z)
	\right\}\,\mathrm{d}P_n^*(u,\delta,z)=
	O_p^*(b^{-1}n^{-2/3}),
	\]
	Finally, after proving the bootstrap version of Lemma~\ref{le:sup}, i.e.,
	\[
	\E^*\left[
	n^{4/9}\1_{E_n^*}\sup_{t\in(\epsilon',M'] }\left(\tilde\lambda_n^s(t)-\tilde{\lambda}_n^*(t)\right)^2
	\right]
	\leq C,
	\]
	we obtain the bootstrap version of Lemma~\ref{le:result5} for
	\[
	\eta_{n,x}^*(u,\delta,z)
	=
	\1_{E_n^*}
	\left(\delta\,a_{n,x}^*(u)-\mathrm{e}^{\hat\beta'_nz}\,\int_0^u a_{n,x}^*(v)\,\mathrm{d}\Lambda_n^s(v)\right),\quad u\in[0,\tau_H],
	\]
	by using arguments similar to those in the proof of Lemma~\ref{le:result5}.
	Next, the proof of~\eqref{eqn:asymptotic_normality_bootstrap} for~$\tilde\lambda_n^{SG,*}$ is the same as that of Theorem~\ref{theo:distr}
	for $\tilde\lambda_n^{SG}$.
	
	%

\end{document}